\documentclass[9pt,amstex]{book}
\usepackage{amssymb,amsmath} 
\usepackage{latexsym,color}
\usepackage{latexsym}

\newtheorem{dfn}{Definition}[section] 
\newtheorem{rmk}{Remark}[section]
\newtheorem{rmks}{Remarks}[section]
\newtheorem{thm}{Theorem}[section] 
\newtheorem{cor}{Corollary}[section]
\newtheorem{prop}{Proposition}[section] 
\newtheorem{lem}{Lemma}[section]
\newtheorem{exs}{Examples}[section]
\newtheorem{ex}{Example}[section] 
 
\newtheorem{nt}{Notation}[section] 

\def\cyclic{\mathop{\kern0.9ex{{+}
\kern-2.2ex\raise-.28ex\hbox{\Large\hbox
{$\circlearrowright$}}}}\limits}

\def\buildrel#1_#2^#3{\mathrel{\mathop{\kern 0pt#1}\limits_{#2}^{#3}}}

\newcommand{\Pf}{{\em Proof}. }
\newcommand{\EPf}
{%
\mbox{}%
\nolinebreak%
\hfill%
\rule{2mm}{2mm}%
\medbreak%
\par%
}

\newcommand{\trianglehorizontal}
{\mbox{\begin{picture}(40,40)\put(0,0){\line(1,0){20}}
\put(-1,-1.5){$\bullet$}
\put(37.5,-1.7){$\bullet$}\put(17.3,27.3){$\bullet$}
\put(40,0){\vector(-1,0){20}}\put(0,0){\vector(2,3){10}}
\put(10,15){\line(2,3){10}}\put(20,30){\vector(2,-3){10}}
\put(30,15){\line(2,-3){10}}\end{picture}}}
\newcommand{\Triangle}[3]{\mbox{\begin{picture}(60,100)
\put(10,20){$#1$}\put(60,20){$#3$}\put(35,60){$#2$}
\put(15,25){\begin{picture}(30,40)\put(0,0){\line(1,0){20}}
\put(-1,-1.5){$\bullet$}\put(37.5,-1.7){$\bullet$}
\put(17.3,27.3){$\bullet$}\put(40,0){\vector(-1,0){20}}
\put(0,0){\vector(2,3){10}}\put(10,15){\line(2,3){10}}
\put(20,30){\vector(2,-3){10}}\put(30,15){\line(2,-3){10}}
\end{picture}}\end{picture}}}

\newcommand{\Trianglefat}[3]{\mbox{\begin{picture}(60,100)
\put(10,20){$#1$}\put(60,20){$#3$}\put(35,60){$#2$}
\put(15,25){\begin{picture}(30,40)\thicklines

\put(0,0){\line(1,0){20}}
\put(-1,-1.5){$\bullet$}\put(37.5,-1.7){$\bullet$}
\put(17.3,27.3){$\bullet$}\put(40,0){\vector(-1,0){20}}
\put(0,0){\vector(2,3){10}}\put(10,15){\line(2,3){10}}
\put(20,30){\vector(2,-3){10}}\put(30,15){\line(2,-3){10}}
\end{picture}}\end{picture}}}

\newcommand{\carrehoriz}[1]{\mbox{\begin{picture}(100,100)
\put(10,10){\line(1,0){70}}\put(10,80){\line(1,0){70}}\put(10,10){\line(1,1){70}}

\put(10,80){\line(1,-1){70}}\put(30,30){\vector(1,1){2}}
\put(60,60){\vector(-1,-1){2}}\put(45,10){\vector(-1,0){2}}
\put(45,80){\vector(1,0){2}}\put(60,30){\vector(1,-1){2}}
\put(30,60){\vector(-1,1){2}}\put(8,8){$\bullet$}\put(8,78){$\bullet$}
\put(77,8){$\bullet$}\put(78,78){$\bullet$}\put(42,42){$\bullet$}\put(0,5){$d$}

\put(85,5){$c$}\put(85,80){$b$}\put(0,80){$a$}\put(50,43){#1}\end{picture}}}

\newcommand{\carrevertic}[1]{\mbox{\begin{picture}(100,100)
\put(10,10){\line(0,1){70}}\put(80,10){\line(0,1){70}}
\put(10,10){\line(1,1){70}}\put(10,80){\line(1,-1){70}}
\put(30,30){\vector(-1,-1){2}}\put(60,60){\vector(1,1){2}}
\put(10,45){\vector(0,1){2}}\put(80,45){\vector(0,-1){2}}
\put(60,30){\vector(-1,1){2}}\put(30,60){\vector(1,-1){2}}
\put(8,8){$\bullet$}\put(8,78){$\bullet$}\put(77,8){$\bullet$}
\put(78,78){$\bullet$}\put(42,42){$\bullet$}\put(0,5){$d$}\put(85,5){$c$}

\put(85,80){$b$}\put(0,80){$a$}\put(50,43){#1}\end{picture}}}
\newcommand{\carreverticfat}[1]{\mbox{\begin{picture}(100,100)
\thicklines
\put(10,10){\line(0,1){70}}\put(80,10){\line(0,1){70}}
\put(10,10){\line(1,1){70}}\put(10,80){\line(1,-1){70}}
\put(30,30){\vector(-1,-1){2}}\put(60,60){\vector(1,1){2}}
\put(10,45){\vector(0,1){2}}\put(80,45){\vector(0,-1){2}}
\put(60,30){\vector(-1,1){2}}\put(30,60){\vector(1,-1){2}}
\put(8,8){$\bullet$}\put(8,78){$\bullet$}\put(77,8){$\bullet$}
\put(78,78){$\bullet$}\put(42,42){$\bullet$}\put(0,5){$d$}\put(85,5){$c$}

\put(85,80){$b$}
\put(0,80){$a$}
\put(50,43){#1}
\end{picture}}}

\def\dref#1{Definition~\ref{#1}}

\newcommand{\Aut}{\mbox{$\mathtt{Aut}$}}
\newcommand{\aut}{\mbox{$\mathfrak{aut}$}}

\renewcommand{\mid}{\mbox{$\mathtt{mid}$}}
\newcommand{\End}{\mbox{$\mathtt{End}$}}
\newcommand{\homo}{\mbox{$\mathtt{Hom}$}}
\renewcommand{\span}{\mbox{$\mathtt{span}$}}
\newcommand{\Der}{\mbox{$\mathfrak{Der}$}}
\newcommand{\Id}{\mbox{$\mathtt{Id}$}}
\newcommand{\id}{\mbox{$\mathtt{id}$}}

\newcommand{\Ad}{\mbox{$\mathtt{Ad}$}}
\newcommand{\ad}{\mbox{$\mathtt{ad}$}}

\newcommand{\ddto}{\left.\frac{\rm d}{{\rm d}t}\right|_0}

\newcommand{\ddtd}{\frac{{\rm d}^{2}}{{\rm d}t^{2}}}

\newcommand{\vf}{\varphi}
\newcommand{\C}{\mathbb C}

\renewcommand{\L}{\mathbb L} 
\newcommand{\bC}{{\bf C}} 
\newcommand{\bO}{{\bf O}}

\newcommand{\bT}{{\bf T}} 
 
\newcommand{\bS}{{\bf S}}

\newcommand{\bff}{{\bf f}} 
\newcommand{\bg}{{\bf g}} 
\newcommand{\bA}{{\bf A}} 
\newcommand{\bs}{{\bf s}}

\newcommand{\bE}{{\bf E}}
\newcommand{\bD}{{\bf D}}
\newcommand{\bU}{{\bf U}}
\newcommand{\bSigma}{{\bf \Sigma}}
\newcommand{\bsigma}{{\bf \sigma}}

\newcommand{\bOmega}{{\bf \Omega}}
\newcommand{\bR}{{\bf R}}

\newcommand{\bm}{{\bf m}}

\newcommand{\D}{\mathbb D} 
\newcommand{\F}{{\cal F}{}}
\newcommand{\A}{\mathbb A}

\renewcommand{\S}{\mathbb S} 
\newcommand{\M}{\mathbb M}

\newcommand{\fW}{{\mathfrak W} {}}

\newcommand{\Z}{\mathbb Z} 

\newcommand{\R}{\mathbb R}

\newcommand{\N}{\mathbb N}

\newcommand{\arcsinh}{\mbox{\rm arcsinh}} 
 
\newcommand{\Tr}{\mbox{\rm Tr}} 

\newcommand{\g}{{\mathfrak{g}}{}}

\renewcommand{\k}{{\mathfrak{k}}{}} 
\newcommand{\p}{{\mathfrak{p}}{}} 
 
\newcommand{\fB}{{\mathfrak{B}}{}} 
 
\newcommand{\fF}{{\mathfrak{F}}{}} 
\newcommand{\fL}{{\mathfrak{L}}{}} 
\newcommand{\fS}{{\mathfrak{S}}{}} 
 
\newcommand{\fZ}{{\mathfrak{Z}}{}}

\newcommand{\n}{{\mathfrak{n}}{}} 
\newcommand{\s}{{\mathfrak{s}}{}} 
 
\newcommand{\Dto}{{\frac{{\rm d}}{{\rm d}t}|_{0}}{}}

\newcommand{\fd}{{\mathfrak{d}}{}} 
\newcommand{\h}{{\mathfrak{h}}{}} 
\renewcommand{\a}{{\mathfrak{a}}{}} 
 
\newcommand{\z}{{\mathfrak{z}}{}}
\newcommand{\CO}{{\cal O}{}}
\newcommand{\CG}{{\cal G}{}}

\newcommand{\CM}{{\cal M}{}}
\newcommand{\sldr}{{{\mathfrak{sl}(2,\mathbb R)}}{}}
\newcommand{\SLdr}{{{\mbox{\rm SL}(2,\mathbb R)}}{}}

\newcommand{\CP}{\mathcal P}
\newcommand{\CS}{\mathcal S}
\newcommand{\CE}{\mathcal E}

\newcommand{\CL}{{\cal L}{}} 
\newcommand{\CR}{{\cal R}{}} 
\newcommand{\CD}{\mathcal D}
\newcommand{\CJ}{\mathcal J}
\newcommand{\CI}{\mathcal I}
\newcommand{\CH}{\mathcal H}
\newcommand{\CB}{\mathcal B}
\newcommand{\CU}{\mathcal U}
\newcommand{\CF}{\mathcal F}

\newcommand{\CA}{\mathcal A}

\newcommand{\uC}{\underline{C}}

\newcommand{\ulambda}{\underline{\lambda}}

\newcommand{\us}{\underline{s}}

\newcommand{\Exp}{\mbox{\rm Exp}}

\title{Symmetric spaces, non-formal star products and Drinfel'd twists}

\author{{\bf Pierre Bieliavsky} (UCLouvain, Belgium)}

\date{} 

\begin{document}

\maketitle

\tableofcontents

\newpage

\chapter*{Introduction}

These notes refer to a mini-course I gave at the occasion of the conference meeting "Applications of NonCommutative Geometry to Gauge Theories, Field Theories, and Quantum Space-Time" to be held from 7 April to 11 April 2025 at the Centre International de Rencontres Mathématiques in Luminy. 

\noindent They consist in a review of some work of mine and collaborators (\cite{Bithese,B1,BM1,BG,BDS}) in the field of non-formal deformation quantization admitting a large group of symmetries. But they also contain new material and results. 

\vspace{3mm}

\noindent The motivations for this work is threefold. First, it concerns the solution of a problem raised by A. Weinstein in 1994 about geometrically describing the three-point Schwartz kernel of an invariant star-product on a symplectic symmetric space $G/K$ in terms of symplectic areas of its geodesic triangles. The problem was first solved in \cite{B1} for solvable symplectic symmetric spaces (i.e. $G$ is a solvable Lie group). Later, a solution was given in \cite{BDS} in the case of the hyperbolic plane. The method used there is called the \emph{retract method}. Chapter \ref{CHAPRETRACT} presents a wide generalisation of that method. 
\noindent However, the question of interpreting these solutions as composition formulae in the context of a $G$-equivariant symbolic calculus (``quantization map'') between classical observable defined as smooth functions on  the phase space $G/K$ and operators acting in a unitary irreducible representation of $G$ remained unclear. In particular, so was the relation between our work and A. and J. Unterberger' s work in such a framework where $G=\SLdr$ and $K=\mbox{SO}(2)$.

\noindent Second, in the same years (1993-1994), M. Rieffel introduced a way to render \emph{non-formal} the usual Abelian Moyal-Drinfel'd twist based on the vector group (or Lie algebra) $\R^{2n}$ (\cite{Ri1}). In other words, Rieffel defined a Universal Deformation Formula which, when applied to $C^{\star}$ or Fr\'echet $\R^{2n}$ algebras, produces deformations of them \emph{within the same topological class} as opposed to formal power series expansions. 
\noindent The question of defining such Universal Deformation Formulas for action of non-Abelian group was raised by Rieffel in \cite{Ri2}. Based on the explicit non-formal star-product formulae of \cite{B1}, a solution to that problem was given in \cite{BG} for actions of the Iwasawa factor $AN$ of every real simple Lie group of Hermitean type. In particular, this work provides
an interpretation of the formulae of \cite{B1} as symbolic composition formulae in a non-Abelian equivariant pseudo-differential calculus on a solvable Lie group or symmetric space. The fact that non-Abelian UDF's derive from symplectic areas of geodesic triangles in symplectic symmetric space is remarkable.

\noindent Third, the above question addressed in the solvable case remained open in the case of the hyperbolic plane. From both sides: in my context as well as in Unterberger's where no closed operator algebras were defined. In Chapter \ref{HYPSCHW} of the present review, I give a solution to that question by $\SLdr$-equivariantly quantizing  a Schwartz-type non-commutative function algebra on the hyperbolic plane by subalgebras of compact operators in the projective holomorphic discrete series of $\SLdr$.

\chapter{Symplectic symmetric spaces}
\section{Symmetric spaces}
\begin{dfn}\label{SS}\cite{L} A \emph{symmetric space}
is a pair $(M,s)$ where $M$ is a connected differentiable manifold and 
$$
s:M\times M\to M:(x,y)\mapsto s_xy
$$
is a smooth map such that

\vspace{2mm}

(1) for every point $x$ in $M$, the map $s_x:M\to M$ is involutive (i.e. $s_x^2\;=\;\id_M$) and admits the point $x$ as an isolated fixed point.
The map $s_x$ is called the \emph{symmetry} centered at $x$.

\vspace{2mm}

(2) For all $x$ and $y$ in $M$, one has
$$
s_xs_ys_x\;=\;s_{s_xy}\;.
$$
\end{dfn}
\begin{exs}\label{exss}
{\rm 
(i) Every finite dimensional real vector space $V$ is naturally a symmetric space for 
$$
s_xy\;:=\;2x-y\quad(x,y\in V)\;.
$$
(ii) More generally, every Lie group $G$ is canonically a symmetric space when equipped with the symmetries given by
$$
s_xy\;:=\;xy^{-1}x\quad(x,y\in G)\;.
$$
(iii) Let $B$ be a non-trivial scalar product (of arbitrary signature and possibly degenerated) on a finite dimensional real vector space $V$. Then every (non-null) sphere
$$
S_{B,R}\;:=\;\{x\in V\;|\;B(x,x)\;=\;R\}\quad (R\in\R_0)
$$
is naturally a symmetric space when equipped with the ``induced" symmetric structure. Indeed, for every $x$ in the sphere, the orthogonal symmetry around the
vector axis directed by $x$ is well defined and stabilises the sphere. Its restriction to the sphere then defines the symmetry centered at $x$:
$$
s_xy\;:=\;\frac{2}{R}B(x,y)x-y\quad(y\in S_{B,R})\;.
$$
}
\end{exs}
A symmetric space carries a canonical geometry. More precisely, there exists (\cite{L}) a symmetry-compatible covariant derivative on it, which admits a very simple explicit formula (\cite{Bithese, BB}):
\begin{prop}
Let $(M,s)$ be a symmetric space. Let $X,Y$ be smooth tangent vector fields on $M$ and $x$ be a point in $M$. Then, 

(i) the following 
formula
\begin{equation}\label{LC}
\left.\nabla_XY\right|_x\;:=\;\frac{1}{2}\left[X\,,\,Y\,+\,s_{x\star}Y\right]_x
\end{equation}
defines an affine connection $\nabla$ on $M$.

(ii) $\nabla$ is the unique affine connection on $M$ that is invariant under the symmetries.

(iii) It is torsion-free and its Riemann curvature tensor is parallel.

(iv) The symmetries extend the local geodesic symmetries of $\nabla$. More precisely, if $x$ is a point in $M$, and if we denote by 
$\mbox{\rm Exp}^\nabla_x: U_x\subset T_xM\to M$ the corresponding exponential map ($U_x$ denotes a symmetric open neighborhood of zero in $T_xM$ where 
the map is well defined), one has
$$
s_x\mbox{\rm Exp}^\nabla_x(v)\;=\;\mbox{\rm Exp}^\nabla_x(-v)\;.
$$
\end{prop}
\Pf First, we observe that
\begin{equation}\label{SMI}
s_{x\star x }\;=\;-\Id_{T_xM}\;.
\end{equation}
Indeed $s_{x\star x }$ is an involution of $T_xM$. Hence $T_xM$ decomposes into the direct sum of $\pm1$-eigenspaces.
Because $x$ is fixed isolately, the subspace of fixed vectors is trivial.

\noindent As a map that associates to two vector fields $X, Y\in\mbox{\rm Vec}(M)$ a third one $\nabla_XY$ it is 
clearly $\R$-bilinear. Now let $f\in C^\infty(M)$ be a smooth function on $M$. We first have
\begin{eqnarray*}
&&
\nabla_{fX}Y|_x\;=\;\frac{1}{2}\left[fX\,,\,Y\,+\,s_{x\star}Y\right]_x\\
&=&\frac{1}{2}\left(f\left[X\,,\,Y\,+\,s_{x\star}Y\right]-s_{x\star}Y(f)X-Y(f)X\right)_x\\
&=&f\nabla_XY|_x-\frac{1}{2}\left(s_{x\star}Y(f)X+Y(f)X\right)_x
\end{eqnarray*}
which equals $f\nabla_XY|_x$ in view of Property (\ref{SMI}).

\noindent Also, one has 
\begin{eqnarray*}
&&
\nabla_{X}(fY)|_x\;=\;\frac{1}{2}\left[X\,,\,fY\,+\,s_x^\star f s_{x\star}Y\right]_x\\
&=&\frac{1}{2}\left(f\left[X\,,\,Y\,\right]+(s_x^\star f )\left[X\,,\,s_{x\star}Y\right]
+X(f)Y+X(s_x^\star f )s_{x\star}Y\right)_x\\
&=&f\nabla_XY|_x+\frac{1}{2}\left(X(f)Y+X(s_x^\star f )s_{x\star}Y\right)_x\;=\;f\nabla_XY|_x
\end{eqnarray*}
as a consequence of Formula (\ref{SMI}). This proves item (i).

\noindent Regarding item (ii), we first observe that for every $y\in M$, one has 
\begin{eqnarray*}
&&
\nabla_{s_{y\star}X}(s_{y\star}Y)|_x\;=\;\frac{1}{2}\left(s_{y\star}\left[X\,,\,s_{x\star}s_{y\star}Y\,+\,Y\right]\right)_x\\
&=&
\frac{1}{2}\left(s_{y\star}\left[X\,,\,s_{y\star}s_{x\star}s_{y\star}Y\,+\,Y\right]\right)_x\;=\;
\frac{1}{2}s_{y\star s_yx}\left[X\,,\,s_{s_yx\star}Y\,+\,Y\right]\\
&=&
\left(s_{y\star}\nabla_{X}(Y)\right)_x\;.
\end{eqnarray*}

\noindent Now consider another symmetry invariant connection $\nabla'$. And 
form the tensorial quantity:
$$
S(X,Y)\;=\;\nabla_XY-\nabla'_XY\;.
$$
One then has:
\begin{eqnarray*}
&&
S(X,Y)_x\;=\;S(s_{x\star}X,s_{x\star}Y)_x\;=\;\nabla_{s_{x\star}}Xs_{x\star}Y|_x-\nabla'_{s_{x\star}X}s_{x\star}Y|_x\\
&=&
(s_{x\star}S(X,Y))_x\;=\;-S(X,Y)_x\;=\;0\;.
\end{eqnarray*}
This prove item (ii).

\noindent Expanding all the terms involved in the torsion and in the covariant derivative 
of the Riemann tensor using (\ref{SMI}), a long but straightforward computation yields item
(iii).

\noindent Item (iv) follows from the fact that, since $s_{x\star_x}^2=\id_{T_xM}$, the tangent space at $x$ splits into $(\pm1)$-eigenspaces of $s_{x\star_x}$.
But, since $s_x$ fixes $x$ isolately, the sub-space of fixed vectors must be trivial. Indeed, considering a Riemannian metric preserved
by $s_x$, every non-trivial geodesic (w.r.t. the metric) through $x$ must be flipped under $s_x$ because $s_x$ fixes 
$x$ isolately. Hence no non-trivial geodesic is tangent to the space of fixed vectors in $T_xM$.
\EPf
\begin{rmk}

\noindent The fact that the symmetries are globally defined then entails that the connection $\nabla$ is geodesically complete.
\end{rmk}

\begin{dfn}
A \emph{morphism} between two symmetric spaces $(M_i,s^{(i)})$ ($i=1,2$) is a smooth map $\phi:M_1\to M_2$
that intertwines the symmetries:
$$
\phi\circ s^{(1)}_x\;=\;s^{(2)}_{\phi(x)}\circ\phi\quad(\forall x\in M_1)\;.
$$
In the case of a diffeomorphism, one speaks about an \emph{isomorphism} between $M_1$ and $M_2$ and about 
an \emph{automorphism} in the case $M_1\;=\;M_2$.
\end{dfn}
A symmetric space turns automatically into a  homogeneous space. More precisely, one has

\begin{prop}
The group $\Aut(M,s)$ of automorphisms of a symmetric space is a Lie group of transformations of $M$ that 
acts transitively on $M$.
\end{prop}
\Pf We will first prove that the group $\mbox{\rm Aff}(\nabla)$ of affine transformations
of $(M,\nabla)$ equals $\Aut(M,s)$. 

\noindent let $\phi$ be an automorphism of $(M,s)$. Then one has 
\begin{eqnarray*}
&&
\phi_\star\nabla_XY|_x\;=\;\frac{1}{2}\phi_{\star_{\phi^{-1}(x)}}\left[X,Y+s_{\phi^{-1}(x)\star}Y\right]\\
&=&
\frac{1}{2}\phi_{\star_{\phi^{-1}(x)}}\left[X,Y+\phi^{-1}_\star s_{x\star}\phi_\star Y\right]
\;=\;\frac{1}{2}\phi_{\star}\phi^{-1}_\star\left[\phi_{\star}X,\phi_{\star}Y+s_{x\star}\phi_\star Y\right]_x
\\
&=&
\nabla_{\phi_\star X}\phi_\star Y|_x\;.
\end{eqnarray*}
The inclusion $\Aut(M,s)\subset\mbox{\rm Aff}(\nabla)$ then follows.

\noindent Regarding the other inclusion, we observe that any affine transformation 
$\psi$ of $(M,\nabla)$ locally intertwines the geodesic symmetries around 
any point $x$ of $M$:
$$
\psi s_x(\Exp_x(v))\;=\;\psi\Exp_x(-v)\;=\;\Exp_{\psi(x)}(-\psi_{\star x}(v))\;=\;s_{\psi(x)}\psi(\Exp_x(v))\;.
$$
We therefore have that there exists an open subset $U$ of $M$ where the affine transformations
$\psi s_x$ and $s_{\psi(x)}\psi$ coincides. From the fact that the curvature and torsion tensors are 
parallel, one knows \cite{KN1} [page 263 Theorem 7.7] that the manifold $M$ is analytic with respect to 
the normal coordinates. Accordingly the affine mappings are then analytic mappings.
Since the manifold is connected, $\psi s_x$ and $s_{\psi(x)}\psi$ must therefore coincide 
globally.

\noindent Since $\mbox{\rm Aff}(\nabla)$ is a Lie group of transformations of $M$ \cite{KN1} [Theorem 1.5 page 229],
so is $\Aut(M,s)$.

\noindent Let us now prove the transitivity of the action. Let $x,y$ be points in $M$ and let 
$\gamma:[0,1]\to M$ be a continuous path joining them. Since $\gamma[0,1]$ is compact,
there exists a finite number of points $\{x\;=\;x_0,...,x_N\;=\;y\}$ such that 
$U_i\;:=\;\Exp_{x_i}(T_{x_i}M)\cap\Exp_{x_{i+1}}(T_{x_{i+1}}M)\neq\emptyset$.
This yields a broken geodesic with $2N$ consecutive geodesic segments, $\{\gamma_j\}$ joining $x$ to $y$. 
Denoting by $m_j$ the midpoint of $\gamma_j$ on has $s_{m_{2N}}\circ..\circ s_{m_1}(x)\;=\;y$,
proving the assertion of transitivity.
\EPf

The affine symmetric space $M$ is therefore a homogeneous space: fixing a base point $o\in M$, and denoting 
by $H$ the stabilizer of $o$ in $\Aut(M,s)$, one has the $H$-principal fibration:
$$
\Aut(M,s)\to M\simeq\Aut(M,s)/H:g\mapsto g(o)\;.
$$
We observe that
the map 
$$
\tilde{\sigma}:\Aut(M,s)\to\Aut(M,s):g\mapsto s_ogs_o
$$
is an involutive automorphism of the Lie group $\Aut(M,s)$.

\noindent Denoting by $\aut(M,s)\subset\mbox{\rm Vec}(M)$ the Lie algebra of $\Aut(M,s)$, one has 
the decomposition 
$$
\aut(M,s)\;=\;\h\oplus\p
$$
into $\pm 1$-eigenspaces of $$\sigma\;:=\;\tilde{\sigma}_{\star e}\;=:\;\id_{\h}\oplus(-\id_\p)\;.$$
\begin{prop}
The Lie algebra of $H$ equals $\h$ and the vector space $\p$ naturally identifies with $T_o(M)$.
\end{prop}
\Pf
We let $Z$ be an element of the Lie algebra $\tilde{\h}$ of $H$. One has, for every $t\in \R$:
$$
\exp(-tZ)s_o\exp(tZ)\;=\;s_{\exp(-tZ)o}\;=\;s_o\;.
$$
Hence 
$$
\ddto\exp(-t\sigma Z)\exp(tZ)\;=\;0\;=\;Z-\sigma Z\;.
$$
that is; $Z\in\h$.

\noindent For the other inclusion, one observes that if $X\in\h$, one has 
$$
\exp(t\sigma X).o\;=\;s_o\exp(tX)s_oo\;=\;s_o\exp(tX).o\;=\;\exp(tX).o\;,
$$
hence, since $s_o$ fixes $o$ isolately: $\exp(tX)o=o$ for $t$ small i.e. $X\in\tilde{\h}$.

\noindent Denoting by 
$$
\pi:\Aut(M,s)\to M:g\mapsto g.o\;,
$$
on has that the restriction to $\p$ of $\pi_{\star e}:\aut(M,s)\to T_oM$
$$
\pi_{\star e}|_\p:\p\to T_oM
$$
is a linear isomorphism.
\EPf

\begin{dfn}
A pair $(\g,\sigma)$ where $\g$ is a finite dimensional real Lie algebra and $\sigma$ is an involutive automorphism of $\g$ is called 
an \emph{involutive Lie algebra} (briefly \emph{iLa}). 

\vspace{2mm}

\noindent A \emph{morphism} between two iLa's $(\g_i,\sigma^{(i)})$ ($i=1,2$) is a homomorphism of Lie algebras $\psi:\g_1\to\g_2$
that intertwines the involutions:
$$
\psi\circ \sigma^{(1)}\;=\;\sigma^{(2)}\circ\psi\;.
$$
When invertible, one speaks about an \emph{isomorphism}.
\end{dfn}
\begin{rmk}
Given an involutive Lie algebra $(\g,\sigma)$, the vector space $\g$ decomposes into a direct sum of the $(\pm1)$-eigenspaces
of $\sigma$:
$$
\g\;=\;\k\oplus\p\quad\quad\sigma\;=\;\id_\k\oplus(-\id_\p)\;.
$$
One then observes in this context the following relations:
\begin{eqnarray*}
\left[\k,\k\right]&\subset&\k \\
\left[\k,\p\right]&\subset&\p \\
\left[\p,\p\right]&\subset&\k \;.
\end{eqnarray*}
In particular, $\k$ turns into a Lie sub-algebra of $\g$ that acts its canonical vector direct summand $\p$.
\end{rmk}
\begin{prop}
The pair $(\aut(M,s),\sigma\;:=\;\tilde{\sigma}_{\star e})$ is an involutive Lie algebra.
\end{prop}
Given a Lie group, we denote by $\exp$ the associated exponential map.
\begin{prop}
(i) The subspace 
$$
\g(M,s)\;:=\;\k\oplus\p\quad\mbox{\rm with}\quad\k\;:=\;[\p,\p]
$$
is a Lie sub-algebra of $\aut(M,s)$.

(ii) Denoting by $G(M,s)$ the connected Lie sub-group 
of $\Aut(M,s)$ tangent to $\g(M,s)$ and by 
$$K\;:=\;G(M,s)\cap H\;,$$
the mapping
$$
\Psi:G(M,s)/K\to M:g\mapsto g.o
$$
is a diffeomorphism.
\end{prop}
\Pf
We first prove that, for every $X\in\p$, the mapping $\Exp_o:\R\to M:t\mapsto\exp(tX).o$ is the geodesic 
by $o$ tangent to $X$: 
$$
\Exp_o(tX)\;=\;\exp(tX).o\quad\mbox{\rm where}\quad\p\;=\;T_oM\;.
$$ 
To see this, we consider the vector field 
$$
X^\star_x\;:=\;\ddto\exp(-tX).x\;.
$$
We have
\begin{eqnarray*}
\nabla_{X^\star} X^\star|_{\exp(tX).o}\;=\;\frac{1}{2}\left[X^\star,s_{\exp(tX).o_\star}X^\star\right]_{\exp(tX).o}\;.
\end{eqnarray*}
But, 
\begin{eqnarray*}
&&
-\left(s_{\exp(tX).o_\star}X^\star\right)_x\;=\;s_{\exp(tX).o_\star s_{\exp(tX).o}(x)}X^\star\\
&=&
\left.\frac{\rm d}{{\rm d}s}\right|_{s=0}s_{\exp(tX).o}\left(\exp(sX)s_{\exp(tX).o}(x)\right)\\
&=&
\left.\frac{\rm d}{{\rm d}s}\right|_{s=0}\exp(tX)s_o\exp(-tX)\exp(sX)\exp(tX)s_o\exp(-tX)(x)\\
&=&
\left.\frac{\rm d}{{\rm d}s}\right|_{s=0}\exp(tX)s_o\exp(sX)s_o\exp(-tX)(x)\\
&=&
\left.\frac{\rm d}{{\rm d}s}\right|_{s=0}\exp(tX)\exp(-sX)\exp(-tX)(x)\;=\;X^\star_x\;.
\end{eqnarray*}
Hence the assertion.

\noindent Now let $x\in M$ and $\{\gamma_j\}_{j=0,...,N}$ be a finite sequence of consecutive geodesics joining $o$ to $x$.
At $\gamma_0(1)\;=\;\exp(X_0).o$ every geodesic $\gamma_1$ is a translate by $\exp(X_0)$ of a geodesic at $o$ of the 
form $\exp(tY).o$ ($Y\in\p$):
$$
\gamma_1(t)\;=\;\exp(X_0)\exp(tY).o\;=\;\exp(t\Ad_{\exp X_0}Y)\exp(X_1).o\;.
$$
In particular, there exists an element $X_1\in\g(M,s)$ such that $$\gamma_1(t)=\exp(tX_1)\exp(X_0).o\;.$$
Iterating, we find elements $X_0,...,X_N\in\g(M,s)$ such that
$$
x\;=\;\exp(X_N)...\exp(X_0).o\;.
$$
The group $G(M,s)$ therefore acts transitively on $M$.
\EPf
\begin{dfn}
The group $G(M,s)$ is called the \emph{transvection group} of the symmetric space $(M,s)$.
\end{dfn}
\begin{rmk}
{\rm 
The transvections generalize to symmetric spaces the notion of translation: think about the example of a vector space $V$, where the 
product of two Euclidean centered symmetries indeed defines a translation.
}
\end{rmk}
From all this, one encodes the geometry of the symmetric space $M$ by the purely Lie theoretic notion  of involutive Lie algebra.
We  describe this process below.

\begin{dfn}
A \emph{transvection} iLa is an iLa $(\g,\sigma)$ enjoying the following two properties:

(1) $[\p,\p]=\k$ and

(2) the representation $\k\times\p\to\p:(X,Y)\mapsto[X,Y]$ is faithful. 
\end{dfn}
\begin{prop}
Let $(M,s)$ be a symmetric space. The pair $(\g(M,s),\sigma\;:=\;\tilde{\sigma}_{\star e}|_{\g(M,s)})$
is a transvection iLa. In particular one has a correspondance between symmetric spaces and 
transvection iLa's: $$(M,s)\mapsto (\g(M,s),\sigma)\;.$$.
\end{prop}
\Pf
The only thing to prove is faithfulness, which follows immediately from the fact that $\Aut(M,s)$
is constituted of transformations of $M$ which by definition do act faithfully. 
\EPf

\begin{thm}\label{CORRSS}
The correspondance described above that associates to every symmetric space a transvection iLa induces a bijection between the isomorphism classes
of simply connected symmetric spaces and the isomorphism classes of transvection iLa's.
\end{thm}
\Pf
We let $(M,s)$ and $(M',s')$ be two simply connected symmetric spaces that are isomorphic under the isomorphism $\phi:M\to M'$.
We observe the  isomorphism between their automorphism groups:
$$
\hat{\phi}: \Aut(M,s)\to\Aut(M',s'):\varphi\mapsto\phi\varphi\phi^{-1}\;.
$$
We choose a base point $o\in M$ and denote by $\phi(o)=:o'$ the corresponding one in $M'$. These choices 
define involutions $\tilde{\sigma}$ and $\tilde{\sigma}'$ where
$$
\tilde{\sigma}'\;=\;\hat{\phi}\tilde{\sigma}\hat{\phi}^{-1}\;.
$$
One also has the  decompositions:
$$
\aut(M,s)\;=\;\h\oplus\p\quad\mbox{\rm and}\quad\aut(M',s')\;=\;\h'\oplus\p'
$$
according respectively to
$$
\sigma\quad\mbox{\rm and}\quad\sigma'\;:=\;\tilde{\sigma}'_{\star\id}\;=\;
\hat{\phi}_{\star\id}{\sigma}\hat{\phi}^{-1}_{\star\id}\;.
$$
Now, setting $\underline{\phi}\;:=\;\hat{\phi}_{\star\id}$, one has 
$$
\p'\;=\;\underline{\phi}(\p)\;.
$$
Hence, since $\underline{\phi}$ is an isomorphism of Lie algebras, its restriction
to $\g(M,s)$ realizes an isomorphism with $\g(M',s')$ that intertwines the involutions.

\vspace{2mm}

\noindent We now revert the correspondance. We fix a transvection iLa $(\g\;=\;\k\oplus\p,\sigma)$.
We will associate to it a simply connected symmetric space. In order to do so, we consider 
the simply connected Lie group $G$ whose Lie algebra is $\g$ as well as the connected 
Lie subgroup $K$ of $G$ tangent to $\k$. The long homotopy sequence of the fibration 
$K\hookrightarrow G\to G/K$ tells us that, since $K$ is connected and $G$ is simply connected, $G/K$
is simply connected. This cosetspace is canonically a symmetric space. Indeed, considering 
$\tilde{\sigma}$, the involutive automorphism of $G$ tangent to $\sigma$, one sets
$$
s_{gK}(g'K)\;:=\;g\tilde{\sigma}(g'g^{-1})K\;.
$$
The above formula defines a structure of symmetric space on $G/K$. Indeed, denoting by 
$\pi:G\to G/K$ the canonical projection, we have for every $X\in\p$:
$$
s_{K_{\star K}}\pi_{\star e}X\;=\;\ddto s_K\pi(\exp(tX))\;=\;\ddto\pi\tilde{\sigma}(\exp(tX))\;=\;-\pi_{\star e}X\;.
$$
Hence $s_{K_{\star K}}\;=\;-\id_{T_KM}$ which implies by homogeneity that $s_{gK}$ fixes $gK$
isolately. The other axioms are easy to check. 

Now, we consider the Lie group homomorphism:
$$
\tau:G\to\Aut(G/K,s):g\mapsto\tau_g\;:=\;[g'K\mapsto gg'K]\;.
$$
In general, this homomorphism is not injective. However, its differential is.
Indeed, setting, at the level of $\Aut(G/K,s)$,  $\hat{\sigma}(\varphi)\;:=\;s_K\varphi s_K$, one has 
\begin{eqnarray*}
&&
\hat{\sigma}_{\star\id}\tau_{\star e}(X)|_{gK}=\ddto s_K\tau_{\exp(tX)}s_K(gK)=
\ddto \tilde{\sigma}({\exp(tX)}\tilde{\sigma}(g))\\
&=&
\tau_{\star e}\sigma(X)|_{gK}\;.
\end{eqnarray*}
Which implies that, denoting by $\aut(G/K,s)\;=\;\hat{\h}\oplus\hat{\p}$ the iLa 
decomposition, one has 
$$
\tau_{\star e}\p\subset\hat{\p}\;.
$$
A dimensional argument on the dimension of the homogeneous space $G/K$ implies 
that 
$$
\tau_{\star e}\p=\hat{\p}\;.
$$
The fact that the action of $\k\;=\;[\p,\p]$ on $\p$ is faithful and the fact that $\tau_{\star e}$
is a homomorphism of Lie algebras imply that $\tau_{\star e}|_\k$ is injective as well.
Hence $\tau_{\star e}$ defines an isomorphism between $(\g,\sigma)$ and the transvection
iLa of $(G/K,s)$.

\noindent The rest is obvious. \EPf
\begin{ex}\label{exlgss}
{\rm 
Viewing a Lie group $\L$ as a symmetric space, one associates to it the so-called ``exchange case" iLa. Denoting by $\fL$
the Lie algebra of $\L$, one sets 
$$
\fL^2\;:=\;\fL\oplus\fL\quad\mbox{\rm and}\quad\sigma:\fL^2\to\fL^2:(X,Y)\mapsto(Y,X)\;.
$$
The pair $(\fL^2,\sigma)$ is then an iLa with canonical eigensubspaces
\begin{eqnarray*}
\k_2\;=\;\{(X,X)\}\;\simeq\;\fL\\
\p_2\;:=\;\{(X,-X)\}\;.
\end{eqnarray*}
The action of $\k_2$ on $\p_2$ is generally not faithful and one generally does not have $[\p_2,\p_2]=\k_2$. However, one may recover the transvection iLa
by considering
$$
\g\;:=\;\k\oplus \p
$$
with 
$$
\k\;:=\;\fL'/\z\quad\mbox{\rm and}\quad\p\;:=\;\p_2
$$
where $\fL'\;:=\;[\fL,\fL]$ denotes the first derivative of $\fL$  and where $\z$ the centralizer of $\fL$ in $\fL'$ (w.r.t. the adjoint representation):
$$
\z\;:=\;\fL'\cap\z(\fL)\;.
$$
}
\end{ex}

\begin{dfn}
The symmetric space $M$ is called \emph{group type} when there exists a Lie subgroup $\S$ of $G$ that 
acts simply transitively on $M$ i.e. one has the $\S$-equivariant diffeomorphism:
$$
\S\to M: x\mapsto xK\;.
$$
One then talk about a group type symmetric space on the abstract Lie group $\S$.
\end{dfn}
We now consider a Lie group $\S$ equipped with a smooth one-parameter family of group type symmetric space structure (hence left-invariant):
$$
s^{(\kappa)}:\S\times\S\to\S \quad(\kappa\in I)
$$
where $I$ is a connected interval of $\R$.

\noindent We denote by $\g_\kappa$ the associated family of transvection Lie algebras and by
$$
\g_\kappa\;=\;\k_\kappa\oplus\p_\kappa
$$
the corresponding symmetric space decomposition. In that context we make the following definition.
\begin{dfn}
For $\kappa$ and $\kappa'$, we say that $s^{(\kappa')}$ is a \emph{curvature contraction} of $s^{(\kappa)}$ if 
$$
\dim\fZ_{\p_{\kappa'}}(\k_{\kappa'})\;>\;\dim\fZ_{\p_\kappa}(\k_\kappa)
$$
where we denote by $\fZ_\p(\k)$ the centralizer of $\k$ in $\p$.
\end{dfn}
This essentially says that the curvature of $s^{(\kappa')}$ is obtained from the one of $s^{(\kappa)}$ by letting tend to 
zero part of it curvature components.

\section{Symplectic symmetric spaces}
\begin{dfn}\cite{Bithese,BCG}
A \emph{symplectic symmetric space} (briefly SSS) is a triple $(M,s,\omega)$ where $(M,s)$ is a symmetric space and where $\omega$ is a non-degenerate 2-form on $M$ that is symmetry invariant:
$$
s_x^\star\omega\;=\;\omega\quad(\forall x\in M)\;.
$$
Two such symplectic symmetric spaces are \emph{isomorphic } if there exists a symplectomorphism 
between them which is at the same time an isomorphism of symmetric spaces between the underlying 
symmetric spaces. 
\end{dfn}
The following observation is an immediate consequence of \cite{KN1} (page ??)
\begin{prop}
The two form $\omega$ is parallel w.r.t. the canonical affine connection $\nabla$:
$$
\nabla\omega\;=\;0\;.
$$
In particular, it is closed:
$$
{\rm d}\omega\;=\;0
$$
hence symplectic.
\end{prop}

\begin{ex} {\bf Hermitean symmetric spaces}

\vspace{3mm}

{\rm
We start with a simple Lie group $G$ with finite center whose Lie algebra $\g$
is absolutely simple (i.e. $\g^\C\;:=\;\g\otimes\C$ is simple). We consider a maxi lam compact subgroup $K\subset G$ whose Lie algebra $\k$ is assumed to admit a non-trivial center $\z(\k)$.
Such a Lie group $G$ is called of \emph{Hermitean type}.

\noindent In that case, $\z(\k)$ is known (\cite{K,Bithese}) to be one-dimensional and generated by an element 
$Z_0$ that acts as a complex structure on $\g/\k$.

\noindent More precisely, denoting by $B$ the Killing form on $\g$, the 
 decomposition
$$
\g\;=\;\k\,\oplus\,(\p\;:=\;\k^{\perp_B})
$$
underlies an iLa $(\g,\theta)$ whose involution $\theta$ is a Cartan involution 
of $\g$.

\noindent The action of $Z_0$ on $\p$:
$$
\CJ\;:\p\to\p:X\mapsto[Z_0,X]
$$
is then a complex structure:
$$
\CJ^2\;=\;-\id_{\p}\;.
$$
One then gets the symplectic element $\Omega\in\bigwedge^2(\p^\star)$ defined by
$$
\Omega(X,Y)\;:=\;B(\CJ X,Y)\;,
$$
which on the coset space $M\;:=\;G/K$ defines a $G$-invariant symplectic structure 
$\omega$ throughout the formula:
\begin{equation}\label{SYMPF}
\omega_{gK}(X^\star,Y^\star)\;:=\;\Omega\left((Ad_{g^{-1}}X)_\p\,,\,(Ad_{g^{-1}}Y)_\p\right)\quad(X,Y\in\p)
\end{equation}
where we write for every $T\in\g$
$$
T\;=\;T_\k+T_\p
$$
according to the Cartan decomposition.

\noindent As seen above, the homogeneous space $M$ is canonically equipped with a symmetric space 
structure
$$
s_{gK}g'K\;:=\;g\theta(g^{-1}g')K
$$
where we (abusively) denote by $\theta$ the involutive automorphism of $G$ that integrates
the Cartan involution $\theta$ of $\g$. The triple $(M,s,\omega)$ is then a symplectic symmetric space.
}
\end{ex}
\begin{rmks}
(i) The Iwasawa decomposition of the semi-simple Lie group $G\;=\;ANK$ implies that these
Hermitean symmetric spaces are all of group type.

(ii) In the above example, the formula:
$$
J_{gK}(X^\star)\;:=\;\CJ(Ad_{g^{-1}}X)_\p\quad(X\in\g)
$$
defines a $G$-invariant complex structure on $M=G/K$.
The element
$$
\beta(X^\star,Y^\star)\;:=\;-\omega(\CJ X,Y)
$$
then defines a $G$-invariant Kahler metric on $M$.

\end{rmks}

\noindent The infinitesimal structure, analogue to the one of iLa, that encodes the notion of \emph{symplectic} 
symmetric space is the following one.
\begin{dfn}
A \emph{symplectic triple} (briefly ST) is a triple $(\g,\sigma,\Omega)$ where $(\g\;=\;\k\oplus\p,\sigma)$ is an iLa and where $\Omega$ is a Chevalley 2-cocycle
 (w.r.t the trivial representation of $\g$ on $\R$):
 $$
 \Omega\in\bigwedge^2(\g^\star): \delta\Omega\;=\;0
 $$
 that vanishes identically on $\k\times\k$ and $\k\times\p$ and
 restricts to $\p\times\p$ as a non-degenerate 2-form.
 
 \vspace{2mm}
 
 \noindent Two such symplectic triples are \emph{isomorphic} when there exists an isomorphism 
 of iLa's between the underlying iLa's that intertwines the 2-cocycles.
\end{dfn}
\begin{rmk}
{\rm 
In the above definition, the Chevalley coboundary operator is defined as usual by
$$
\delta\Omega(X,Y,Z)\;:=\;\Omega([X,Y],Z)+ \Omega([Z,X],Y)+\Omega([Y,Z],X)\;.
$$
}
\end{rmk}

\noindent One then has
 \begin{thm}
The correspondence between the symmetric spaces and the iLa's (c.f. Theorem \ref{CORRSS})  induces a bijection between the set of isomorphisms classes of simply connected SSS's and the set of isomorphism classes of transvection ST's. 
 \end{thm}
 \Pf
 To a SSS $(M,s,\omega)$ one associates a ST  in the obvious way:
$(\g,\sigma)$ is defined as the transvection iLa associated to the symmetric space $(M,s)$, while,
through the isomorphism
$$
\pi_{\star_0}|_\p:\p\to T_oM
$$
(with $o\in M$ such that $\sigma_{o_{\star e}}\;=\;\sigma$), one defines $\Omega$ to be the extension
by zero on $\k\times\k$ and $\k\times\p$ of $(\pi_{\star_0}|_\p)^\star\omega_o$. One readily checks that the 
obtained triple $(\g,\sigma,\Omega)$ is a SST.

\noindent Conversely, to a ST $(\g,\sigma,\Omega)$ the proof of Theorem \ref{CORRSS} associates
a symmetric space modelled on $G/K$. The formula (\ref{SYMPF}) then defines the symplectic structure
needed to from a SSS $(G/K,s,\omega)$.

\noindent The rest is immediate following the lines of the proof of Theorem \ref{CORRSS}.\EPf

Due to an old result of B. Kostant, a symplectic symmetric space is, generally, locally a co-adjoint orbit. But in general not of its transvection group. However, in this book, we will be concerned with the following sub-class of SSS's.
\begin{dfn}
A SSS $(M,s,\omega)$ si called \emph{Hamiltonian} when its transvection group $G$ acts 
on $(M,\omega)$ is a Hamiltonian way, that is, when equipping $C^\infty(M)$ with the Poisson structure 
associated with the symplectic structure, there exists a $G$-equivariant homomorphism of Lie algebras
$$
\lambda:\g\to C^\infty(M):X\mapsto\lambda_X
$$
such that
$$
{\rm d}\lambda_X\;=\;i_{X^\star}\omega\;.
$$
where $i$ denotes the usual interior product of differential forms by vector fields.
\end{dfn}
\begin{rmk} In that case, a general result due to Kostant about symplectic homogeneous spaces
states that the map
$$
m:M\to\g^\star:x\mapsto[X\mapsto\lambda_X(x)]
$$
realizes a $G$-equivariant symplectic covering of $M$ over a co-adjoint orbit $\CO$ of $G$ in $\g^\star$.
\end{rmk}

\begin{rmk}
For \emph{semi-simple} symplectic symmetric space i.e. when $G$ is semi-simple, an easy application of Whitehead's lemmas implies that's it always the case.
\end{rmk}
\begin{rmk}
To a Hamiltonian SSS $(M,s.,\omega)$ corresponds such a transvection ST $(\g,\sigma,\Omega)$
with \emph{exact} symplectic element
$$
\Omega\;=\;\delta\xi \quad(\xi\in \g^\star)\;.
$$
\end{rmk}
\begin{ex}\label{PLP} {\bf Contracting the Lobatchevsky plane}

{\rm 

\noindent In this paragraph, we will express the geometry of the group $\S\;=\;ax+b$ as a curvature contraction 
of the hyperbolic plane. 

\vspace{2mm}

\noindent We now consider the following realization of the connected component of 2-di\-men\-sio\-nal Lie group $\S$ of affine transformations
of the real line.  We consider the two dimensional Lie algebra $\s$ generated by two elelments $H$ and $E$, with Lie bracket:
$$
[H,E]\;:=\;2E\;.
$$
The corresponding connected and simply connected  Lie group $\S$ is then coordinated by the global chart
$$
\s\to\S:aH+nE\mapsto \exp(aH).\exp(nE)\;=:\;(a,n)
$$
in which the group law is given by
$$
(a,n).(a',n')\;:=\;(a+a'\,,\,e^{-2a'}n\,+\,n')
$$
with unit
$$
e\;:=\;(0,0)
$$
and inverse
$$
(a,n)^{-1}\;=\;(-a,-e^{2a}n)\;.
$$
Within this chart, a left-invariant Haar volume is given by
$$
\omega\;:=\;{\rm d}a\wedge{\rm d}n\;.
$$
Within this context we have
$$
\tilde{H}\;=\;\partial_a\,-\,2n\partial_n\;\quad\tilde{S}\;=\;\partial_n\;.
$$
\begin{prop}
(i) Setting
\begin{eqnarray*}
&&
s^{(\kappa)}_e:\S\to\S:(a,n)\mapsto(\,-a-\frac{1}{2}\log(1+\kappa n^2)\,,\,-n \,)\\
&&\mbox{and}\quad s^{(\kappa)}_x\;:=\;L_x s^{(\kappa)}_eL_{x^{-1}}\quad \kappa\in[0,\infty[
\end{eqnarray*}
defines  a smooth family of symplectic symmetric spaces on $(\S,\omega)$.

\vspace{3mm}

\noindent (ii) For every $\kappa>0$, the symplectic symmetric space $(\S,s^{(\kappa)},\omega)$
is isomorphic (as SSS) to the hyperbolic plane
$$
\D\;:=\;\SLdr/\mbox{\rm SO}(2)
$$
with (normalized) Kahler symplectic form $\omega$.

\vspace{3mm}

\noindent (iii) For $\kappa>0$, the symplectic symmetric space $(\S,s^{(0)},\omega)$ 
is a curvature contraction of the symplectic symmetric space $(\S,s^{(\kappa)},\omega)$.

\vspace{3mm}

\noindent (iv) As symmetric space,  $(\S,s^{(0)})$ is isomorphic to the canonical Lie group symmetric space structure on $\S$.

\vspace{3mm}

\noindent (v) The transvection Lie algebra $\g_0$ of $(\S,s^{(0)})$ is isomorphic to the Lie algebra of the Poincar\'e group $SO(1,1)\ltimes\R^2$. 

\vspace{3mm}

\noindent (vi) As symplectic symmetric space,  $(\S,s^{(0)},\omega)$ is canonically isomorphic to the symplectic symmetric space structure on a co-adjoint orbit of  $SO(1,1)\ltimes\R^2$ in 
$\g_0^\star$.

\end{prop}
\Pf
We consider the family of Lie algebra structures $\{\g_\kappa\}_{\kappa\in\R^+}$
over $\R^3\;:=\;\span_\R\{H,E,F\}$ given by the Lie bracket table:
$$
\begin{array}{ccc}
  \left[H,E\right]& =  &2E   \\
  \left[H,F\right]& =  & -2F  \\
  \left[E,F\right]& =  &\kappa H\;.  
\end{array}
$$
When $\kappa\neq0$, $\g_\kappa$ is isomorphic to $\g_0\;\simeq\;\sldr$. 

\noindent There, the element
$Z\;:=\;E-F$ generates a line $\k_\kappa$ that consists in a maximal compact Lie subalgebra
of $\g_\kappa$. 

\noindent In basis $\{H,E,F\}$, the Killing form $B_\kappa$ takes the following matrix form:
$$
\left[B_\kappa\right]\;=\;
\left(
\begin{array}{ccc}
  8&0   &0   \\
  0& 0  &   4\kappa\\
  0&4\kappa  &0   
\end{array}
\right)\;.
$$
The dual metric $B^\star_\kappa$ on $\g^\star_\kappa$ defined by transporting $B_\kappa$
from $\g_\kappa$ to $\g^\star_\kappa$ under the equivariant identification
$$
{}^\sharp:\g^\star_\kappa\to\g_\kappa:\Xi\mapsto{}^\sharp\Xi
$$
$$
B_\kappa({}^\sharp\Xi\,,\,X)\;:=\;<\Xi,X>
$$
has, within the dual basis $\{H^*,E^*,F^*\}$ of $\{H,E,F\}$, the following matrix form:
$$
\left[B^\star_\kappa\right]\;=\;\frac{1}{4\kappa}
\left(
\begin{array}{ccc}
  \frac{\kappa}{2}&0   &0   \\
  0& 0  &   1\\
  0&1   &0   
\end{array}
\right)\;.
$$
Rescaling, we therefore obtain of the dual $\g^\star_\kappa$ a co-adjoint-invariant scalar product $\eta_\kappa$ whose matrix form in basis $\{H^*,E^*,F^*\}$ is given by
$$
\left[\eta_\kappa\right]\;:=\;
\left(
\begin{array}{ccc}
  \frac{\kappa}{2}&0   &0   \\
  0& 0  &   1\\
  0&1   &0   
\end{array}
\right)\;.
$$
The non-null spheres $S_{t,R}$ for $\eta_\kappa$ are therefore symmetric spaces
whose common transvection group is realized by the adjoint group
$G_\kappa$ of $\g_\kappa$ (isomorphic to $\mbox{\rm P}\SLdr$) acting on $S_{\kappa,R}$
by the restricted co-adjoint action.

\noindent Each $S_{\kappa,R}$ then turns into a SSS when equipped with its canonical symplectic 
form ---the so-called \emph{KKS form}--- associated to its nature of co-adjoint orbit.

\noindent Observing that the element $o\;:=\;E^*-F^*$ is fixed by $K_\kappa$ (for $(E^*-F^*)\circ\ad_{Z}\;=\;0$), we get that, within this picture, the hyperbolic plane $\D$ is realized by the connected component of  $S_{\kappa,-2}$.

\noindent Setting
$$
\begin{array}{ccc}
 \n & :=  & \R E  \\
 \a & :=  & \R H  \\
 A & :=   &\exp\a\\   
 N & :=   &\exp\n\\   
 \S & :=   &A.N \\
 K_\kappa & :=   &\exp\k_\kappa\;,\\      
\end{array}
$$
the Iwasawa decomposition
$$
G_\kappa\;=\;\S K_\kappa
$$
yields a global diffeomorphism
$$
\S\to\D\subset\g^\star_\kappa:x\mapsto x.o\;.
$$
The above diffeomorphism is $\S$-equivariant when endowing the source with 
the action of $\S$ on itself by left-translations.

\noindent Coordinating the Iwasawa factor $\S$ via
$$
\R^2\to\S:(a,n)\mapsto\exp(aH)\exp(nE)
$$
or, equivalently, its orbit $\D\subset S_{\kappa,-2}$ under
$$
\R^2\to\D:(a,n)\mapsto \exp(aH)\exp(nE).o\;,
$$
a small computation realizes the symmetry $s_o$ at the level of $\S$ as
$$
\S\to\S:(a,n)\mapsto(\,-a-\frac{1}{2}\log(1+\kappa n^2)\,,\,-n \,)\;.
$$
The coordinate system $(a,n)$ enjoys the property of being \emph{Darboux}, in the sense
that the following 2-form:
$$
\omega\;:=\;{\rm d}a\wedge{\rm d}n\;.
$$
is left-invariant. It therefore corresponds to a normalization of both a Haar volume and the KKS form.

\vspace{3mm}

\noindent The limit $\kappa\to0$ yields now a curvature contraction of the hyperbolic plane $(\D,\us,\omega)$ into a symplectic symmetric space denoted hereafter $(\M,s,\omega)$  whose transvection
Lie algebra $\g_0$ is isomorphic to the Poincar\'e Lie algebra $so(1,1)\ltimes\R^2$.

\noindent The spheres $S_{\kappa,-2}$ contract to spheres $S_{0,-2}$ in the inner product space $(\g_0^\star,\eta_0)$. Warn the fact that in this limit $\kappa=0$ the scalar product $\eta_0$ is now degenerated. 
\noindent As co-adjoint orbits in $\g_0^\star$, these spheres $S_{0,-2}$ corresponds to ``massive" orbits of the Poincar\'e group $SO(1,1)\ltimes\R^2$. 

\noindent Regarding the symplectic symmetric space structure on $S_{\kappa,-2}$, in the limit $\kappa=0$, one has the transvection iLa $(\g_0,\sigma)$ with non-vanishing brackets:
$$
\begin{array}{ccc}
  \left[H,E\right]& =  &2E   \\
  \left[H,F\right]& =  & -2F
\end{array}
$$
 and ``contracted Cartan involution":
 $$
 \sigma H=-H\quad\sigma E=-F\;.
 $$
 The associated eigenspace decomposition is given by
 $$
 \g_0\;=\;\k_0\oplus\p_0
 $$
 with 
 $$
 \k_0\;=\;\span\{H\,,\,E+F\}\quad\p_0\;=\;\span\{Z\;:=\;E-F\}\;.
 $$
 At last, we observe that the transvection iLa $(\g_0,\sigma)$ is isomorphic to the transvection
iLa of the Lie group $\S\;=\;AN$ of affine transformations of the real line.
Indeed, denoting by $\s\;:=\;\a\oplus\n$ its Lie algebra, from what was described above in Example
\ref{exss}, we 
have 
$$
\fL\;=\;\s\quad,\,\z\;=\;\{0\}\quad,\,\fL'\;=\;\R(E,E)\;.
$$
The iLa isomorphism is then given by
$$
\begin{array}{ccc}
\fL'\oplus\p_2&\longrightarrow&\k_0\oplus\p_0\\
(E,E)&\mapsto&E-F\\
(E,-E)&\mapsto&E+F\\
(H,-H)&\mapsto&H\;.\\
\end{array}
$$ \EPf

}
\end{ex}

\chapter{The retract method in Formal Deformation Quantization}\label{CHAPRETRACT}
\section{Star-products on symplectic manifolds}
\begin{dfn}\label{STARPROD}
Let $(M,\omega)$ be a symplectic manifold with associated symplectic Poisson bracket denoted by $\{\,,\,\}^\omega$. Let us denote by $C^\infty(M)[[\nu]]$
the space of formal power series with coefficients in $C^\infty(M)$ in the formal parameter $\nu$. A \emph{star-product}
on $(M,\omega)$ is a $\C[[\nu]]$-bilinear associative algebra structure on 
$C^\infty(M)[[\nu]]$:
$$
\star_\nu: C^\infty(M)[[\nu]]\times C^\infty(M)[[\nu]]\to C^\infty(M)[[\nu]]
$$
such that, viewing $C^\infty(M)$ embedded in $C^\infty(M)[[\nu]]$ as 
the zero order coefficients, and writing for all $u,v\in C^\infty(M)$:
$$
u\star_\nu v\;=:\;\sum_{k=0}^\infty\nu^kC_k(u,v)\;,
$$
one has

\vspace{3mm}

(1) $C_0(u,v)\;=\;uv$ (pointwise multiplication of functions).

\vspace{3mm}

(2) $C_1(u,v)-C_1(v,u)\;=\;\{u,v\}^\omega$.

\vspace{3mm}

(3) For every $k$, the coefficient $C_k:C^\infty(M)\times C^\infty(M)\to C^\infty(M)$
is a bi-differential operator. These are called the \emph{cochains} of the 
star-product.

\vspace{3mm}

\noindent A star-product is called \emph{natural} when $C_2$ is 
a second order bi-differential operator (cf. \cite{L}).
\end{dfn}

\begin{dfn}\label{EQUIVSTARDEF}
Two such star-products $\star^j_\nu$ ($j=1,2$) are called \emph{equivalent}
if there exists a formal power series of differential operators of the form
\begin{equation}\label{EQUIVSTAR}
T\;=\;\id\;+\;\sum_{k=1}^\infty\nu^kT_k\;:\;C^\infty(M)\to C^\infty(M)[[\nu]]
\end{equation}
such that for all $u,v\in C^\infty(M)$:
\begin{equation}\label{EQUIVSTAR2}
T(u)\star^2_\nu T(v)\;=\;T(u\star^1_\nu v)\;,
\end{equation}
where we denote by $T$ the $\C[[\nu]]$-linear extension of (\ref{EQUIVSTAR}) to 
$C^\infty(M)[[\nu]]$.
\end{dfn}
Since the zeroth order of such an equivalence operator $T$ is the identity, hence invertible, the followingl argument implies that $T$ itself is invertible
 as a $\C[[\nu]]$-linear operator on $C^\infty(M)[[\nu]]$. 
 \begin{lem}\label{INV}
Let $(\A,I)$ be a unital associative algebra over $\R$ and  $A\;:=\;I\,+\,\sum_{k=1}^\infty\nu^k A_k$ be an element of $\A[[\nu]]$.
Then $A$ admits an inverse $A^{-1}$ in $\A[[\nu]]$.
\end{lem}
\Pf
One readily check that $A^{-1}\;=\;I\,+\;\sum_{k=1}^\infty \nu^pA_p'$ where 
the $A_p'$ are defined as solutions of the recurrence:
$$
A'_p\;=\;-A_p\,-\,\sum_{n+k=p;\,n,k\geq1}A_kA'_n\;.
$$
\EPf
 \noindent One denotes its inverse by
 $$
 T^{-1}:C^\infty(M)[[\nu]]\to C^\infty(M)[[\nu]]\;.
 $$
 We the adopt the following notation.
 \begin{dfn}
 Let $\star^j_\nu$ ($j=1,2$) be two star-products on $(M,\omega)$ that are equivalent 
 to each other under an equivalence operator $T$ as in Definition (\ref{EQUIVSTARDEF}). We encode equation (\ref{EQUIVSTAR2}) by the notation
 $$
\star^2_\nu\;=\;T(\star_\nu^1)\;.
 $$
 \end{dfn}
 Let now $G$ be a group of transformations of $M$.
 \begin{dfn}
 A star-product $\star_\nu$ on $(M,\omega)$
 is called $G$-\emph{invariant} if its cochains are $G$-invariant bi-differential
 operators; for every $g\in G$, we write:
 $$
 g^\star(u\star v)\;=\;g^\star(u)\star_\nu g^\star(v)\quad (u,v\in C^\infty(M))\;.
 $$
 Two $G$-invariant star-products $\star^j_\nu$ ($j=1,2$) are called 
 $G$-\emph{equivariantly equivalent} if there exists a $G$-commuting 
 equivalence $T$ between them:
 $$
 T(\star^1_\nu)\;=\;\star^2_\nu\quad \mbox{\rm and}\quad T(g^\star u)\;=\;g^\star T(u)
 $$
 for all $g\in G$ and $u\in C^\infty(M)$.
 \end{dfn}
 \begin{thm}\label{EXISTSTAR}
 Let $(M,\omega)$ be a symplectic manifold. Let $G$ be a group of symplectic transformations of $(M,\omega)$ that preserves an affine connection on $M$. Then, the set of $G$-equivariant
 equivalence classes of $G$-invariant star-products on $(M,\omega)$ is 
 canonically parametrized by the space of formal power series, $H^2_{dR}(M)^G[[\nu]]$,
with coefficients in the $G$-invariant second de Rham cohomology space of $M$ 
 \end{thm}
 In the above statement, we mention the $G$-invariant de Rham complex.
 By this, in degree two, we mean the following. The de Rham differential
 in degree one restricts to the differentiable one-forms that are $G$-invariant and 
 ranges in the closed two-forms that are also $G$-invariant. The cohomology 
 space associated to this restriction defines $H^2_{dR}(M)$.
 
 \begin{rmk}
 {\rm 
 Theorem \ref{EXISTSTAR} summarizes a long series of works by many people  in deformation quantization of symplectic manifolds. The existence of star-products on symplectic manifolds has essentially three sources : \cite{DL,OMY,F}. Where, to my opinion, the third one
 (\cite{F}) is the most constructive and suitable for applications. Classification under equivalence resulted from a combination of \cite{F}, \cite {G} and \cite{Be}. Published later and 
 independently, one finds it also in \cite{NT}.
 
 \noindent For the affine $G$-invariant case, see \cite{BBG}.
 
 }
 \end{rmk}
 We end this section with
 \begin{dfn}
 A \emph{derivation} of a star-product $\star$ on $(M,\omega)$ is a
 formal series
 $$
 D\;=\;\sum_{k=0}^\infty\nu^kD_k
 $$
  of differential operators $\{D_k\}$ on $C^\infty(M)$ such that for all $u,v\in C^\infty(M)$, one has
  $$
  D(u\star v)\;=\;(Du)\star v+u\star(Dv)\;.
  $$
  The (Lie algebra) of derivations of $\star$ is denoted by $\Der(\star)$.
 \end{dfn}
 \begin{prop}\label{DERVGUTT}\cite{RG}
 When $M$ is simply connected, every derivation $D$ of $\star$ is \emph{interior} in the sense
 that there exists an element $\lambda_D\in C^\infty(M)[[\nu]]$ such that
 $$
 \nu Du\;=\;[\lambda_D\,,\,u]_{\star} \quad(\forall u)
 $$
 where we set
 $$
 [u,v]_{\star}\;=\;u\star v-v\star u\;.
 $$
 \end{prop}

\section{Drinfel'd twists and symplectic Lie groups}

We consider a Lie group $\L$ with Lie algebra $\fL$ whose we denote by $\CU(\fL)$
the enveloping algebra. The latter underlies a Hopf algebra structure whose we consider the
associative product 
$$
.\,\,:\CU(\fL)\otimes\CU(\fL)\to\CU(\fL):X\otimes Y\to X.Y
$$
and co-product
$$
\Delta:\CU(\fL)\to\CU(\fL)\otimes \CU(\fL): X\mapsto \Delta(X)\;.
$$
These structures naturally $\C[[\nu]]$-linearly extend to the space of formal power series 
$\CU(\fL)[[\nu]]$, endowing it with a bi-algebra structure.
\begin{dfn}
A \emph{Drinfel'd twist based on $\fL$} is an element $F\in\CU(\fL)\otimes\CU(\fL)[[\nu]]$ of the form
\begin{equation}\label{DFT}
F\;=\;1\otimes1\,+\,\sum_{k=1}^\infty\nu^kF_k
\end{equation}
such that the following \emph{cocycle relation} holds:
\begin{equation}\label{DFCOCYCLE}
(\Delta\otimes\id)(F).(F\otimes1)\;=\;(\id\otimes\Delta)(F).(1\otimes F)\;.
\end{equation}
\end{dfn}
In the above relation (\ref{DFCOCYCLE}), the $``\,.\,"$ stands for the algebra
structure on $\CU(\fL)\otimes\CU(\fL)\otimes\CU(\fL)[[\nu]]$ naturally induced by 
the product on $\CU(\fL)$.

\noindent Another way to formulate the notion of Drinfel'd twist is the following one.
Exactly along the same lines that a left-differential operator on $\L$ defines an element of $\CU(\fL)$, every element $F_k\in\CU(\fL)\otimes\CU(\fL)$ occurring in a series such as 
(\ref{DFT}) may be interpreted as the value at the unit $e$ of $\L$ of a uniquely associated left-invariant
\emph{bi-differential} operator $$\widetilde{F}_k:C^\infty(\L)\times C^\infty(\L)\to C^\infty(\L)\;.$$
To the series $F$ is therefore uniquely associated a formal power series of left-invariant differential
operators:
\begin{equation}\label{DFTINV}
\widetilde{F}\;=\;m_0\,+\,\sum_{k=1}^\infty\nu^k\widetilde{F}_k
\end{equation}
where $m_0$ denotes the pointwise multiplication of functions on $C^\infty(\L)$.
One then has that the cocycle condition (\ref{DFCOCYCLE}) exactly corresponds to associativity, namely one has
\begin{prop}\label{DTSP}[Drinfel'd]
An element $F\in\CU(\fL)\otimes\CU(\fL)[[\nu]]$ is a Drinfel'd twist if and only if the associated element 
$\widetilde{F}$ (cf. (\ref{DFTINV})) consists in the cochain series of a left-invariant formal star-product $\tilde{\star}^F_\nu$ on $\L$ i.e. 
an left-invariant associative algebra law on $C^\infty(\L)[[\nu]]$ that satisfies conditions (1) and (3) of Definition \ref{STARPROD}.
\end{prop}
Given such a star-product, one readily checks that the associativity at order one in the deformation parameter $\nu$
corresponds to the property that the skewsymmetrization of the first order cochains is a \emph{Poisson bracket}. Namely, 
one has in our case
\begin{lem}
The element 
$$
\{\,,\,\}^F:C^\infty(\L)\times C^\infty(\L)\to C^\infty(\L):(u,v)\mapsto\widetilde{F}_1(u,v)-\widetilde{F}_1(v,u)
$$
defines a structure of Lie algebra on $C^\infty(\L)$ with the additional property that for every 
$u\in C^\infty(\L)$, the element 
$$
\Xi^F_u\;:=\;\{\,u\,,\,.\,\}^F
$$
is a vector field on $\L$.
\end{lem}
From left-invariance and Jacobi, one easily proves
\begin{prop}
The tangent distribution $\tilde{\s}^F$ on $\L$ defined by
$$
x\mapsto\tilde{\s}^F_x\;:=\;\{\Xi^F_u|_x\;:\,u\in C^\infty(\L)\}
$$
is left-invariant i.e. preserved by the left-translations. In particular, it has constant rank $r$.

\noindent Moreover, it is smooth and \emph{involutive} in the sense that the bracket of vector fields is internal. 
\end{prop}
The distribution $\tilde{\s}^F$ is therefore integrable in the sense that the group manifold $\L$ foliates into a collection
of $r$-dimensional immersed leaves all tangent to $\tilde{\s}^F$. 
From all this, one then immediately gets
\begin{prop}
(i) The leaf $\S^F$ tangent to $\tilde{\s}^F$ through the unit $e$ consists in an immersed Lie sub-group of $\L$.

\vspace{2mm}

(ii) The Poisson bi-vector field 
$$
\tilde{w}^F(u,v)\;:=\;\{u,v\}^F
$$
when restricted to $\S^F$ is \emph{symplectic} in the sense that 
for every $x\in \S^F$, the map
\begin{equation}\label{MUSIC}
T^\star_x(\S^F)\to T_x(\S^F):\alpha\mapsto \iota_\alpha\tilde{w}^F_x\;=:\;{}^\sharp\alpha
\end{equation}
is a linear isomorphism.
\end{prop}
In the above statement, the symbol $\iota$ denotes the interior product 
of bi-vectors by one-forms.

\begin{rmk} The inverse of the above isomorphism (\ref{MUSIC})
:
$$
T_x(\S^F)\to T^\star_x(\S^F):v\mapsto{}^\flat v
$$
corrseponds to an element $\omega^F_x$ of $T^\star_x(\S^F)\otimes T^\star_x(\S^F)$:
$$
\omega^F_x(v,w)\;:=\;<{}^\flat v\,,\,w>\;.
$$
The field 
$$
\omega^F:x\mapsto\omega^F_x
$$
turns out to be a symplectic 2-form on $\S^F$ which is invariant under the left-translations.
\end{rmk}
The above discussion motivates the following definition (see \cite{Li2}) which singles out
the class of Lie groups that are semi-classical limits of Drinfel'd twists.
\begin{dfn}
A \emph{symplectic Lie group} is a pair $(\S,\omega)$  where $\S$ is a Lie group and $\omega$ is a symplectic structure on $\S$ that is invariant under the left-translations:
$$
L^\star_x\omega\;=\;\omega\quad (\forall x \in \S)\;.
$$
\end{dfn}
One may wonder whether any symplectic Lie group comes from a Drinfel'd twist.
This is indeed the case.
\begin{prop}
Every symplectic Lie group admits a natural left-invariant
star-product.
\end{prop}
\Pf
We start by considering the canonical affine symmetric connection $\nabla^\S$ on $\S$
associated with the structure of symmetric space of the Lie group $\S$.
This connection is not compatible with $\omega$ in general in the sense that
only when $\S$ is Abelian, the 2-form $\omega$ is parallel \cite{Bithese}. 
But the usual trick (due to Lichn\'erowicz) produces a symplectic connection $\nabla$ from the 
data of $\nabla^\S$:
$$
\nabla_XY\;:=\;\nabla_XY\,+\,\frac{1}{3}{}^\sharp\left(\iota_X\nabla^\S_Y\omega
\,+\,\iota_Y\nabla^\S_X\omega
\right)\;.
$$
This connection $\nabla$ is clearly left-invariant. Now, we use Fedosov's method to uniquely associate
a star-product $\star_\nu$ to the data of $(\omega,\nabla)$ \cite{F}. The constructed
star-product then answers the question.
\EPf
As a consequence of the above proof, or more precisely, of the fact that every symplectic Lie group
admits a left-invariant symplectic connection, we are in the situation where we can apply
Theorem \ref{EXISTSTAR} to this specific situation:
\begin{prop}
On a simply connected symplectic Lie group $\S$ with Lie algebra $\s$, the left-equivariant equivalence classes of left-invariant (symplectic) star-products are parametrized by the sequences of elements
of the second Chevalley cohomology space $H^2_{Chev}(\s)$ associated to  the trivial representation of $\s$ on $\R$.
\end{prop}
\section{Intertwining Drinfel'd twists -- the Retract method}
\subsection{Symplectic twists and naturality}
We consider a simply connected \emph{solvable}\footnote{This hypothesis of solvability is in fact not essential. However, it brings proofs that are more Lie theoretic.} symplectic Lie group $(\S,\omega)$ with Lie algebra $\s$. The left-invariant symplectic structure
$\omega$ determines a bilinear two-form 
$$
\Omega\;:=\;\omega_e
$$
on $\s$.
\begin{dfn}
A Drinfel'd twist $F$ based on $\s$ is said to be \emph{based on} $(\s,\Omega)$ 
if its associated left-invariant Poisson structure on $\S$ coincides with the Poisson 
bracket associated to $\omega$.
\end{dfn}
\begin{dfn} Let $F^{(j)}$ ($j=1,2$) be two Drinfel'd twists based on $(\s,\Omega)$. 
We say that they are \emph{equivalent} when their associated left-invariant formal star-products
$\tilde{\star}_\nu^{F^{(j)}}$ live in the same $\S$-left-equivariant equivalence class.
\end{dfn}
In the present section, we give a description of the equivalence classes of Drinfel'd twists.
\begin{dfn}
Let $F$ be a Drinfel'd twist based on $(\s,\Omega)$ with associated left-invariant star-product
$\tilde{\star}_\nu^{F}$. The \emph{internal symmetry Lie algebra} $\g_F$ of $F$ is 
the Lie algebra of vector fields on $\S$ that derivate $\tilde{\star}_\nu^{F}$:
$$
\g_F\;:=\;\Gamma^\infty(T(\S))\cap\Der(\tilde{\star}_\nu^{F})\;.
$$
\end{dfn}
\begin{rmk}
One has the canonical inclusion
$$
\s\hookrightarrow\g_F:X\mapsto {X}^\star
$$
with
$$
X^\star_x\;:=\;\Dto\exp(-tX)x\;.
$$

\end{rmk}
\begin{lem}
When a twist $F$ based on $(\s,\Omega)$ is natural, in the sense that $\tilde{\star}_\nu^{F}$ is natural, then
its internal symmetry Lie algebra is finite dimensional.
\end{lem}
\Pf Such a natural star-product determines a an affine connection whose $\g_F$ is a Lie sub-algebra 
of derivations \cite{RG}. \EPf
\subsection{The retract method I: describing equivalences as evolution solutions}
\noindent We now consider, within a given equivalence class of Drinfeld's twists based on $(\s,\Omega)$,
a natural element $F$ (from what we know so far, such an element always exists). We denote by
$$\sharp\;:=\;\tilde{\star}_\nu^{F}$$ the associated left-invariant star-product on $\S$ and by $$\g\;:=\;\g_F$$
its internal symmetry Lie algebra. 

\noindent We also fix another Drinfel'd twist equivalent to $F$ with associated left-invariant star-product
denoted by $\star$,  and no further requirement.
\begin{prop}\label{HOMS} {\rm [Relevance of the retract method]}
Let $\homo_\s(\g,\Der(\star))$ the space of Lie algebra homomorphisms from $\g$ to the Lie algebra $\Der(\star)$ of derivations of $(C^\infty(\S)[[\nu]],\star)$ that are \emph{$\s$-relative}
in the sense that every element $D$ of $\homo_\s(\g,\Der(\star))$ restricts to $\s$ as the identity:
$$
D_X\;=\;{X}^\star\;.
$$
Then, $\homo_\s(\g,\Der(\star))$ is finite dimensional over the formal field $\R[[\nu]]$.
\end{prop}
Before passing to the proof, we observe
\begin{lem}
Let $V$ be a finite dimensional real vector space. Consider $f\in C^\infty(\S)\otimes V$
and $A\in\mbox{\rm Hom}_\R(\s,\End(V))$. Then, the space of solutions of the equation
$$
X^\star\ulambda\;=\;A(X)\ulambda\,+\,f
$$
is finite dimensional.
\end{lem}
\Pf
Consider an element $0\neq X_1\in\s\backslash[\s,\s]$ and a vector sub-space
$\s_2\;<\;\s$ such that
$$
\s\;=\;\R X_1\oplus\s_2\quad\mbox{\rm and}\quad [\s,\s]\subset \s_2\;.
$$
The space $\s_2$ is then automatically a Lie sub-algebra of $\s$.

\noindent Denote by $\S_2$ the analytic sub-group of $\S$ with Lie algebra $\s_2$ and
consider the local coordinate system
$$
\R\times\S_2\to\S:(x^1,y)\mapsto\exp(-x^1X_1).y\;.
$$
Within this coordinate system, the equation 
$$
X_1^\star\ulambda\;=\;A(X_1)\ulambda\,+\,f
$$
reads
$$
\partial_{x^1}\ulambda\;=\;A_1\ulambda+f\quad(A_1\;:=\;A(X_1))\;,
$$
whose every solution is of the form
$$
\ulambda(x^1,y)\;=\;\exp(x^1A_1)\ulambda_2(y)\,+\,\mu_1(x^1,y)
$$
with
$$
\mu_1(x^1,y)\;:=\;\int_{x_0^1}^{x^1}\exp(-tA_1)f(\exp(tX_1),y)\,{\rm d}t
$$
and 
$\ulambda_2\in C^\infty(\S_2)\otimes V$ arbitrary.

\noindent In particular, in the case $m=1$ the equation admits an affine
space of solutions of dimension $\dim V$.

\noindent We now proceed by induction on $\dim\S\;=:m$.

\noindent The Lie sub-algebra $\s_2$ being solvable as well, one fixes an 
element $0\neq X_2\in\s_2\backslash[\s_2,\s_2]$ and a vector sub-space
$\s_3$ of $\s_2$ such that
$$
\s_2\;=\;\R X_2\oplus\s_3\quad\mbox{\rm and}\quad [\s_2,\s_2]\subset \s_3\;.
$$
\noindent Denote by $\S_3$ the analytic sub-group of $\S_2$ with Lie algebra $\s_3$ and 
consider the local coordinate system
$$
\R^2\times\S_2\to\S:(x^1,x^2,z)\mapsto\exp(-x^1X_1).\exp(-x^2\Ad_{x_1^{-1}}X_2).z\;,
$$
where 
$$
x_1\;:=\;\exp(-x^1X_1)\;.$$
Within this coordinate system, setting
$$
x_2\;:=\;\exp(-x^2\Ad_{x_1^{-1}}X_2)\;,
$$
one has
$$
X^\star_2\;=\;\partial_{x^2}\;.
$$
Indeed:
\begin{eqnarray*}
&&
\ddto\exp(-tX_2)x_1x_2z\;=\;\ddto x_1x_1^{-1}\exp(-tX_2)x_1x_2z\\
&=&
L_{x_1\star}\ddto\exp(-t\Ad_{x_1^{-1}}X_2)\exp(-x^2\Ad_{x_1^{-1}}X_2)z\\
&=&L_{x_1\star}\ddto\exp(-(t+x^2)\Ad_{x_1^{-1}}X_2)z\;.
\end{eqnarray*}
Therefore:
$$
X^\star_2\ulambda(x^1,x^2,z)\;=\;\exp(x^1A_1)\partial_{x^2}\ulambda_2(x^2,z)+\mu_2(x^1,x^2,y)
$$
with
$$
\mu_2\;:=\;\partial_{x^2}\mu_1\;.
$$
This leads, with $A_2\;:=\;A(X_2)$, to
$$
A_2\ulambda+f\;=\;A_2\exp(x^1A_1)\ulambda_2+A_2\mu_2+f\;.
$$
Hence, the equation with $X\;:=\;X_2$ reads:
$$
\partial_{x^2}\ulambda_2\;=\;A'_2\ulambda_2+f_2
$$
with
$$
A'_2\;:=\;\exp(-x^1A_1).A_2.\exp(x^1A_1)\quad\mbox{\rm and}\quad f_2\;:=\;\exp(-x_1A_1)(A_2\mu_2+f)\;.
$$
The induction then yields the assertion.
\EPf
{\sl Proof of the Proposition.}
Consider a vector sub-space $W$ of $\g$ supplementary to $\s$:
$$
\g\;=\;\s\oplus W
$$
which is assumed to be non-trivial. Let $\{Y_a\}_{a\in\{1,...,n\}}$ be a basis of $W$. Furthermore, we will consider a basis $\{X_j\}_{j\in\{1,...,m\}}$
of $\s$.

\noindent Denote by $C$ and $\uC$ the structure constants associated to the action of $\s$:
$$
[X_j,Y_a]\;=:\;C^k_{ja}X_k\,+\,\uC^b_{ja}Y_b\;.
$$
Let $D\in\homo_\s(\g,\Der(\star))$.

\noindent Since $\S$ is simply connected, we have functions $\{\lambda_j\}$ and $\{\ulambda_a\}$ such that
$$
\nu X_j^\star u\;=\;[\lambda_j\,,\,u]_\star\quad\mbox{\rm and}\quad \nu D_{Y_a}u\;=\;[\ulambda_a\,,\,u]_\star\quad(u\in C^\infty(\S))\;.
$$
Moreover, the Lie algebra homomorphism condition on $D$ implies that  for all $X_j$ and $Y_a$, the 
element 
$$
c_{ja}\;:=\;C^k_{ja}\lambda_k\,+\,\uC^b_{ja}\ulambda_b-[\lambda_j\,,\,\ulambda_a]_\star
$$
is central in $(C^\infty(\S)[[\nu]],\star)$ hence a formal constant (because symplectic).

\noindent Therefore for every $j$, we have two matrices $C_j$ and $\uC_j$:
$$
[C_j]^k_a\;:=\;-\nu^{-1}C_{ja}^k\quad\mbox{\rm and}\quad-[\uC_j]^b_a\;:=\;\nu^{-1}\uC_{ja}^b
$$
such that
\begin{equation}\label{COMM}
X_j^\star\ulambda\;=\;-C_j\lambda\,+\,-\uC_j\ulambda\,+\,c_{j}
\end{equation}
where $c_j$ is a vector in $\R^n[[\nu,\nu^{-1}]$.

\noindent The assertion is the an immediate consequence of the Lemma.
\EPf

\noindent Now, let us discuss the intertwiners between $\star$ and $\sharp$. We start by observing
\begin{lem}
Because it commutes with the left-regular action, any (formally invertible) intertwiner $T$ between $\star$ and $\sharp$ (i.e. $T(\star)=\sharp$) must be a convolution operator of the form
$$
T^{-1}\varphi(x)\:=\;<u\,,\,L^\star_x\varphi>
$$
where $u\in\CD'(\S)[[\nu]]$ denotes the (formal) distributional kernel of $T^{-1}$.
\end{lem}
\Pf By the kernel theorem of Schwartz, every differential operator $A$ on $C^\infty(\S)$
restricts to $\CD(\S)$ under the form
$$
A\varphi\;=:\;[x\mapsto<u_A(x),\varphi>]
$$
where $u_A\in C^\infty(\S,\CD'(\S))$. When left-invariant, one obtains
$$
A\varphi(x)\;=\;<u\,,\,L_x^\star\varphi>\quad\mbox{\rm with}\quad u\;:=\;u_A(e)\;.
$$
This applies of course to the formal power series case. \EPf

\noindent A similar statement
of course holds for $T$ but we will use it only later on. 

\noindent Back to $T^{-1}$, we will adopt the 
following notation in accordance with the distributions that are defined by locally summable 
functions:
 $$
 T^{-1}\varphi(x)\;:=:\;<L^\star_{x^{-1}}u\,,\,\varphi>\;=\;\int_\S u(x^{-1}y)\,\varphi(y)\,{\rm d}y
 $$
 where ${\rm d} y$ denotes a left-invariant Haar measure on $\S$.
 \begin{lem}\label{SMOOTHDISTR} Let $u\in\CD'(\S)$ and $A$ be a differential operator on $C^\infty(\S)$.
 Then, for every $\varphi\in\CD(\S)$, the map $$x\mapsto <L^\star_{x^{-1}}u,\varphi>$$
 belongs to $C^\infty(\S)$, and one has
 $$
 <A_x[L^\star_{x^{-1}}u]\,,\,\varphi>\;=\;A_x(<[L^\star_{x^{-1}}u]\,,\,\varphi>)\;.
 $$
 \end{lem}
 \Pf Set
 $$
\tau\varphi(x)\;:=\;<[L^\star_{x^{-1}}u]\,,\,\varphi>\;=\;<u,L_x^{\star}\varphi>\;.
 $$
Consider the case where $A$ is a vector field with local flow $\phi$:
$$
A_x\;=:\;\ddto\phi_t(x)\;.
$$
One has
$$
\frac{1}{t}\left(<u\,,\,L^\star_{{\phi_t(x)}}\varphi>-
<u\,,\,\varphi>\right)\;=\;<u\,,\,\frac{1}{t}\left(L^\star_{{\phi_t(x)}}\varphi-
\varphi\right)\,>\;.
$$
The limit $t\to0$ of $\frac{1}{t}\left(L^\star_{{\phi_t(x)}}\varphi-
\varphi\right)$ exists in $\CD(\S)$. Hence the limit $t\to0$ of the right hand side to the above equality exists since $u\in\CD'(\S)$. Which proves, that $[x\mapsto\tau\varphi(x)]$ is smooth.

\noindent Iterating one gets the assertion for higher order differential operators.
\EPf
\begin{prop}\label{RETRACT1}
\noindent Consider such an intertwiner $T$  with distributional kernel $u$. Then.

\vspace{3mm} 

\noindent (i) the map $$D^T:\g\to\Der(\star):X\mapsto T^{-1}X^\star T$$ belongs to $\homo_\s(\g,\Der(\star))$.

\vspace{3mm} 

\noindent (ii) The kernel $u$ is a (weak) joint solution of the following evolution equations:
$$
-X^\star_y[u(x^{-1}y)]\;=\;D^T_X|_x[u(x^{-1}y)]\quad(\forall X\in \g)\;.
$$

\end{prop}
\Pf With the notations introduced above and based on the previous lemmas, we have
\begin{eqnarray*}
&& 
T^{-1}X^\star(\varphi)(x)\;=\;\int u(x^{-1}y)\,X^\star_y(\varphi)\,{\rm d}y
\;=\;-\int X^\star_y[u(x^{-1}y)]\,\varphi(y)\,{\rm d}y\\&=&D^T_X|_x\int  u(x^{-1}y)\,\varphi(y)\,{\rm d}y
\;=\;\int D^T_X|_x u(x^{-1}y)\,\varphi(y)\,{\rm d}y
\end{eqnarray*}
because $X^\star$ is symplectic and thus preserves the (Liouville) Haar measure and because $\varphi$ is compactly supported.
\EPf
\subsection{The retract method II: the general context}
In this section, we consider 
\begin{enumerate}
\item[$\bullet$]
a solvable, simply connected symplectic Lie group 
$(\S,\omega)$ with symplectic Lie algebra $(\s, \Omega)$.
\item[$\bullet$] A fixed finite dimensional Lie algebra $\g$ that contains $\s$
as Lie sub-algebra.
\item[$\bullet$] A fixed Drinfel'd twist based on $(\s, \Omega)$ or, equivalently,
a left-invariant formal star-product $\star$ on $(\S,\omega)$.
\end{enumerate}
\begin{lem}
Let $D$ be an element in $\homo_\s(\g,\Der(\star))$. 
\begin{enumerate}
\item[(i)] For every $X\in\g$, the map
$$
X_D^\star:C^\infty(\S)\to C^\infty(\S);f\mapsto D_X(f)\mod(\nu)\;=:\;X_D^\star f
$$
is a symplectic vector field.
\item[(ii)] The map 
\begin{equation}\label{ZEROACT}
\g\to\Gamma^\infty(T(\S)):X\mapsto X_D^\star
\end{equation}
is a homomorphism of Lie algebras.
\item[(iii)] For every $x\in\S$, we set
$$
\k_x\;:=\;\{X\in\g\;|\;\left.X^\star_D\right|_x\;=\;0\}
$$
and, at the unit element $e\in\S$:
$$
\k\;:=\;\k_e\;.
$$
Then, at every $x\in\S$, the stabilizer $\k_x$ is a Lie sub-algebra
of $\g$ and one has
$$
\k_x\;=\;\Ad_x\k\;.
$$
\end{enumerate}
\end{lem}
\Pf
As previously, we consider the map
$$
\ulambda:\g\to C^\infty(\S)[[\nu]];X\mapsto\ulambda_X
$$
determined (up to additive formal constants) by the conditions:
$$
D_X\varphi\;=\;\frac{1}{2\nu}\left[\ulambda_X,\varphi\right]_\star\;.
$$
Taking zero orders then yields items (i) and (ii).

\noindent Item (iii) is obvious. \EPf
\begin{dfn}
Let $(M,{\rm d}\mu)$ be a smooth manifold endowed with a volume form and denote by 
$L^2(M,{\rm d}\mu)$ the associated Hilbert space of square integrable function classes. Let 
$A$ be a differential operator in $C^\infty(M)$. Denote by $A^\dagger$ its adjoint on $\CD(M)$
with respect to the Hilbert structure:
$$
\int_M\overline{A^\dagger\varphi}\,\psi\,{\rm d}\mu\;:=\;\int_M\overline{\varphi}\,A\psi\,{\rm d}\mu\quad(\varphi\,,\,\psi\in\CD(M))\;.
$$
We say that a distribution $u\in\CD'(M)$ belongs to the \emph{weak kernel} of $A$ when
it satisfies the condition
$$
<u,A^\dagger\varphi>\;=\;0
$$
for every $\varphi\in\CD(M)$.
\end{dfn}
\begin{dfn}
Let $D$ be an element in $\homo_\s(\g,\Der(\star))$. 

\noindent An element $u\in\CD'(\S)[[\nu]]$
such that
\begin{enumerate}
\item[(a)] it is a joint (weak) solution of $$D_Zu\;=\;0\quad\forall Z\in\k\;.$$
\item[(b)] Its singular support is reduced to the unit element:
$$
\mbox{\rm supp}(u)\;=\;\{e\}\;.
$$
\item[(c)] It satisfies the following semi-classical condition:
$$u\;\mbox{\rm mod}\;\nu\;=\;\delta_e\;.$$
\end{enumerate}
is called a \emph{$D$-retract}.
\end{dfn}
The following statement is classical.
\begin{prop}\label{UGSUPP}
The space of left-invariant differential operators on $\S$ (canonically isomorphic to $\CU(\s)$)
identifies with the sub-space of distributions in $\CD'(\S)$ supported at the unit $e$.
\end{prop}
\Pf
We first recall Peetre's theorem that realizes the space of differential operators on a manifold $M$ as the space of endomorphisms
of $C^\infty(M)$ that reduce the supports \cite{P}.

\noindent Let $A$ be an element in $\CU(\s)$ and $\tilde{A}$ the corresponding left-invariant differential
operator on $C^\infty(\S)$. Denote by $u_A$ its distributional kernel:
$$
\tilde{A}_x\varphi\;=:\;<u_A\,,\,L_x^\star\varphi>\;.
$$
Assume that $\mbox{\rm supp}\varphi$ does not contain the unit $e$. Then, since $\tilde{A}$ is differential,
one has the inclusion $\mbox{\rm supp}\tilde{A}\varphi\subset\mbox{\rm supp}\varphi$. Hence 
$$
\tilde{A}_e\varphi\;=\;0\;=\;<u_A\,,\,\varphi>\;.
$$
In other words, the support of $u_A$ consists in the unit. 

\noindent Reciprocally, assume $u$ to be supported on the unit $e$ only and consider the left-equivariant
mapping $\tilde{U}:\CD(\S)\to C^\infty(\S)$ defined by
$$
\tilde{U}\varphi(x)\;:=\;<u,L_x^\star\varphi>\;.
$$
Observing, for every $x\in\S$, the equality
$$
\mbox{\rm supp}L_x^\star\varphi\;=\;L_{x^{-1}}\mbox{\rm supp}\varphi\;,
$$
one deduces that, for every point $x$ belonging to $\mbox{\rm supp}\tilde{U}\varphi$:
$$
\tilde{U}\varphi(x)\;\neq0\;\Rightarrow\;\mbox{\rm supp}L_x^\star\varphi\ni e
\;\Rightarrow\;L_{x^{-1}}\mbox{\rm supp}\varphi\ni e
\;\Rightarrow\;\mbox{\rm supp}\varphi\ni x\;,
$$
hence the inclusion
$$
\mbox{\rm supp}\tilde{U}\varphi\subset\mbox{\rm supp}\varphi
$$
which implies differentiability by virtue of Peetre's theorem. \EPf
\begin{cor} Let $u$ be a $D$-retract, and consider the modular function $\Delta\in C^\infty(\S)$
defined by
$$
{\rm d}(y^{-1})\;=:\;\Delta(y)\,{\rm d}y\;.
$$
Then, the convolution operator
$$
\tau_u\varphi(x)\;:=\;<u^\vee\,,\,L_x^\star\varphi>\;:=\;<u\,,\,\Delta R_{x^{-1}}^\star(\varphi^\vee)> \quad(\varphi\in\CD(\S))
$$
is differential and satisfies the semi-classical
condition:
$$
\tau_u\;=\;\Id\;\mbox{\rm mod}\;\nu\;.
$$
It therefore extends to an invertible endomorphism
$$
\tau_u:C^\infty(\S)[[\nu]]\to C^\infty(\S)[[\nu]]
$$
called \emph{inverse retract operator}. Its formal inverse
$$
T_u\;:=\;
\tau_u^{-1}
$$
is called \emph{direct retract operator}.
\end{cor}
\Pf
Within the same setting as in the proof of Proposition \ref{UGSUPP},
we observe that
$$
<u_A\,,\,\Delta R_{x^{-1}}^\star(\varphi^\vee)>\;=\;\tilde{A}_e\left(\Delta R_{x^{-1}}^\star(\varphi^\vee)\right)\;=:\;\tilde{A}^\dagger_x\varphi
$$
which is a left-invariant differential operator as well. What remains is obvious. \EPf

\begin{thm}

\noindent (i) Let $u$ be a $D$-retract with associated direct retract operator $T_u:C^\infty(\S)[[\nu]]\to C^\infty(\S)[[\nu]]$. Then,
$$
T_u(\star)\;=:\;\sharp^u
$$
is a formal star-product on $(\S,\omega)$ that is $\g$-invariant with respect to the
infinitesimal action (\ref{ZEROACT}).

\vspace{3mm}

\noindent (ii) Given any symplectic action of $\g$ on $(\S,\omega)$ that
restricts to $\s$ as the regular one, every $\g$-invariant star-product that is $\S$-equivariantly equivalent to $\star$ is of the form 
$\sharp^u$ for a certain $D$-retract $u$.
\end{thm}
\Pf The $\g$-invariance is implied by the following condition on $u$:
$$
D_X\tau_u\varphi\;=\;\tau_uX^\star_D\varphi
$$
which, setting $v\;:=\;u^\vee$, amounts to
$$
<v\,,\,D_X|_x[L_x^\star\varphi]-L_x^\star X^\star_D\varphi>\;=\;0\;.
$$
We now determine, for every $x$, the adjoint of the operator:
$$
A_{X,x}:\varphi\mapsto D_X|_x[L_x^\star\varphi]-L_x^\star X^\star_D\varphi\;.
$$
Let $w\in\CD(\S)$ and observe that left-invariance of the Haar measure and the fact
 that $X^\star_D$ is a symplectic vector field yield:
\begin{eqnarray*}&&
\int w(y)\,\left(D_X|_x[L_x^\star\varphi](y)-L_x^\star X^\star_D\varphi(y)\right)\,{\rm d}y\\
&=&
\int\left(D_X|_x[w(x^{-1}y)]+X_D^\star|_y[w(x^{-1}y)]\right)\,\varphi(y)\,{\rm d}y\;.
\end{eqnarray*}
Writing every $X$ in $\g$ as
$$
X\;=:\;X_\k+X_\s
$$
according to the decomposition 
$$
\g\;=\;\k\oplus\s\;,
$$
we observe that
$$
X^\star_D|_e\;=\;(X_\s)^\star_e\;.
$$
Now, the $\S$-equivariance of $D$ yields:
\begin{eqnarray*}
&&
D_X|_x[w(x^{-1}y)]+X_D^\star|_y[w(x^{-1}y)]\;=\;D_X|_x[w^\vee(y^{-1}x)]+X_D^\star|_y[w(x^{-1}y)]\\
&=&L^\star_{y^{-1}}\left(D_{\Ad_{y^{-1}}X}w^\vee\right)(x)+\left(\Ad_{y^{-1}}X_D\right)^\star|_{e}(L^\star_{x^{-1}y}w)
\\
&=&L^\star_{y^{-1}}\left(D_{\Ad_{y^{-1}}X}w^\vee\right)(x)+\left(\Ad_{y^{-1}}X\right)_\s^\star|_{e}(L^\star_{x^{-1}y}w)
\\
&=&L^\star_{y^{-1}}\left(D_{\Ad_{y^{-1}}X}w^\vee\right)(x)+\left(\Ad_{y^{-1}}X\right)_\s^\star|_{y^{-1}x}(w^\vee)
\\
&=&L^\star_{y^{-1}}\left(D_{\Ad_{y^{-1}}X}w^\vee\right)(x)-D_{\left(\Ad_{y^{-1}}X\right)_\s} w^\vee(y^{-1}x)
\\
&=&L^\star_{y^{-1}}\left(D_{\left(\Ad_{y^{-1}}X\right)_\k}w^\vee\right)(x)\;=\;
\left(D_{\left(\Ad_{y^{-1}}X\right)_\k}w^\vee\right)^\vee(x^{-1}y)\;.
\end{eqnarray*}
Hence, considering a basis $\{Z_j\}$ of $\k$ and denoting by $[A(y)]$ the matrix
of the linear map $\k\to\k:X\mapsto\left(\Ad_{y^{-1}}X\right)_\k$, we 
get for every $j$:
\begin{eqnarray*}
&&
D_{Z_j}|_x[w(x^{-1}y)]+{Z_j}_D^\star|_y[w(x^{-1}y)]
\;=\;
[A(y)]_j^k\left(D_{Z_k}w^\vee\right)^\vee(x^{-1}y)\;.
\end{eqnarray*}
In other words, we have
$$
<w,A_{Z_j,x}\varphi>\;=\;\sum_k<\left(D_{Z_k}w^\vee\right)^\vee\,,\,L^\star_x([A]_j^k\varphi)>\;.
$$
We conclude that, whenever $u=v^\vee$ is in the weak kernel of $D_{Z_k}$ for all $k$, the star product
$\sharp^u$ is $\g$-invariant, proving item (i).

\noindent Item (ii) follows from Proposition \ref{RETRACT1} (i) combined with the above discussion.
\EPf

\chapter{Admissible phases, Weinstein's and \v{S}evera's areas}

The question addressed in this section is~: {\sl given an orientable symmetric space $M$ endowed with 
a symmetry invariant volume form $\mu$, what are the  conditions on a (symmetry invariant) three-point function 
$$S\in C^{\infty}(M\times M\times M,\R)$$ which would guarantee associativity 
of the product 
$$
u\star v(x)=\int_{M\times M}e^{iS(x\, , . \, , \, . \,)}u\otimes v \,
\mu\otimes \mu\quad ?
$$}
Above, by \emph{symmetry invariant}, we mean:
$$
S(s_wx,s_wy,s_wz)\;=\;S(x,y,z)\quad (\forall x,y,z,w\in M)\;.
$$
When writing a (continuous) multiplication on functions via a kernel formula 
of the type~:  
$$
u\star v(x)=\int_{M\times M}K(x,\, .\, , \,. \,)u\otimes v\quad , 
$$
a computation shows that associativity for the multiplication $\star$ is 
(at least formally) equivalent to the following condition~: 

\begin{equation}\label{ASSK}
\int_{M}K(a,b,t)K(t,c,d)\mu(t)=\int_{M}K(a,\tau,d)K(\tau,b,c)\mu(\tau),
\end{equation}
for every quadruple of points $a,b,c,d$ in $M$. Equality~(\ref{ASSK}) obviously 
holds if one can pass from one integrand to the other using a change a 
variables $\tau=\varphi(t)$. This motivates 

\begin{dfn}\label{coucou}
Let $(M,\mu)$ be a symmetric space endowed with a symmetry invariant volume form $\mu$. 
A three-point kernel $K\in C^{\infty}(M\times M\times M)$ is called {\bf geometrically 
associative} if for every quadruple of points $a,b,c,d$ in $M$ there exists 
a volume preserving diffeomorphism
$$
\varphi:(M,\mu)\to(M,\mu),
$$
such that for all $t$ in $M$~:
$$
K(a,b,t)K(t,c,d)=K(a,\varphi(t),d)K(\varphi(t),b,c).
$$
\end{dfn}

We will prove in Proposition~\ref{GEOM} that the following structure leads (under an additional condition) to a 
geometrically associative kernel.
\begin{dfn}\label{WT} Let $(M,\mu)$ be an oriented symmetric space.
An {\bf admissible function} $S:M\times M\times M \to {\R}$ is a smooth symmetry-invariant  three-point function 
on $M$ such that~:
\begin{enumerate}
\item[(i)] $S(x,y,z)=S(z,x,y)=-S(y,x,z)$;
\item[(ii)] for all $x \in M$, one has: 
$$
S(x,y,z)=-S(x,s_x(y),z)\qquad \forall y,z\in M.
$$
\end{enumerate}
\end{dfn}
Property (i) in Definition~\ref{WT} naturally leads us to adopt the following 
``oriented graph" type notation for~$S$~:

$$\begin{picture}(200,50)
\put(0,0){$x$}
\put(50,0){$z$}
\put(25,40){$y$}
\put(5,5){\trianglehorizontal}
\put(57,19){$\stackrel{\mbox{def.}}{=}\,S(x,y,z).$}
\end{picture}$$
A change of orientation in such an ``oriented triangle'' leads to a change 
of sign of its value. However, the value represented by such a ``triangle''  
does not depend on the way it ``stands'', only the data of the vertices 
and the orientation of the edges matters.

Now,  given two compactly supported functions $u$ and $v\in 
C^{\infty}_{c}(M)$, consider the following ``product''~:  
$$
u\star v(x)\;:=\int_{M\times M} u(y) v(z) e^{iS(x,y,z)}\mu(y) \,\mu(z).
$$

With the above notation for $S$, associativity for $\star$ then formally reads as follows~:

$$\begin{array}{l}
\begin{array}{ccccc}
\begin{picture}(60,50)
\put(0,25){
$u\star(v\star w)(a)=$}
\end{picture}
&
\begin{picture}(60,80)
\put(0,25){\mbox{$\displaystyle \int \exp i ($}}

\put(30,10){
\begin{picture}(50,35)
\put(0,0){$d$}
\put(30,15){$t$}
\put(0,30){$a$}
\put(5,5){$\bullet$}
\put(5,25){$\bullet$}
\put(25,15){$\bullet$}
\put(7,17){\vector(0,1){2}}
\put(17,22){\vector(2,-1){2}}
\put(17,12){\vector(-2,-1){2}}
\put(7,7){\line(0,1){20}}
\put(7,7){\line(2,1){20}}
\put(27,17){\line(-2,1){20}}
\end{picture}
}
\end{picture}
&

\begin{picture}(60,10)
\put(0,25){\mbox{$) \displaystyle \int \exp i ($}}

\put(40,10){
\begin{picture}(50,35)
\put(25,0){$c$}
\put(-2,15){$t$}
\put(25,30){$b$}
\put(20,5){$\bullet$}
\put(20,25){$\bullet$}
\put(2,15){$\bullet$}
\put(23,17){\vector(0,-1){2}}
\put(13,22){\vector(2,1){2}}
\put(13,12){\vector(-2,1){2}}
\put(23,7){\line(0,1){20}}
\put(23,7){\line(-2,1){20}}
\put(3,17){\line(2,1){20}}
\end{picture}
}
\end{picture}
&
\begin{picture}(70,50)
\put(0,25){
$)\;u(b)\,v(c)\,w(d)\;$}
\end{picture}
&
\end{array} \\

\begin{array}{cccc}
\begin{picture}(60,10)
\put(-7,25){\mbox{$=\displaystyle \int \exp i ($}}

\put(43,15){

\begin{picture}(50,35)
\put(22,25){$b$}
\put(10,-2){$\tau$}
\put(-8,25){$a$}
\put(16,20){$\bullet$}
\put(-4,20){$\bullet$}
\put(6,2){$\bullet$}
\put(8,23){\vector(1,0){2}}
\put(3,13){\vector(-1,2){2}}
\put(13,13){\vector(-1,-2){2}}
\put(18,23){\line(-1,0){20}}
\put(18,23){\line(-1,-2){10}}
\put(8,3){\line(-1,2){10}}
\end{picture}
}
\end{picture}
&

\begin{picture}(55,80)
\put(0,25){\mbox{$) \displaystyle \int \exp i ($}}

\put(43,12){

\begin{picture}(60,35)
\put(21,1){$c$}
\put(6,30){$\tau$}
\put(-10,1){$d$}
\put(16,5){$\bullet$}
\put(-4,5){$\bullet$}
\put(6,25){$\bullet$}
\put(8,7){\vector(-1,0){2}}
\put(3,17){\vector(1,2){2}}
\put(13,17){\vector(1,-2){2}}
\put(18,7){\line(-1,0){20}}
\put(18,7){\line(-1,2){10}}
\put(8,27){\line(-1,-2){10}}
\end{picture}
}
\end{picture}
&

\begin{picture}(160,50)
\put(0,25){
$) \;u(b)\,v(c)\,w(d)\;$}
\end{picture}
&
\begin{picture}(40,50)
\put(-95,25){$=(u\star v)\star w(a),$}
\end{picture}
\end{array}
\end{array}$$
the $\mu$-integration being taken over variables $b,c,d,t$ and $\tau$.
This leads, for $K=e^{iS}$, to an equality between two ``distribution valued  
functions'' on $M\times M\times M\times M$ (cf. formula~(\ref{ASSK})): for every 
quadruple of points $a,b,c,d$ in $M$ associativity for $\star$ reads
\begin{equation}\label{ASS}
\begin{array}{ccccc}
\begin{picture}(50,100)
\put(0,50){${\displaystyle \int}_M \exp\;i($}
\end{picture} & \hspace{-6mm}\carrehoriz{$t$} & 
 \hspace{-8mm}
\begin{picture}(80,100)
\put(0,50){$)\;\mu(t)={\displaystyle \int}_M \exp\;i($}
\end{picture} &  \hspace{-4mm}\carrevertic{$\tau$} & 
 \hspace{-6mm}
\begin{picture}(50,100)
\put(0,50){$)\;\mu(\tau).$}
\end{picture}
\end{array}
\end{equation}
In the above formula, the diagram in the argument of the exponential in 
the LHS (respectively the RHS) stands for $S(a,b,t)+S(t,c,d)\;$ (respectively 
$S(a,d,\tau)+S(\tau,b,c)$).  

\begin{prop}\label{GEOM}
Let $S$ be admissible on $M$. Then, provided the following ``cocycle
condition" holds:
$$
S(m,x,y)+S(m,y,z)+S(m,z,x)\;=\;S(x,y,z) \quad(\forall x,y,z,m\in M)\;,
$$
 the associated three-point kernel 
$K=e^{iS}$ is geometrically associative.  
\end{prop}
\Pf
Fix four points $a,b,c,d$. Regarding \dref{coucou} and 
formula~(\ref{ASS}), one 
needs to construct our volume preserving diffeomorphism 
$\varphi:(M,\mu)\to(M,\mu)$ in such a way that for all $t$,
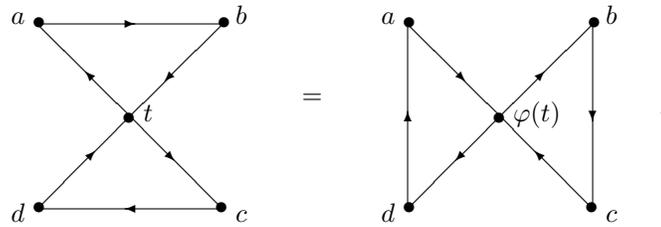
$$
\begin{array}{ccc}
\carrehoriz{$t$} & 
\begin{picture}(20,100)
\put(0,50){=}
\end{picture} & \carrevertic{$\varphi(t)$} 
\end{array}.
$$
We first observe 
that the data of four points $a,b,c,d$ determines what we call
an ``$S$-barycenter", that is a point $g=g(a,b,c,d)$ such that

$$
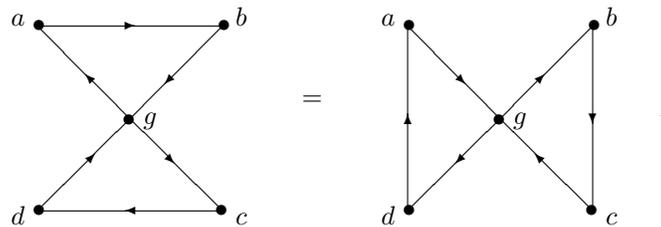
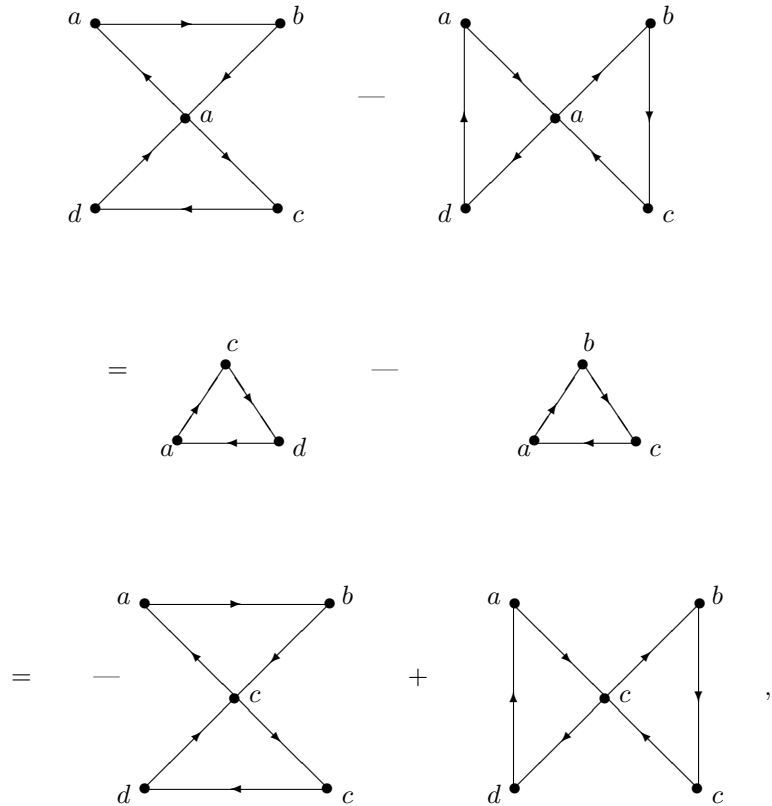
\begin{array}{ccc}
\carrehoriz{$g$} & 
\begin{picture}(20,100)
\put(0,50){=}
\end{picture} & \carrevertic{$g$} 
\end{array}.
$$
Indeed, since
$$
\begin{array}{ccc}
\carrehoriz{$a$} &
\begin{picture}(20,100)
\put(0,50){---}
\end{picture} & \carrevertic{$a$} \\
\begin{picture}(10,100)
\put(0,50){=}
\end{picture}  \Triangle{a}{c}{d} & \begin{picture}(10,100)
\put(0,50){---}
\end{picture} & \Triangle{a}{b}{c}
\end{array}
$$

$$
\begin{array}{ccccc}
\begin{picture}(30,100)
\put(0,50){= \qquad ---}
\end{picture} &
\carrehoriz{$c$} &
\begin{picture}(20,100)
\put(0,50){+}
\end{picture} & \carrevertic{$c$}
\end{array},
$$
any continuous path joining $a$ to $c$ contains such a point $g$.

Now, we fix once for all such an $S$-barycenter $g$ for $\{ a,b,c,d\}$ and 
we adopt the following notation. For all $x$ and $y$ in $M$, the value 
of $S(g,x,y)$ is denoted by a ``thickened arrow''~:
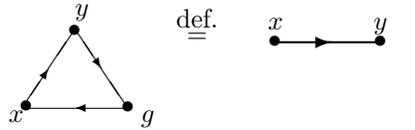
$$
\begin{array}{ccc}
\Triangle{x}{y}{g} & 
\begin{picture}(30,50)
\put(0,50){\mbox{   $\stackrel{\mbox{def.}}{=}$}}
\end{picture} & 
\begin{picture}(38,55)
\thicklines
\put(-2,55){$x$}
\put(38,55){$y$}
\put(-2,48){$\bullet$}
\put(38,48){$\bullet$}
\put(20,50){\vector(1,0){2}}
\put(0,50){\line(1,0){40}}
\end{picture}
\end{array}\quad.
$$
Again, a change of orientation in such an arrow changes the sign of its 
value. Also, property (ii) of admissibility (Definition~\ref{WT}) which 
reads  
$$
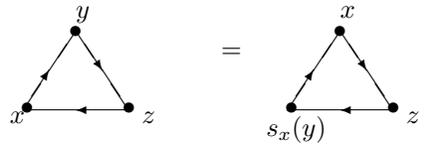
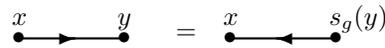
\begin{array}{ccc}
\begin{picture}(60,50)
\put(0,0){$x$}
\put(50,0){$z$}
\put(25,40){$y$}
\put(5,5){\trianglehorizontal}
\end{picture} & 
\begin{picture}(20,50)
\put(10,25){=}
\end{picture} & 
\begin{picture}(60,50)
\put(-3,-5){$s_x(y)$}
\put(50,0){$z$}
\put(25,40){$x$}
\put(5,5){\trianglehorizontal}
\end{picture}
\end{array},
$$
implies
$$
\begin{array}{ccc}
\begin{picture}(50,100)
\thicklines
\put(-2,55){$x$}
\put(38,55){$y$}
\put(-2,48){$\bullet$}
\put(38,48){$\bullet$}
\put(20,50){\vector(1,0){2}}
\put(0,50){\line(1,0){40}}
\end{picture}
&
\begin{picture}(10,90)
\put(0,50){=}
\end{picture} 
&
\begin{picture}(50,100)
\thicklines
\put(-2,55){$x$}
\put(38,55){$s_g(y)$}
\put(-2,48){$\bullet$}
\put(38,48){$\bullet$}
\put(20,50){\vector(-1,0){2}}
\put(0,50){\line(1,0){40}}
\end{picture}
\end{array}
$$
for all $x$ and $y$ in $M$.
While, from the cocycle condition (with $m=g$), one gets
$$
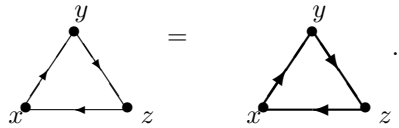
\begin{array}{ccc}
\Triangle{x}{y}{z} &
\begin{picture}(10,90)
\put(0,50){=}
\end{picture} &
\Trianglefat{x}{y}{z}
\end{array}.
$$
Moreover, the barycentric property of $g$ reads
$$
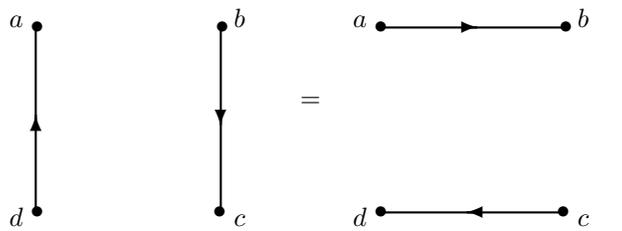
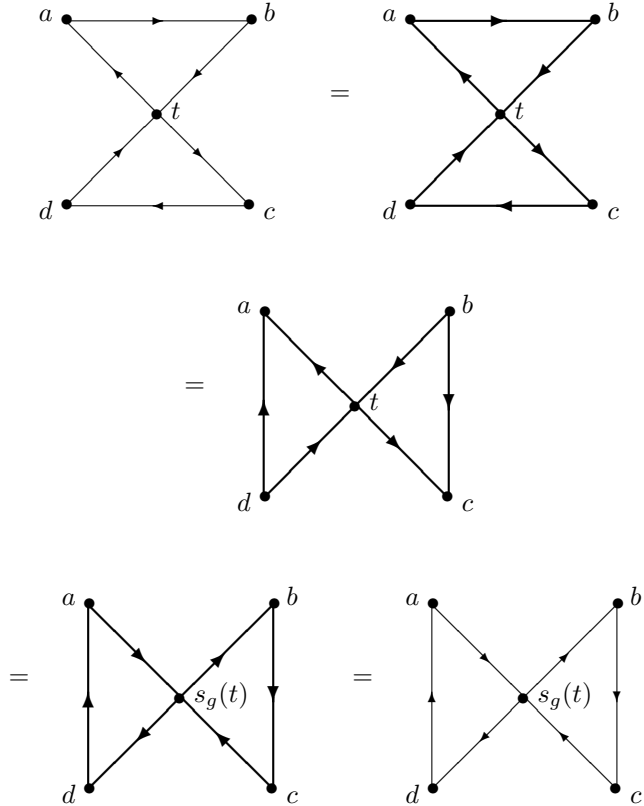
\begin{array}{ccc}
\begin{picture}(100,100)
\thicklines
\put(10,10){\line(0,1){70}}
\put(80,10){\line(0,1){70}}
\put(10,45){\vector(0,1){2}}
\put(80,45){\vector(0,-1){2}}
\put(8,8){$\bullet$}
\put(8,78){$\bullet$}
\put(77,8){$\bullet$}
\put(78,78){$\bullet$}
\put(0,5){$d$}
\put(85,5){$c$}
\put(85,80){$b$}
\put(0,80){$a$}
\end{picture}
&
\begin{picture}(10,90)
\put(0,50){=}
\end{picture} 
&
\begin{picture}(100,100)
\thicklines
\put(10,10){\line(1,0){70}}
\put(10,80){\line(1,0){70}}
\put(45,10){\vector(-1,0){2}}
\put(45,80){\vector(1,0){2}}
\put(8,8){$\bullet$}
\put(8,78){$\bullet$}
\put(77,8){$\bullet$}
\put(78,78){$\bullet$}
\put(0,5){$d$}
\put(85,5){$c$}
\put(85,80){$b$}
\put(0,80){$a$}
\end{picture}.
\end{array}
$$
Hence
$$
\begin{array}{ccc}
\carrehoriz{$t$}
&
\begin{picture}(10,90)
\put(0,50){=}
\end{picture} 
&
\begin{picture}(100,100)
\thicklines
\put(10,10){\line(1,0){70}}
\put(10,80){\line(1,0){70}}
\put(10,10){\line(1,1){70}}
\put(10,80){\line(1,-1){70}}
\put(30,30){\vector(1,1){2}}
\put(60,60){\vector(-1,-1){2}}
\put(45,10){\vector(-1,0){2}}
\put(45,80){\vector(1,0){2}}
\put(60,30){\vector(1,-1){2}}
\put(30,60){\vector(-1,1){2}}
\put(8,8){$\bullet$}
\put(8,78){$\bullet$}
\put(77,8){$\bullet$}
\put(78,78){$\bullet$}
\put(42,42){$\bullet$}
\put(0,5){$d$}
\put(85,5){$c$}
\put(85,80){$b$}
\put(0,80){$a$}
\put(50,43){$t$}
\end{picture}
\end{array}
$$

$$
\begin{array}{cc}
\begin{picture}(10,90)
\put(0,50){=}
\end{picture} 
&
\begin{picture}(100,100)
\thicklines
\put(10,10){\line(0,1){70}}
\put(80,10){\line(0,1){70}}
\put(10,10){\line(1,1){70}}
\put(10,80){\line(1,-1){70}}
\put(30,30){\vector(1,1){2}}
\put(60,60){\vector(-1,-1){2}}
\put(10,45){\vector(0,1){2}}
\put(80,45){\vector(0,-1){2}}
\put(60,30){\vector(1,-1){2}}
\put(30,60){\vector(-1,1){2}}
\put(8,8){$\bullet$}
\put(8,78){$\bullet$}
\put(77,8){$\bullet$}
\put(78,78){$\bullet$}
\put(42,42){$\bullet$}
\put(0,5){$d$}
\put(85,5){$c$}
\put(85,80){$b$}
\put(0,80){$a$}
\put(50,43){$t$}
\end{picture}
\end{array}
$$

$$
\begin{array}{cccc}
\begin{picture}(10,90)
\put(0,50){=}
\end{picture} 
&
\carreverticfat{$s_g(t)$}
&
\begin{picture}(10,90)
\put(0,50){=}
\end{picture} 
&
\carrevertic{$s_g(t)$}.
\end{array}
$$
One can therefore choose our diffeomorphism $\varphi$ as 
$$
\varphi=s_g.
$$
\EPf

\begin{rmk}\label{rmk:Ravel}\hspace{-.1cm}
{\rm In the case of the hyperbolic plane endowed with its natural structure 
of symmetric space, it is tempting to ask whether the three-point function 
defined by the symplectic area $SA(x,y,z)$ of the geodesic triangle ${\Delta 
xyz}$ would satisfy properties (i)---(iii) in \dref{WT}. The answer is 
negative. Indeed, properties (ii) and (iii) would imply unboundedness of 
$SA$. This indicates that, when requiring some compatibility between $S$ 
and the symmetries, property (iii) is somehow too strong.
}\end{rmk}

\section{The Weinstein area}

In \cite{W1}, Weinstein proved that, in the case of a Hermitian symmetric 
space $(M,\omega,s)$ of the noncompact type, the phase function $S_W$, 
occurring in the expression of a given invariant WKB-quantization of 
$(M,\omega)$ defined by an oscillatory integral formula of the type
\begin{equation}\label{SW}
u\star v(x)=\int u(y)v(z) a_\hbar(x,y,z) 
e^{\frac{i}{\hbar}S_W(x,y,z)}dy\,dz,
\end{equation}
must be as follows. 

\begin{enumerate}
\item[$\bullet$]
Let $x,y$ and $z$ be three points in $M$ such that 
the following equation admits a solution $t$~:  
$$
t=s_xs_ys_z(t)
$$
($t$ is unique if it exists). 
\item[$\bullet$] 
Let $\Sigma$ be a surface in $M$ 
bounded by the geodesic triangle ${\Delta tAB}$ 
where $A=s_x(t)$ and $B=s_y(A)=s_ys_x(t)\quad$ ($t=s_z(B)$). 
\end{enumerate}

Then, if for 
some formal amplitude of the form 
$$
a_\hbar(x,y,z)=a_0(x,y,z) +\hbar a_1(x,y,z) + \hbar^2 a_2(x,y,z) + 
...\quad,
$$
the product (\ref{SW}) defines an invariant deformation quantization of 
$(M,\omega)$, the value of the ``WKB"-phase function 
$S_W$ on $(x,y,z)$ is given by 

$$
S_W(x,y,z)=-\int_\Sigma\omega.
$$

Practically, the function $S_W$ is hard to compute explicitly; however, 
some two dimensional examples have been treated in \cite{Q}. The problem 
of finding the amplitude $a_\hbar$ is open.

\begin{dfn}
Let $(M,s)$ be a symmetric space with transvection group $G$. A two-point function $u\in C^\infty(M\times M,\R)$
is called \emph{admissible} iff the following two  properties hold

\vspace{2mm}

(a) there exists a point $o\in M$ such that the function $u$ is invariant under the isotropy group $K_o\subset G$ at point $o$: for all $k\in K_o$ and $x,y\in M$:
$$
u(k.x\,,\,k.y)\;=\;u(x,y)\;.
$$

\vspace{2mm}

(b) For all $x,y\in M$:
$$
u(s_ox,y)\;=\;-u(x,y)\;.
$$

\vspace{2mm}

(c) There exists a $\nabla$-normal neighbourhood $u_0$ centered at point $o$ in $M$ such that
for all $x\in u_0$ and $y\in M$:
$$
u(x,s_{\frac{x}{2}}y)\;=\;-u(x,y)
$$
where $\frac{x}{2}$ denotes the unique mid-point on the geodesic line between $o$                                                                                                                                                          
and $x$ in $u_0$.
\end{dfn}

the value of the ``WKB"-phase function 
$S_W$ on $(x,y,z)$ is given by 

$$
S_W(x,y,z)=-\int_\Sigma\omega.
$$

Practically, the function $S_W$ is hard to compute explicitly; however, 
some two dimensional examples have been treated in \cite{Q}. The problem 
of finding the amplitude $a_\hbar$ is open.   

\begin{prop}
The function $S_W$ is admissible.
\end{prop}
\Pf Let $x,y$ and $z$ be three points in $M$, and set
$$
\left\{
\begin{array}{ccl}
Y & = & s_zs_ys_x(Y)\qquad(*)\\
Z & = & s_x(Y)\\
X & = & s_ys_x(Y)=s_y(Z)
\end{array}
\right.
$$
and 
$$
\left\{
\begin{array}{ccl}
\zeta & = & s_{s_x(y)}s_zs_x(\zeta)\qquad(**)\\
\xi & = & s_x(\zeta)\\
\eta & = & s_zs_x(\zeta)=s_z(\xi),
\end{array}
\right.
$$
(provided equations (*) and (**) admit solutions). The only thing we 
really need to show is that the function $S_{W}$ satisfies property 
(iii); that is,
$$
SA(X,Y,Z) = SA(\xi, \eta, \zeta),
$$
where the function $SA$ is defined as follows. For a sequence of points $\{x_i\}_{0\leq i\leq N}$, the 
expression $SA(x_0,...,x_N)$ means $SA(\gamma) = \int_{\gamma}\alpha$, 
where $\gamma$ is the piecewise geodesic  path whose $i^{\mbox{th}}$ geodesic 
segment starts from point $x_i$ and ends at point $x_{i+1\mbox{mod}(N+1)}$, 
and where $\alpha$ is a $1$-form such that $\omega = d\alpha$.

On the first hand, one has~: 
$
s_xs_{s_x(y)}s_z(X)=s_xs_xs_ys_xs_z(X)=s_ys_xs_z(X)=X.
$
Hence $X=\xi$ and $s_{s_x(z)}(Z)=\zeta$. On the second hand, the 
invariance of $SA$ under the symmetries yields~:
$SA(\xi,\eta,\zeta)=SA(s_x(\xi),s_x(\eta),s_x(\zeta))=SA(\zeta,Z,\xi)$. 
Moreover, $SA(\xi,\eta,\zeta)+SA(\zeta,Z,\xi) = SA(\xi,\eta,\zeta,Z)$, 
hence 
$$S_W(x,s_x(y),z)=\frac{1}{2}SA(\xi,\eta,\zeta,Z).$$ Similarly, one 
gets $$S_W(x,y,z)=-\frac{1}{2}SA(\xi,\eta,\zeta,Z).$$

$$
\begin{picture}(220,180)
\put(20,140){\line(1,-2){60}}
\put(80,140){\line(1,-2){60}}
\put(140,140){\line(1,-2){60}}
\put(20,140){\line(1,0){120}}
\put(80,20){\line(1,0){120}}
\put(50,80){\line(1,0){120}}
\put(20,140){\line(3,-2){180}}
\put(80,20){\line(1,2){60}}
\put(50,80){\line(1,2){30}}
\put(80,140){\line(3,-2){90}}
\put(0,150){$X=\xi$}
\put(70,150){$z$}
\put(130,150){$Y=\eta$}
\put(30,75){$y$}
\put(100,72){$x$}
\put(175,80){$s_x(y)$}
\put(65,5){$Z$}
\put(125,5){$s_x(z)$}
\put(210,10){$\zeta$}
\end{picture}
$$

\EPf

\begin{lem}\label{GEOD}
Let $S\in C^\infty(M\times M\times M,\R)$ be an admissible three-point function.
Let $o$ be a point in $M$. Then, the function $u\in C^\infty(M\times M,\R)$ defined by
$$
u(x,y)\;:=\; S(o,x,y)
$$
is two-point admissible.
\end{lem}
\Pf
Let us first choose a normal neighbourhood $u_0$ centered at $o$ with well-defined
smooth mid-point map:
$$
u_0\to u_0:x\mapsto\frac{x}{2}\;.
$$
Then one has for all point $x\in u_0$ and $y\in M$
\begin{eqnarray*}
&&S(x,y,o)\;=\;S(s_{\frac{x}{2}}o,y,o)\;=\;S(s_{\frac{x}{2}}s_{\frac{x}{2}}o,s_{\frac{x}{2}}y,s_{\frac{x}{2}}o)\;=\;S(o,s_{\frac{x}{2}}y,x)\\
&=&u(s_{\frac{x}{2}}y,x)\;=\;S(o,x,y)\;=\;u(x,y)\;.
\end{eqnarray*}
The rest is obvious.
\EPf
Reciprocally, one has
\begin{prop}
Let $(M,s)$ be symmetric space. For simplicity, let us assume it to be a strictly geodesically convex i.e. any two points can be joined
by a unique geodesic segment. Let $u\in C^\infty(M\times M,\R)$ be two-point admissible.
Then, the element $S\in C^\infty(M\times M\times M,\R)$ defined by
$$
S(x,y,z)\;:=\;u(s_{\frac{x}{2}}y,s_{\frac{x}{2}}z)
$$
is three-point admissible.
\end{prop}

\Pf We first note the $s_o$-invariance of $u$:
$$
u(s_ox,s_oy)\;=\;-u(s_ox,s_os_oy)\;=\;u(y,s_ox)\;=\;-u(y,s_os_ox)=u(x,y)\;.
$$
We continue by  observing  that for all point $x$ and $z$ in $M$, the element
$$
s_os_{\frac{x}{2}}s_{\frac{z}{2}}s_{\left(\frac{s_{\frac{z}{2}}(x)}{2}\right)}
$$
belongs to $K_o$.

\noindent Now, 
\begin{eqnarray*}
&&
S(x,y,z)\;=\;u(s_{\frac{x}{2}}(y),s_{\frac{x}{2}}(z))\\
&=&-u(s_{\frac{x}{2}}(z),s_{\frac{x}{2}}(y))=-u(s_{\frac{x}{2}}s_{\frac{z}{2}}\left(s_o
s_{\frac{x}{2}}s_{\frac{z}{2}}s_{\left(\frac{s_{\frac{z}{2}}(x)}{2}\right)} 
\right)^{-1}(o),s_{\frac{x}{2}}(y))\\
&=&
-u(s_{\frac{x}{2}}s_{\frac{z}{2}}s_{\left(\frac{s_{\frac{z}{2}}(x)}{2}\right)}
s_{\frac{z}{2}}
s_{\frac{x}{2}}(o),s_{\frac{x}{2}}(y)) \\
&=& -u(s_{\frac{x}{2}}s_{\frac{z}{2}}s_{\left(\frac{s_{\frac{z}{2}}(x)}{2}\right)}
s_{\frac{z}{2}}(x),
s_{\frac{x}{2}}s_{\frac{z}{2}}s_{\left(\frac{s_{\frac{z}{2}}(x)}{2}\right)}s_{\left(
\frac{s_{\frac{z}{2}}(x)}{2}\right)}
s_{\frac{z}{2}}(y))\\
&=& -u(s_{\frac{z}{2}}(x),s_{\left(\frac{s_{\frac{z}{2}}(x)}{2}\right)}s_{\frac{z}{2}}(y))= 
u(s_{\frac{z}{2}}(x),s_{\frac{z}{2}}(y))\\
&=& S(z,x,y)\;=\;-u(s_{\frac{z}{2}}(y),s_{\frac{z}{2}}(x))\;=\;-S(z,y,x)\;.
\end{eqnarray*}
Similarly, for every $x,w\in M$, the element
$$
s_{\frac{x}{2}}s_ws_{\frac{s_wx}{2}}
$$
stabilizes $o$. Hence:
\begin{eqnarray*}
&&
S(s_wx,s_wy,s_wz)\;=\;u(s_{\frac{s_wx}{2}}s_wy,s_{\frac{s_wx}{2}}s_wz)\\
&=&
u(s_{\frac{x}{2}}s_ws_{\frac{s_wx}{2}}s_{\frac{s_wx}{2}}s_wy,s_{\frac{x}{2}}s_ws_{\frac{s_wx}{2}}s_{\frac{s_wx}{2}}s_wz)\\
&=&u(s_{\frac{x}{2}x}y,s_{\frac{x}{2}x}z)\;=\;S(x,y,z)\;.
\end{eqnarray*}
At last:
\begin{eqnarray*}
&&
S(x,s_xy,z)\;=\;u(s_{\frac{x}{2}}s_xs_{\frac{x}{2}}s_{\frac{x}{2}}y,s_{\frac{x}{2}}z)\;=\;
u(s_{s_{\frac{x}{2}}x}s_{\frac{x}{2}}y,s_{\frac{x}{2}}z)\\
&=&u(s_os_{\frac{x}{2}}y,s_{\frac{x}{2}}z)\;=\;-u(s_{\frac{x}{2}}y,s_{\frac{x}{2}}z)\;=\;-S(x,y,z)\;.
\end{eqnarray*}
\EPf
We end this section by a remark and a consequent definition. The notion of 2-point admissibility
on  a function $u$ as above entails that for every integer $r$:
$$
u(x,y)\;=\;u(x,s_{\frac{x}{2}}s_oy)\;=\;u(x,\left(s_{\frac{x}{2}}s_o\right)^ry)
$$
for all point $x$ and $y$ in $M$. Combining this fact with Lemma \ref{GEOD} and setting
$$
x\;=\;\mbox{\rm Exp}_o^\nabla(X)\quad (X\in\p\;:=\;\p_o)\;,
$$
we get
$$
u(\exp(X).o\,,\,\exp(rX).y)
$$
for every $y\in M$.
This motivates the following slightly stronger notion of admissibility.
\begin{dfn} Let $(M,s)$ be a (strictly geodesically complete) symmetric space.
A $K_o$-invariant skewsymmetric element $u\in C^\infty(M\times M,\R)$ is called \emph{regular admissible}
whenever for every $X\in\p$, every $t\in\R$ and $x\in M$:
$$
u(\exp(X).o\,,\,x)\;=\;-u(\exp(X)\,,\,s_ox)\;=\;u(\exp(X).o\,,\,\exp(tX).x)\;.
$$
\end{dfn}
\begin{exs}
(i) Given a real finite dimensional vector space $V$ (viewed as a flat symmetric space) together 
with a bilinear 2-form $\Omega\in\bigwedge^2(V^\star)$, the element
$$
u_0(x,y)\;:=\;\Omega(x,y)\quad(x,y\in V)
$$
is regular admissible. The associated admissible three-point function
$$
S_0(x,y,z)\;=\;\Omega(x,y)+\Omega(y,z)+\Omega(z,x)
$$
satifies the cocycle condition of Proposition \ref{GEOM}.

\vspace{2mm}

(ii) Given a non-compact irreducible Riemannian symmetric space $M\;=\;G/K$. One observes that the Iwasawa decomposition of its (necessarily reductive) transvection group $G\;=\;ANK$ implies that a regular admissible
element $u\in C^\infty(M\times M,\R)$ is entirely determined by the data of a smooth element $f\in C^\infty(M)$ through the formula
$$
f(anK)\;:=\;u(aK\,,\,nK)
$$
where $a\in A$ and $n\in N$.
\end{exs}
Based on this remark, we formulate
\begin{dfn}
We call $A$-\emph{orbit} in $M$ any path of the form
$$
t\mapsto\exp(tX).x
$$
with $x\in M$ and $X\in\p$.
\end{dfn}

\section{The \v{S}evera area}

Let us consider a Hamiltonian SSS $(M,s,\omega)$ with transvection group $G$.
We adopt the notations introduced in Section \ref{}.
For every $x$ in $M$, the linear injection
$$
i_x:\p_x\hookrightarrow\g
$$
into the transvection Lie algebra induces dually a $K_x$-equivariant linear projection
$$
\g^\star\to\p_x^\star\;.
$$
Composing the latter with the symplectic isomorphism
$
\p_x^\star\to\p_x
$
induced by the $\omega_x$ and pre-composing with the moment map $m:M\to\g^\star$, one gets the map:
\begin{equation}\label{SEVERAPROJ}
p_x:M\to\p_x\;.
\end{equation}
Given two points $y$ and $z$ in $M$, the \v{S}evera area is defined as the flux of $\omega_x$
through the Euclidean triangle with vertices $\{0\,,\,p_xy\,,\,p_xz\}$ in $\p_x$. Formaly, one makes
\begin{dfn}
Given three points $x,y$ and $z$ in $M$, one denotes by
\begin{equation}\label{SEVERAAREA}
\mbox{\rm Triang}(0\,,\,p_xy\,,\,p_xz)
\end{equation}
the oriented Euclidean triangle with vertices $(0\,,\,p_xy\,,\,p_xz)$.

\noindent Choosing any primitive one-form $\alpha_x\in\Omega^1(\p_x)$ for $\omega_x$:
$$
{\rm d}\alpha_x\;=\;\omega_x\;,
$$
the \emph{\v{S}evera area} of $x, y$ and $z$ is defined by
$$
S_M(x,y,z)\;:=\;\oint_{\mbox{\rm Triang}(0\,,\,p_xy\,,\,p_xz)}\alpha_x\;.
$$
The smooth function
$$
S_M:M\times M\times M\to\R:(x,y,z)\mapsto S_M(x,y,z)
$$
is called the \v{S}evera area of $M$.
\end{dfn}
\begin{prop} Provided $\dim M=2$, 
the \v{S}evera area of $M$ is regular admissible.
\end{prop}
\Pf Let $\g=\k\oplus\p$ be the transvection iLa of $M$ corresponding to the choice of a base point $o\in M$. Since $\dim\p=2$ and $[\p,\p]=\k$, one has
$\dim\k=1$ (it cannot be zero by the Hamitonian condition) i.e. $$\dim\g=3\;.$$
Denote by $\CO$ the coadjoint orbit in $\g^\star$ image of the moment map $m$ of $M$.
For every $x\in M$ and every non-zero $X\in \p$, one has, for every $t\in \R$:
$$
<m(\exp(tX).x)\,,\, X>\;=\;<m(x)\,,\,\Ad_{\exp(-tX)}X>\;=\;<m(x)\,,\,X>\;.
$$
In other words, the curve $\R\to\g^\star:t\mapsto m(\exp(tX).x)-m(x)$
is valued in the  2-plane $X^{*\hspace{-1.8mm}\perp}$ orthodual to the element $X$.

\noindent Now, the kernel of the projection $i_o^\star:\g^\star\to\p^\star$ exactly consists 
in the space $\p^{*\hspace{-1.8mm}\perp}$ orthodual to $\p$. In particular, it is a vector line
in $X^{*\hspace{-1.8mm}\perp}$.

\noindent Putting all this together, we conclude that the ``$A$-orbit" $t\mapsto m(\exp(tX).x)$
projects under $i_o^\star$ into the straight line in $\p^\star$ consisting in
$i_o^\star m(x)\,+\,X^{*\hspace{-1.8mm}\perp}\cap\p^\star$.

\vspace{3mm}

\noindent In other words, every $A$-orbit $t\mapsto\exp(tX).x$ projects under $p_o$ into a straight line
in $\p$  that is \emph{parallel} to the one which $\{p_o(\exp(tX).o)\}_{t\in\R}$ lives in. 

\vspace{2mm}

Therefore, regarding assertion (i), since for all $v,w\in\p$, the quantity $\Omega(u,v)$ (i.e. $\omega_o(u,v)$)
computes the (signed) Euclidean area of the Euclidean triangle with vertices $0, v$ and $w$, we have for every $t$:
\begin{eqnarray*}
&&
S_M(o,\exp(X).o,\exp(tX).x)\;=\;\Omega(p_o(\exp(X).o),p_o(\exp(tX).x))
\\
&=&\Omega(p_o(\exp(X).o),p_o(x))\;=\;
S_M(o,\exp(X).o,x)\;.
\end{eqnarray*}
We now prove that the element:
$$
u_M(x,y)\;:=\; S_M(o,x,y)
$$
satisfies the other conditions of admissibility.

\vspace{2mm}

\noindent\underline{$K_o$-invariance.} The map
$$
p_o: M\to\p
$$
is $K_o$-equivariant since the moment map $m$ is $G$-equivariant by
hypothesis and the injection $\p\hookrightarrow\g$ is $K_o$-invariant.
The $K_o$-invariance of the element $u_M$ then follows from the $K_o$
invariance of the symplectic element $\Omega\;=\;\omega_o$ on $\p\;=\;T_oM$.

\vspace{2mm}

\noindent\underline{Skew-symmetry.} It is obvious from the skew-symmetry of $\Omega$.

\vspace{2mm}

\noindent\underline{$s_o$-compatibility.} It follows from the identity
$$
p_o\circ s_o\;=\;-\id_{\p}\;.
$$
To conclude, we observe that the maps $p_x$ ($x\in M$) enjoy the property that
for every $g\in G$:
$$
p_{g.x}(g.y)\;=\;Ad_gp_xy\;.
$$
The \v{S}evera area $S_M$ is therefore invariant under the diagonal action of $G$ on $M\times M\times M$.
Hence it is determined by $u_M$ via $S_M(g.o,x,y)\;=\;u_M(g^{-1}.x,g^{-1}.y)$. \EPf

\section{The hyperbolic plane}
In this section we treat the case of the hyperbolic plane. We prove in particular that
admissible regular phases are far from being unique.

\noindent Let $S$ be a three-point regular admissible function on the Lobatchevsky plane
$\D\;=\;G/K$ ($G=\SLdr$, $K=SO(2)$).
It determines a regular admissible two-point function $u$ through
$$
u(x,y)\;:=\;S(o,x,y)\quad(o\;:=\;eK).
$$
The two-point function $u$ being $K$-invariant, it is itself determined by its restriction
to $A.o\times \D$ where $G=ANK$ denotes the Iwasawa decomposition considered earlier.
Moreover, Property (\ref{}) of regular admissibility that it is in fact completely determined by 
its restriction to $A.o\times N.o$
which naturally identifies with $\D\;\simeq\;\S=AN$. In other words, the element 
$f\in C^\infty(\S)$ defined by
\begin{equation}\label{ONEPADM}
f(an)\;:=\;u(a.o,n.o)\;=\;S(o,a.o,n.o)
\end{equation}
determines $S$.

\noindent Proposition \ref{} explains which 2-point functions on $\D$ come through this restriction procedure 
from regular admissible three-point phases on $\D$. Pushing one step further: what one-point functions $f$
do come from two-point regular admissible ones?
This is the question we answer to in the present paragraph.

\vspace{3mm}

\noindent We start with a few lemmata.

\begin{lem}\label{KAK}
For every $n\in N$, there is a unique $a_n\in A$ and a unique $k_n\in K$ such that
the following equality holds in $\D$:
$$
nK\;=\;k_na_nK\;.
$$
\end{lem}
\begin{dfn}
We define the \emph{dressing actions} of $K$ on $\S$ and of $\S$ on $K$ by
$$
K\times\S\to\S:(k,s)\mapsto s^k
$$
and 
$$
\S\times K\to K:(s,k)\mapsto k^s
$$
through the relations
$$
ks\;=:\;s^kk^s\;\in G\;.
$$
Also for every $g\in G$, denoting by 
$$
g\;=:\;g_Ag_Ng_K
$$
its Iwasawa decomposition w.r.t. $G\;=\;ANK$,
we define the map
$$
\Phi:\S\to\S:an\mapsto a_n\left((a^{-1})^{k_n^{-1}}\right)_N\;.
$$
\end{dfn}

\begin{lem}\label{PHINVAR}
$$
\Phi(a,n)\;=\;\left(\,\mbox{\rm arcsinh}(n/2)\,,\,\frac{\sinh(2a)}{\sqrt{1+n^2/4}}\,\right)\;.
$$
\end{lem}
\Pf
We start by determining, for every  $n\in N$, $k_n$ and $a_n$.
Which amounts to find $k_0\in K$ such that 
$$
n\;=\;k_na_nk_0
$$
as dictated by the so-called ``$KA^+K$-decomposition".
Setting 
$$
k_n\;=:\;\left(\begin{array}{cc}\cos(\theta)&-\sin(\theta)\\
\sin(\theta)&\cos(\theta)\end{array}\right)\;;\;
k_0\;=:\;\left(\begin{array}{cc}\cos(\theta_0)&\sin(\theta_0)\\
-\sin(\theta_0)&\cos(\theta_0)\end{array}\right)\;;\;
$$
and
$$
a_n\;=:\;\left(\begin{array}{cc}e^r&0\\
0&e^{-r}\end{array}\right)\quad(r>0)\;,
$$
this amounts to the linear system
$$
\left(
\begin{array}{cccc}
1&0&-e^{r}&0\\
0&1&0&-e^{-r}\\
n&1&0&-e^r\\
1&-n&-e^{-r}&0
\end{array}
\right)
\left(
\begin{array}{c}
\cos(\theta)\\
\sin(\theta)\\
\cos(\theta_0)\\
\sin(\theta_0)
\end{array}
\right)\;=\;0\;.
$$
The latter is equivalent under Gauss' algorithm to
$$
\left(
\begin{array}{cccc}
1&0&-e^{r}&0\\
0&1&0&-e^{-r}\\
0&0&ne^r&-2\sinh(r)\\
0&0&2\sinh(r)&-ne^{-r}
\end{array}
\right)
\left(
\begin{array}{c}
\cos(\theta)\\
\sin(\theta)\\
\cos(\theta_0)\\
\sin(\theta_0)
\end{array}
\right)\;=\;0\;.
$$
This implies 
$$
\det
\left(
\begin{array}{cc}
ne^r&-2\sinh(r)\\
2\sinh(r)&-ne^{-r}
\end{array}
\right)\;=\;0\;
$$
i.e.
$$
n\;=\;\epsilon2\sinh(r)\quad(\epsilon\;:=\;\mbox{sign}(n))\;.
$$
Using the above condition, the system is then solved by
$$
\tan(\theta)\;=\;\epsilon e^{-r}\;.
$$
Now the Iwasawa decomposition is  given by
$$
\left(
\begin{array}{cc}
a&b\\
c&d
\end{array}
\right)\;=\;\left(
\begin{array}{cc}
\frac{1}{\sqrt{c^2+d^2}}&0\\
0&\sqrt{c^2+d^2}
\end{array}
\right)\left(
\begin{array}{cc}
1&bd+ac\\
0&1
\end{array}
\right)\left(
\begin{array}{cc}
\frac{d}{\sqrt{c^2+d^2}}&\frac{-c}{\sqrt{c^2+d^2}}\\
\frac{c}{\sqrt{c^2+d^2}}&\frac{d}{\sqrt{c^2+d^2}}
\end{array}
\right)\;.
$$
Therefore, applying it to $g=k_n^{-1}a$ yields
$$
(a^{k_n^{-1}})_N\;=\;-\sinh(2a)\sin(2\theta)\;=\;\frac{-\epsilon\sinh(2a)}{\cosh(r)}\;=\;\frac{-\epsilon\sinh(2a)}{\sqrt{1+n^2/4}}\;.
$$
\EPf

\begin{prop}
Let $f_0:\S\to\R$ be a measurable function. Then, it is the one-point function associated 
to a (measurable) regular admissible function on the Lobatchevsky plane $\D$ (cf. (\ref{ONEPADM})) if and only if it satisfies the
two following properties:

(a)  $$\Phi^\star f_0\;=\;f_0\;,$$
and

(b) $$f_0(a,n)\;=\;-f_0(-a,n)\;.$$
\end{prop}
\Pf Starting with any function $f_0$ it is easy to to define a two point function $u_0$ that is $K$-invariant and
that possesses Property (?) of regular admissibility. Indeed, based on both $G=KA^+K$ polar decomposition
and Iwasawa decomposition, one sets, for all points $x_j\;=\;a_jn_j\in\S\subset G$ ($j=1,2$):
$$
x_1\;=:\;k_1a'_1k_1'\quad\mbox{\rm and}\quad a'_2n'_2k_2\;:=\;k_1^{-1}x_2\;.
$$
Then we define:
$$
u_0(x_1.o,x_2.o)\;:=\;f_0(a'_1n'_2)\;.
$$
The above element $u_0$ satisfies the requirement. Indeed, for every $k\in K$, one has
$$
kx_1\;=\;kk_1a_1'k'_1\quad\mbox{\rm and}\quad(kk_1)^{-1}kx_2\;=\;a'_2n'_2k_2\;.
$$
Hence the $K$-invariance. 
Also, for every $a\in A$, one has
$$
u_0(a'_1.o,ax_2.o)\;=\;u_0(a'_1.o,aa_2n_2.o)\;=\;f_0(a'_1n_2)\;=\;u_0(a'_1.o,x_2.o)\;.
$$
Now, we address the question of skew symmetry for $u_0$. It amounts to the property that
$$
u_0(x_1.o,a_2.o)\;=\;-u_0(a_2.o,x_1.o)\;.
$$
Indeed:
$$
u_0(x_2.o,x_1.o)\;=\;u_0(k_1^{-1}x_2.o,a_1'.o)\;=\;-u_0(a_1'.o,k_1^{-1}x_2.o)\;=\;-u_0(x_1.o,x_2.o)\;.
$$
Now, write 
$$
x_1.o\;=:\;k_1a_1.o
$$
as in Lemma \ref{KAK}. Skewsymmetry then amounts to
\begin{eqnarray*}
&&
u_0(k_1a_1.o,a_2.o)\;=\;u_0(a_1.o,(a_2^{k_1^{-1}})_N.o)\;=\;f_0(a_1(a_2^{k_1^{-1}})_N)\\
&=&-u_0(a_2.o,k_1a_1.o)\;=\;-u_0(a_2.o,(a_1^{k_1})_N.o)\;=\;-f_0(a_2(a_1^{k_1})_N)\;.
\end{eqnarray*}
Now, expressing any $n\in N$ as in  Lemma \ref{KAK}, we get 
$$
n.o\;=\;(a_n)^{k_n}_A(a_n)^{k_n}_N.o\;,
$$
and the Iwasawa decomposition yields:
$$
n\;=\;(a_n)^{k_n}_N\;.
$$
Hence, from the above consideration, we get:
$$
f_0(an)\;=\;f_0(a(a_n)^{k_n}_N)\;=\;-f_0(a_n(a_n)^{k_n^{-1}}_N)\;.
$$
We then conclude from 
$$
f_0(an)\;=\;-f_0(s_oan)\;=\;f_0(a^{-1}n)\;.
$$ \EPf
We now give explicit formulae in the parametric models  of $\D$ introduced 
in Example \ref{PLP}.
\begin{prop}
Consider the one-parameter family of Lie groups $\{G_\lambda\}_{\lambda\in[0,\infty)}$
introduced in Example \ref{PLP}.

\noindent (i) Let $\lambda, R>0$ and realize the Lobatchevsky plane $\D$ as the co-adjoint orbit
$S_{\lambda,R}$ living in $(\g_\lambda^\star,\eta_\lambda)$.

\noindent Then, the hyperbolic \v{S}evera area $S_{S_{\lambda,R}}$ is proportional
to the regular admissible three-point phase $S_\D$ associated under (\ref{ONEPADM}) to the 
element $f_\D\in C^\infty(\S)$  given by
$$
f_\D(a,n)\;:=\;\sinh(2a)\,n\;.
$$

\vspace{2mm}

\noindent (ii) The \v{S}evera area $S_\S$ of the canonical symplectic symmetric space structure
on the $ax+b$ group $\S$ is proportional
to the regular admissible three-point phase $S_\S$ associated to the element $f_\S\in C^\infty(\S)$  given by
$$
f_\S(a,n)\;:=\;\sinh(2a)\,n\;.
$$
\end{prop}
\Pf Consider the element $Z=E-F$ as in Example \ref{PLP}. For every element $x\;=:\;an\;=\;\exp(aH)\exp(nE)$
in $\S$, one computes that
\begin{eqnarray*}
 \Ad_xZ\;=\;-\lambda nH+\frac{1}{2}\left(e^{2a}(1+\lambda n^2)-e^{-2a}\right)H^\perp+
 \frac{1}{2}\left(e^{2a}(1+\lambda n^2)+e^{-2a}\right)Z
\end{eqnarray*}
with
$$
H^\perp\;:=\;E+F\;.
$$
Therefore:
$$
(\Ad_aZ)_\p\;=\;\sinh(2a)H^\perp\quad\mbox{\rm and}\quad(\Ad_nZ)_\p\;=\;-\lambda n H+\frac{\lambda n^2}{2}H^\perp\;.
$$
Hence, any skewsymmetric bilinear 2-form on $\p$ evaluates on $((\Ad_aZ)_\p\,,\,(\Ad_nZ)_\p)$
proportional to $f_\D(a,n)\;=\;\sinh(2a)n$. And the same holds for any multiple of $Z$.

\noindent Since the linear isomorphism $\sharp:\g_\lambda^\star\to\g^\lambda$ is $G_\lambda$-equivariant,
one gets item (i). The other assertion then follows from the fact that the element $f_\D$ is $\lambda$-independant, hence it yields the same formula at the contraction limit $\lambda=0$. \EPf
\begin{rmk}
{\rm
Using the explicit formula given in Lemma \ref{PHINVAR}, for $\lambda=1$, 
the invariance property $\Phi^\star f_\D\;=\;f_\D$ is manifest. 
}
\end{rmk}

\section{Elementary group type Hermitean spaces}
\noindent We now pass to the particular case of a given $2d+2$-dimensional elementary normal 
$\bf j$-group $\S$ with associated 
symplectic form $\omega^\S$.
Let $a,t\in\R$ and $v\in V\simeq\R^{2d}$.  The following identification will always be understood:
$$
\R^{2d+2}\ni x\,:=\,(a,v,t)\mapsto aH+v+tE\in\s\,.
$$
The following result is extracted from \cite{BM1}:
\begin{prop}
\label{proprietes}
\noindent (i) The map
\begin{equation}
\label{chartS}
\s\to\S\,,\,\quad(a,v,t)\mapsto \exp(aH)\exp(v+tE)=\exp(aH)\exp(v)\exp(tE)\;,
\end{equation}
is a global Darboux chart on $(\S,\omega^\S)$ in which the symplectic structure  reads:
$$
\omega^\S\;:=\;\;2{\rm d}a\wedge {\rm d}t\,+\,\omega_0\;.
$$
(ii) Setting furthermore
$$
s_{(a,v,t)}(a',v',t')\;:=\;\big(2a-a',2\cosh(a-a')v-v',2\cosh(2a-2a')t+\omega_0(v,v')\sinh(a-a')-t'\big)\;,
$$ 
defines a symplectic symmetric space structure $(\S,s,\omega^\S)$ on the elementary normal 
$\bf j$-group $\S$.

\noindent (iii) The left action $L_x:\S\to\S:x'\mapsto x.x'$ defines a injective Lie group homomorphism
$$
L:\S\to\Aut(\S,s,\omega^\S)\;.
$$
In the coordinates \eqref{chartS}, we have
$$
x.x'\,=\,(a,v,t).(a',v',t')\,=\,\big(a+a',e^{-a'}v+v',e^{-2a'}t+t'+\tfrac12e^{-a'}\omega_0(v,v')\big)\,.
$$
and
$$
x^{-1}\,=\,(a,v,t)^{-1}=(-a,-e^a v,-e^{2a}t)\,.
$$
(iv) The action $\bR\;:\;\mbox{\rm Sp}(V,\omega^0)\times \S\to\S$, 
$(A,(a,v,t))\mapsto\bR_A(a,v,t):=(a,Av,t)$ by automorphisms of the normal $\bf j$-group $\S$
induces an injective Lie group homomorphism:
$$
\bR:\mbox{\rm Sp}(V,\omega^0)\to\Aut(\S,s,\omega^\S)\,,\quad A\mapsto\bR_A\;.
$$
\end{prop}
Note that  in these coordinates the modular function of $\S$, $\Delta_{\S}$, reads $e^{(2d+2)a}$.

\noindent We now pass to the definition of the three-point phase on $\S$. For this we need 
the notion of ``double geodesic triangle" as introduced by A. Weinstein \cite{W1} and Z. Qian
\cite{Q}.
\begin{dfn} Let $(M,s,\omega)$ be a symplectic symmetric space.
A {\bf midpoint map} on $M$ is a smooth map
$$
M\times M\to M\,,\quad(x,y)\mapsto\,\mid(x,y)\,,
$$
such that, for all points $x,y$ in $M$:
$$
s_{\mid(x,y)}(x)\;=\;y\;.
$$
\end{dfn}
\begin{rmk}
{\rm
Observe that in the case where the partial maps
$s^y:M\to M$, $x\mapsto s_x(y)$
are global diffeomorphisms of $M$, a midpoint map exists and is given by:
$$
\mid(x,y)\;:=\;\left(\,s^x\,\right)^{-1}(y)\;.
$$
Note that in this case, every $\vf\in\Aut(M,s,\omega)$ intertwines the midpoints. Indeed, since for all 
$x,y\in M$ we have $\vf(s_y(x))=s_{\vf(y)}\big(\vf(x)\big)$, we get
$$
\vf\big(\mid(x,y)\big)\;=\;\mid\big(\vf(x)\,,\,\vf(y)\big)\;.
$$
}
\end{rmk}
An immediate computation shows that a midpoint map always exists on the symplectic symmetric
space attached to an elementary normal $\bf j$-group:
\begin{lem}
For the symmetric space $(\S,s)$  underlying an elementary normal $\bf j$-group, the associated 
partial maps are global  diffeomorphisms. In the coordinates \eqref{chartS}, we have:
\begin{equation}
\label{midM}
\big(s^{(a_0,v_0,t_0)}\big)^{-1}:(a,v,t)\mapsto\Big(\frac{a+a_0}2,\frac{v+v_0}{2\cosh(\tfrac{a-a_0}2)},
\frac{t+t_0}{2\cosh({a-a_0})}-\omega_0(v,v_0)\tfrac{\sinh(\tfrac{a-a_0}2)}{4\cosh({a-a_0})
\cosh(\tfrac{a-a_0}2)} \Big)\,.
\end{equation}
\end{lem}
The following statement is proven in \cite{B1}.
\begin{prop}
\noindent (i) The affine space $(\S,\nabla)$ is strictly geodesicaly complete, i.e. two points determine 
a unique geodesic arc.

\noindent (ii) The ``double triangle" three-point function
$$
\Phi:\S^3\to\S^3\,,\quad (x_1,x_2,x_3)\mapsto \big(\mid(x_1,x_2),\mid(x_2,x_3),
\mid(x_3,x_1)\big)\;,
$$
is a $\S$-equivariant (under the left regular action) global diffeomorphism. 
\end{prop}
Since our space $\S$ has trivial de Rham cohomology in degree two,  any three points $(x,y,z)$ 
define an oriented
geodesic triangle $T(x,y,z)$ whose symplectic area is well-defined by integrating the two-form 
$\omega$
on any surface admitting $T(x,y,z)$ as boundary. With a slight abuse of notation, we set
$$
{\rm Area}(x,y,z):=\int_{T(x,y,z)}\omega^\S\;.
$$
\begin{dfn}
The {\bf canonical two-point phase} associated to an elementary normal $\bf j$-group  is defined by
$$
S_{\rm can}^\S(x_1,x_2):={\rm Area}\left(\Phi^{-1}(e,x_1,x_2)\right)\;\in\;C^\infty(\S^2,\R)\;,
$$
where $e:=(0,0,0)$ denotes the unit element in $\S$.
In the coordinates \eqref{chartS}, one has the explicit expression:
\begin{equation}\label{SEXPL}
S_{\mathrm{can}}^\S(x_1,x_2)=\sinh(2a_1)t_2-\sinh(2a_2)t_1+\cosh (a_1) \cosh( a_2)\,\omega_0(v_1,v_2)\;.
\end{equation}
\noindent The {\bf canonical two-point amplitude} associated to an elementary normal $\bf j$-group is defined by
$$
A_{\mathrm{can}}^\S(x_1,x_2)
:=\mbox{\rm Jac}_{\Phi^{-1}}(e,x_1,x_2)^{1/2}\;\in\;C^\infty(\S^2,\R)\;.
$$
In the coordinates \eqref{chartS}, it reads
\begin{equation}\label{AEXPL}
A_{\mathrm{can}}^\S(x_1,x_2)
=\big(\cosh (a_1)\,\cosh (a_2)\,\cosh(a_1-a_2)\big)^d\big(\cosh (2a_1)\,\cosh (2a_2)\,
\cosh( 2a_1-2a_2)\big)^{1/2} \;.
\end{equation}
\end{dfn}

\section{Symplectic symmetric geometry of the group $ax+b$}
We will denote by $\S$ the real two-dimensional connected simply connected Lie group whose Lie algebra $\s$ is generated by two elements $H$ and $E$ with table 
$$
[H,E]\;=\;2E\;.
$$
We will refer to $\S$ as the 
\emph{$ax+b$-group}. 

\noindent Throughout the present book, the group $ax+b$ will be coordinated via the map
$$
\R^2\;\simeq\;\s\longrightarrow\S:(a,n)\;:=\;aH+nE\mapsto \exp(aH).\exp(nE)\;=:\;a.n\;.
$$
The above coordinates have the advantage to constitute a global Darboux chart with respect to any left-invariant symplectic volume form $\omega^{(m)}$
on $\S$:
$$
\omega^{(m)}\;=\;m\,{\rm d}a\,\wedge\,{\rm d}n\quad(m\in\R_0)\;.
$$
In these coordinates, the group law reads
$$
(a,n),(a',n')\;=\;(a+a'\,,\,e^{-2a}n+n')
$$
with unit element
$$
e\;=\;(0,0)
$$
and inverse 
$$
(a,n)^{-1}\;:=\;(-a,-e^{2a}n)\;.
$$
\noindent As a Lie group, $\S$ is canonically endowed with the symmetric space structure (cf. Definition \ref{SS}):
$$
\hat{s}:\S\times\S\to\S:(x,y)\mapsto \hat{s}_xy\;:=\;xy^{-1}x\;.
$$
\begin{rmk}
{\rm 
The Haar element $\omega^{(m)}$ is not invariant under the symmetries $\{\hat{s}_x\}_{x\in\S}$. However there is a structure of symplectic symmetric space
on $\S$ canonically associated to its Lie group structure, which we now describe.
}
\end{rmk}
The above symmetric space is associated with the so called ``exchange" iLg $(\S\times\S,\sigma^e)$ with 
$$
\sigma^e(x,y)\;:=\;(y,x)\;;
$$
the projection map being realized by
$$
\pi:\S\times\S\to\S:(x,y)\mapsto xy^{-1}
$$
whose fibers are the left lateral classes of the diagonal sub-group $K^e\simeq\S$ in $\S\times\S$.

\noindent At the infinitesimal level, it corresponds to the  iLa $(\s\oplus\s,\sigma^e)$ defined by
$$
\sigma^e(X,Y)\;:=\;(Y,X)\;.
$$
However, the above iLa does not realize as such the transvection iLa. Indeed, the isotropy Lie algebra in $\s\oplus\s$ consisting in the diagonal
$$
\k^e\;:=\;\{(X,X)\}
$$
is isomorphic to $\s$
while the holonomy Lie algebra $\k$ is one-dimensional. The transvection Lie algebra may be viewed as a sub-iLa of the exchange iLa as follows.
Considering the decomposition into $\sigma^e$-eigenspaces:
$$
\s\oplus\s\;=\;\k^e\oplus\p
$$
with $\p\;:=\;\{(X,-X)\}$, one has
$$
\k\;=\;[\p,\p]\;.
$$
Hence, setting
$$
\a\;:=\;\R\,H\quad\mbox{\rm and}\quad\n\;:=\;\R\,E\;,
$$
the transvection Lie algebra $(\g,\sigma)$ then identifies with following the semi-direct extension of the Abelian $\n\oplus\n$ by $\a$:
$$
\g\;=\;\a\ltimes(\n\oplus\n)
$$
with table
$$
[H\,,\,E]\;=\;2E\quad;\quad[H,F]\;:=\,-2F\quad\mbox{\rm and}\quad[E\,,\,F]\;=\;0
$$
where, with some slight abuse of notation, we have set
$$
H\;:=\;(H,-H)\quad;\quad E\;:=\;(E,0)\quad\mbox{\rm and}\quad F\;:=\;(0, E)\;.
$$
Within this presentation, the involution is given by
$$
\sigma H\;=\;-H\quad,\quad\sigma E\;=\;F\;.
$$
\begin{rmk}
{\rm 
Note that the transvection Lie algebra $\g$ is isomorphic to the Poincar\'e Lie algebra $\mathfrak{so}(1,1)\ltimes\R^2$.
}
\end{rmk}
The Lie algebra embedding 
$$
\s\to\g:aH+nE\mapsto a(H,-H)\,+\,n(E,0)
$$
exponentiates to realizing the $ax+b$ group $\S$ into $\S\times\S$ under
$$
\S\to\S\times\S:an\mapsto(a,a^{-1}).(n,e)\;=\;(an,a^{-1})\;.
$$
This yields a simply transitive action of $\S$ on $(\S\times\S)/K^e$ and associated diffeomorphism
\begin{equation}\label{ANA}
\S\to\S:(a,n)\mapsto\pi(an,a^{-1})\;=\;ana\;.
\end{equation}
\begin{rmk}
{\rm
The map (\ref{ANA}) is not a group homomorphism.
}
\end{rmk}
As announced, we now have
\begin{prop}
Under the diffeomorphism (\ref{ANA}), the group symmetric space structure $\hat{s}$ reads
$$
s_{(a,n)}(a',n')\;=\;(\,2a\,-\,a'\,,\,2\cosh(2(a-a'))n\,-\,n')\;.
$$
The Haar elements $\omega^{(m)}\;=\;m\,{\rm d}a\wedge{\rm d}n$ are the only 2-forms on $\S$ invariant under the transported symmetries $\{s_x\}_{x\in\S}$.
\end{prop}
\Pf Denoting by $\bC$ the conjugate action $\bC_xy\;=\;xyx^{-1}$, one observes:
\begin{eqnarray*}
\hat{s}_{ana}(a'n'a')=anaa^{\prime -1}n^{\prime -1}a^{\prime -1}ana\;=\;a^2a^{\prime -1}\bC_{a'a^{-1}}n.n^{\prime -1}.\bC_{aa^{\prime-1}}n.a^2a^{\prime -1}\;.
\end{eqnarray*}
The announced formula then follows from the fact that $\bC_an\;=\;(0,e^{-2a}n)$. \EPf
\noindent In other words, the triple $(\S,\omega^{(m)},s)$ is a symplectic symmetric space.

\section{Symplectic fluxes trough geodesic triangles}
\begin{thm}\label{SSAXPB}
(i) The symmetric space $(\S,s)$ admits a \emph{smooth mid-point map}, in the sense that given any two points $x$ and $y$ in $\S$, there exists a unique
point $m(x,y)$ such that
$$
s_{m(x,y)}x\;=\;y
$$
and the map $m:\S\times\S\to\S$ is $C^\infty$.

\noindent (ii) The map
$$
\Phi:\S\times\S\times\S\to\S\times\S\times\S: (x,y,z)\mapsto(\,m(x,y)\,,\,m(y,z)\,,\,m(z,x)\,)
$$
is a diffeomorphism.

\noindent (iii) Setting $\Phi^{-1}(x,y,z)\;=:\;(x',y',z')$ and denoting by $\stackrel{\Delta}{x'y'z'}$ the oriented geodesic triangle
with vertices $x'$, $y'$ and $z'$, we consider the flux of $\omega^{(1)}$ through (any surface bounded by) $\stackrel{\Delta}{x'y'z'}$:
$$
S(x,y,z)\;:=\;\int_{\stackrel{\Delta}{x'y'z'}}\omega^{(1)}\;.
$$
One  has the explicit formula:
$$
S(x_0,x_1,x_2)\;=\;\sinh(2(a_0-a_1))n_2\,+\,\sinh(2(a_2-a_0))n_1\,+\,\sinh(2(a_1-a_2))n_0
$$
where for every $j=0,1,2$, we set $x_j=(a_j,n_j)$.
\end{thm}

\chapter{Oscillatory integrals on exponential Lie groups}\label{OSCINT}

\subsection{Symbol spaces}
We consider  a centerless Lie group $\L$ with Lie algebra $\CL$ and associated enveloping algebra $\CU(\CL)$.
For every $X$ in $\CL$,
we denote by $\tilde{X}$ the associated left-invariant vector field on $\L$ and by $X^\star$ the associated right-invariant
vector field:
$$
\tilde{X}_x\;:=\;\ddto x\exp(tX)\,,\quad X^\star_x\;:=\;\ddto\exp(-tX)x\;.
$$
We use the same notations for the elements of the enveloping algebra $\CU(\CL)$.
\begin{dfn}
Choosing a vector norm $|\;|$ on $\fL$, we will consider the element $\fd_\L\in C^\infty(\L)$ defined by
$$
\fd_\L(x)\;:=\;\sqrt{1\,+\,|\Ad_x|^2\,+\,|\Ad_{x^{-1}}|^2}
$$
where for every $x$ in $\L$, $\Ad_x$ stands for the adjoint action of $x$ in $\s$ whose operator norm is denoted by
$|\Ad_x|$.
\end{dfn}

\noindent The main property of the element $\fd_\L$ is
\begin{lem}
(i) There exists a positive $C$ such that, for all $x,y\in\L$:
$$
\fd_\L(xy)\,<\,C\,\fd_\L(x)\,\fd_\L(y)\;.
$$
(ii) For every $X\in\CU(\fL)$, there exists $C$ such that
$$
|\tilde{X}.\fd_\L|\,<\,C\,\fd_\L\quad\mbox{and}\quad|\underline{X}.\fd_\L|\,<\,C\,\fd_\L\;.
$$
\end{lem}

\noindent Denoting by
$$
\CU(\fL)\;\;=\;\sum_{k\in\N}\CU_k(\fL)
$$
the natural filtration, we make the
\begin{dfn}
Let $N\in\R$ and $p\in\N$. We set
$$
\CB^N_p(\L)\;:=\;\{\,F\in C^p(\fL)\,|\,\forall X\in\CU_p(\fL)\, \exists \,C>0\,:\,|\tilde{X}.F|\,<\,C\,\fd^N_\L\,\}\;.
$$
Similarly, we define
$$
\CB^N(\L)\;:=\;\CB^N_\infty(\L)\;:=\;\{\,F\in C^\infty(\fL)\,|\,\forall X\in\CU(\fL)\, \exists \,C>0\,:\,|\tilde{X}.F|\,<\,C\,\fd^N_\L\,\}\;,
$$
and
$$
\CB(\L)\;:=\;\CB^0(\L)\;.
$$
Endowed with the pointwise multiplication of functions, these spaces becomes a commutative Fr\'echet algebras when equipped
with the following seminorms:
\begin{equation}\label{SEMINORM}
|\varphi|_{J}^{N}\;:=\;\sup_{x\in\S}\{\,|\fd_\L^{-N}(x)\,\widetilde{\Xi_J}_x.\varphi|\,\}
\end{equation}
where we fixed a basis $\{\Xi_j\}$ of $\fL$ and where $J=(j_1,j_2,...,j_d)\in\N^{\dim\L}$ is a multi-index such that $d\leq p$ to which we associate the Poincar\'e-Birkoff-Witt (PBW)
basis element
$\Xi_J\;:=\;\Xi^{j_1}_1...\Xi^{j_d}_d$ in $\CU(\fL)$.
\end{dfn}
Endowed with the seminorms (\ref{SEMINORM}), we observe that the $\CB^N_p$'s become Banach spaces, with Fr\'echet projective limits, over $p$ ($N$ being fixed),
$\CB^N$ and $\CB$. 

\noindent Similarly as in the usual Abelian case of $\L\;=\;\R^n$, we have
\begin{prop}\label{SYMBOLS}
(i) For all $N,N'$ and $p,p'$, the pointwise multiplcation of functions yields separately continuous bilinear maps:
$$
\CB^N_p\times\CB^{N'}_{p'}\to\CB^{N+N'}_{\min\{p,p'\}}\;.
$$
(ii) For every $X\in\CU_k(\fL)$, one has the continuous map
$$
\CB^N_p\to\CB^N_{p-k}:F\mapsto\tilde{X}.F\;.
$$
(iii) 
$$
\CB^N_p\;=\;\fd_\L^N\CB^0_p\quad\mbox{and}\quad \CB^N\;=\;\fd_\L^N\CB^0\;.
$$
(iv) Denoting by $\mbox{\rm adh}_{\CB^N_p}$ the topological closure in $\CB^N_p$, one has, for all $N\,>\,N' \,\geq\,0$ and $\infty\,\geq\,p'\,\geq\,p$:
$$
\mbox{\rm adh}_{\CB^N_p}\CD\;\supset\,\CB^{N'}_{p'}
$$
where $\CD=\CD(\L)$ denotes the test space of  compactly supported, $C^\infty$ functions.

\noindent (v) Let us denote by $C_{b}(\L)$ the $C^\star$-algebra of (right-) uniformly continuous bounded functions on $\L$. This $C^\star$-alegbra
is then strongly continuously and isometrically acted on by $\L$ under the right-regular representation. The (Fr\'echet) subalgebra of smooth vectors
in $C_{b}(\L)$ consists in the space $\CB(\L)$.

\noindent In particular, the space $\CB(\L)$ is a smooth left $\L$-module under $R$.
\end{prop}
\noindent If $(\fF,\{|\,\,|\}_{j\in\N})$ is  a Fr\'echet space whose topology is defined a fixed set of semi-norms $\{|\,\,|\}_{j\in\N}$, there is an obvious $\fF$-valued version
version of the above defined spaces $\CB^N_p$. However, it turns out that a slightly more sophisticated definition is actually more relevant in the applications:
\begin{dfn}\label{NJ}
We set
$$
\CB^{{N}}_p(\L,\fF):=\{F\in C^\infty(\L, \fF)\;|\;\forall j\in\N, X\in\CU_p(\fL)\;:\exists C>0\;:|\tilde{X}.F|\,<\,C\,\fd_\L^{N}\}.
$$
\end{dfn}
As in the complex valued case, these spaces are Fr\'echet spaces in a natural way and the obvious analogue of Proposition \ref{SYMBOLS} holds. Namely,
\begin{prop}\label{SYMBOLSF}
(i) For all $N$ and $p,p'$, the pointwise multiplcation of functions yields separately continuous bilinear maps:
$$
\CB^N_p(\L)\times\CB^{{N}}_{p'}(\L,\fF)\to\CB^{2N}_{\min\{p,p'\}}(\L,\fF)\;.
$$

(ii) For every $X\in\CU_k(\fL)$, one has the continuous map
$$
\CB^{{N}}_p\to\CB^{{N}}_{p-k}:F\mapsto\tilde{X}.F\;.
$$
(iii) Denoting by $\mbox{\rm adh}_{\CB^{{N}}_p}$ the topological closure in $\CB^{{N}}_p$, one has, for all ${{N'}}\prec{{N}}$ and $\infty\,\geq\,p'\,\geq\,p$:
$$
\mbox{\rm adh}_{\CB^{{N}}_p}\CD\;\supset\,\CB^{{N'}}_{p'}
$$
where $\CD=\CD(\L,\fF)$ denotes the test space of  $F$-valued compactly supported, $C^\infty$ functions.
\end{prop}

\subsection{Oscillatory integrals}

\noindent In this section, we define a notion of ``improper integral" that will allow us 
to define our non-formal deformation product on $\CB$-type spaces.
We start with the following rather technical definitions, but which will be basic
in the sequel.

\vspace{3mm}

\noindent {\bf From now on, we consider our Lie group $\L$ to be solvable of exponential type.}

\vspace{3mm}

\begin{dfn}
\label{TEMPPAIR}
The pair $(G,S)$ is called {\bf tempered} if the  map
\begin{equation}\label{TEMPCOORD}
\phi:G\,\to\,\g^\star\,:\,x\,\mapsto\,\Big[\,\g\to\R:X\mapsto\,{\rm d}S_x(\widetilde{X})\;=\;\big(\widetilde{X}.\,S\big)(x)\,\Big]\;,
\end{equation}
is a global diffeomorphism.
\end{dfn}

\noindent
Given a tempered pair $(G,S)$, with $\g$ the Lie algebra of $G$, we now consider a vector space decomposition:
\begin{equation}\label{DECOMP}
\g\;=\;\bigoplus_{n=0}^NV_n\;,
\end{equation}
and for every $n=0,\dots, N$,
an ordered basis $\{e^n_{j}\}_{j=1,\dots,\dim (V_n)}$ of $V_n$. We  get global coordinates on $G$:
\begin{equation}\label{COORDINATES}
x_n^{j}\;:=\;\big(\widetilde{e_{j}^n}. S\big)(x)\,,\qquad n=0,\dots,N\,,\quad j=1,\dots,\dim (V_n)\;.
\end{equation}
We choose a scalar product on each $V_n$ and let $|.|_n$ be the associated Euclidean norm.
We will always identify the universal enveloping algebra $\CU(\g)$ 
with  the symmetric algebra $\fS(\g)$ of $\g$, through the Poincar\'e-Birkhoff-Witt
linear isomorphism.
 Given an element $A\in\CU(\g)$, we let  ${\widetilde A}^*$ 
the formal adjoint of the left-invariant differential operator $\widetilde A$, with respect to the inner product of $L^2\big(G,{\rm d}_G\big)$. We make the obvious observation that ${\widetilde A}^*$  is still left-invariant. Indeed, for $\psi,\vf\in C^\infty_c(G)$ and $g\in G$, we have
$$
\langle L^\star_g {\widetilde A}^*\psi,\vf\rangle=\langle  \psi,{\widetilde A}L^\star_{g^{-1}}\vf\rangle=
\langle  \psi,L^\star_{g^{-1}}{\widetilde A}\vf\rangle=\langle  {\widetilde A}^*L^\star_g\psi,\vf\rangle\;.
$$
Moreover, we make the following requirement of compatibility of the adjoint map on $L^2(G,{\rm d}_G)$ with respect to the ordered decomposition \eqref{DECOMP}: 
\begin{equation}
\label{adjoint}
\forall\, n=0,\dots,N,\quad\forall \,A\in
\fS(V_n)\,,\quad\exists\, B\in \prod_{k=0}^n\fS(V_k)\quad\mbox{such that}\quad
\widetilde A^*=\widetilde B\;.
\end{equation}
\noindent We now pass to regularity assumptions regarding the function $S$. 

\begin{dfn}\label{TEMPADM} (\cite{BG} Definition 2.24 page 26) Set 
\begin{equation}
\label{bE}
\bE:=\exp\{iS\}\;.
\end{equation}
 A tempered pair $(G,S)$ is called {\bf admissible}, if there exists a decomposition (\ref{DECOMP})
with associated coordinate system (\ref{COORDINATES}), such that
for every $n=0,\dots,N$, there
exists an element $X_n\in\fS(V_n)\subset\CU(\g)$ whose associated multiplier $\alpha_n$, defined as 
\begin{equation}
\label{multiplier}
\widetilde{X}_n\bE\;=:\;\alpha_n\bE\;,
\end{equation}
satisfies the following properties:

\noindent (i) There exist $C_n>0$ and $\rho_n>0$ such that:
$$
\left|\alpha_n\right|\;\geq\;C_n\big(1\,+\,\left|x_n\right|_n^{\rho_n}\big)\;,
$$
where $x_n\;:=\;(x_n^{j})_{j=1,\dots,{\rm dim}(V_n)}$.\\
\noindent (ii) For all $n=0,\dots, N$, there exists a tempered function $0<\mu_n\in C^\infty(G)$ such that:

(ii.1) For every $A\in
\prod_{k=0}^n\fS(V_k)\subset
\CU(\g)$ there exists $C_A>0$ such that:
\begin{equation}\label{UNIF}
\big|\widetilde{A}\,{\alpha_n}\big|\;\leq\;C_A\,\left|\alpha_n\right|\,\mu_n\;.
\end{equation}

(ii.2) The function $\mu_n$ is independent of the variables $\{x^{j}_r\}_{j=1,\dots,\dim(V_r)}$, for all $r\leq n$:
\begin{equation}\label{INDEP}
\frac{\partial \mu_n}{\partial{x^{j}_r}}\;=\;0\,,\quad\forall r\leq n\,,\,\,\,\forall j=1,\dots,\dim(V_r)\;.
\end{equation}
\end{dfn}

\noindent From the Proposition 2.28 page 30 and its proof in \cite{BG}, we have
\begin{prop}
For every $N\in\N$, there exists a differential operator $\bD$ and an integer 
$K$ such that for every $M$ and $F\in\CB_K^M(\L)$, one has

(i)
$$
|\bD^\star F|\;\leq\; |F|_{K,M}\;\fd^{-N}
$$
and 

(ii) $$
\bD\left(e^{iS}\right)\;=\;e^{iS}\;.
$$
\end{prop}
Letting $(\fB, ||\;||)$ be a Banach space, we deduce:
\begin{prop}
For every $N\in\N$, there exists a differential operator $\bD$ and an integer 
$K$ such that for every $M$ and $F\in\CB_K^M(\L,\fB)$, one has

(i)
$$
||\bD^\star F||\;\leq\; ||F||_{K,M}\;\fd^{-N}
$$
and 

(ii) $$
\bD\left(e^{iS}\right)\;=\;e^{iS}\;.
$$
\end{prop}

\begin{cor} Let  $\mu\in\CB^N(\L)$ and $K$ sufficiently large. 
Then, for every $M\in\N$, the map
\begin{equation}\label{OSC}
\CD(\L,\fB)\to \fB:\varphi\mapsto\int_{\L}\mu\,e^{iS}\,\varphi
\end{equation}
is  continuous w.r.t. the topology induced on $\CD$ by $\CB^{{M+\epsilon}}_{K}(\L,\fB)$.

\noindent Therefore, for every $K$ sufficiently large, the map
$$
{\int\,\mu\,e^{iS}}\,:\,\CD(\L,\fB)\longrightarrow\fB:F\mapsto\langle\,{\int\,\mu\,e^{iS}}\,,\,F\,\rangle\;.
$$
extends in a unique way as a continuous
linear map:
$$
\widetilde{\int\,\mu\,e^{iS}}\,:\,\CB^{M}_K(\L,\fB)\longrightarrow\fB:F\mapsto\langle\,\widetilde{\int\,\mu\,e^{iS}}\,,\,F\,\rangle\;.
$$.
\end{cor} 
\Pf 
One has, for every large $N'$:
\begin{eqnarray*}
&&
\left|\left|\int_\L\mu e^{iS}\varphi\right|\right|\;=\;\left|\left|\int_\L\mu \bD(e^{iS})\varphi\right|\right|\;=\;\left|\left|\int_\L e^{iS}\bD^\star(\mu\varphi)\right|\right|\;\leq\;\int_\L\left|\left|\bD^\star(\mu\varphi)\right|\right|\\
&\leq&||\mu\varphi||_{K,N'}\int\fd^{-N'}\;=\;C'||\mu\varphi||_{K,N'}\;.
\end{eqnarray*}
Now, one has 
\begin{eqnarray*}
&&
||\mu\varphi||_{K,N'}\;=\;\sup\fd^{-N'}||\tilde{X}(\mu\varphi)||\;=\;
\sup\fd^{-N'}||\sum_{(X)}\tilde{X}_{(1)}(\mu)\tilde{X}_{(2)}(\varphi)||
\\
&\leq&\sum_{(X)}\sup\fd^{-N'}\fd^NC_{(1)}||\tilde{X}_{(2)}(\varphi)||\;\leq\;
C^{\prime\prime}||\varphi||_{K,M+\epsilon}\sup\fd^{N-N'+M+\epsilon}
\end{eqnarray*}
which is bounded as soon as $N'>N+M+\epsilon$. This implies that
for every $\epsilon$, one has 
$$
\left|\left|\int_\L\mu e^{iS}\varphi\right|\right|\;\leq\;C||\varphi||_{K,M+\epsilon}\;.
$$

\EPf

\subsection{Deformations}
We consider now $\fB$ to be a Banach algebra.
\begin{dfn}
We consider the following continuous bilinear map
\begin{eqnarray}\label{RA}
\CR\;:&C(\L,\fB)&\times\quad C(\L,\fB)\to C(\L\times\L,C(\F,\fB))\;:\\
&(F_1,F_2)&\mapsto[\,(x_1,x_2)\mapsto\,[x\mapsto F_1(xx_1).F_2(xx_2)\;=\\
&&\quad\quad \mu_\fB(R^\star\otimes R^\star_{(x_1,x_2)}(F_1\otimes F_2))(x)\,]\,]
\end{eqnarray}
$\mu_\fB$ denotes the mutiplication in $\fB$.
\end{dfn}
\begin{lem}
The restriction of the map (\ref{RA}) to $\CB^N_p(\L,\fB)\times\quad \CB^N_p(\L,\fB)$ is separately continuously  valued in $\CB^{q+N,q+N}_{p-q,p-q}(\L\times\L,\CB^{2N}_{q}(\L,\fB))$:

\begin{eqnarray}\label{RAB}
\CR\;:&\CB^N_p(\L,\fB)&\times\quad \CB^N_p(\L,\fB)\to\CB^{q+N,q+N}_{p-q,p-q}(\L\times\L,\CB^{2N}_{q}(\L,\fB))\;:\\
&(F_1,F_2)&\mapsto[\,(x_1,x_2)\mapsto\,[x\mapsto F_1(xx_1).F_2(xx_2)\,]\,]
\end{eqnarray}
for every $q\leq p$.
\end{lem}
\Pf 
The generic semi-norm on $\CB^{N_1}_{q_1}(\L,\fB)\ni f$ is given by ($X\in\CU_{q_1}(\fL)$):
$$
||f||_{q_1}^{N_1}\;:=\;\sup_x\{\fd_\L^{-N_{1}}(x)||\tilde{X}_x(f)||\}\;.
$$
The generic semi-norm on $\CB^{N_2,N'_2}_{q_2,q'_2}(\L\times\L,\CB^{N_1}_{q_1}(\L,\fB))\ni \Phi$ is therefore given by ($Y\in\CU_{q_2+q'_2}(\fL\oplus\fL)$):
$$
||\Phi||_{q_2,q'_2}^{N_2,N'_2}\;:=\;\sup_{(x_1,x_2)}\{\fd_{\L\times\L}^{-N_2,-N'_2}(x_1,x_2)||\tilde{Y}_{(x_1,x_2)}(\Phi)||^{N_1}_{q_1}\}\;.
$$

We now consider the element $\CR(F_1,F_2) \in C^{q_2,q'_2}(\L\times\L,C^{q_1}(\L,\fB))$.
We want to estimate the quantity:
$$
\sup_{(x_1,x_2)}\{\fd_{\L\times\L}^{-N_2,-N'_2}(x_1,x_2)||\tilde{Y}_{(x_1,x_2)}(\CR(F_1,F_2))||^{N_1}_{q_1}\}\;.
$$
is bounded. We start by observing that
for every $\vec{x}\in\L\times\L$ and $\CR(F_1,F_2)\in C^q(\L\times\L,\fB)$:
$$
\tilde{Y}_{\vec{x}}\CR(F_1,F_2)\;=\;\CR(\tilde{Y}(F_1,F_2))(\vec{x})\;.
$$
Hence:
$$
||\tilde{Y}_{(x_1,x_2)}(\CR(F_1,F_2))||^{N_1}_{q_1}\;=\;||\CR(\tilde{Y}(F_1,F_2))(x_1,x_2)||^{N_1}_{q_1}$$
which equals
$$
\sup_x\{\fd_\L^{-N_1}(x)||\tilde{X}_x\left(\CR(\tilde{Y}(F_1,F_2))(x_1,x_2)\right)||\;.
$$
Considering $Y$ of the form $Y\;=\;Y_1\otimes Y_2$, observe that 
$$
\tilde{X}_x\left(\CR(\tilde{Y}(F_1,F_2))(x_1,x_2)\right)\;=\;\left((\Ad_{x_1^{-1}}X)Y_1\otimes (\Ad_{x_2^{-1}}X)Y_2\right)^\sim_{\vec{x}}\CR(F_1,F_2)(x)\;.
$$
Hence
\begin{eqnarray*}&&
 \left|\left|\tilde{X}_x\left(\CR(\tilde{Y}(F_1,F_2))(x_1,x_2)\right)\right|\right|\leq\\ &&  \left|\left|C\fd^{|X|,|X|}_{\L\times\L}(x_1,x_2)\CR((X\otimes X.Y)(F_1,F_2))(x_1,x_2)(x)\right|\right|\leq\\
&& C'   \left|\left|\fd^{q_1,q_1}_{\L\times\L}(x_1,x_2)\CR((X\otimes X.Y)(F_1,F_2))(x_1,x_2)(x)\right|\right|\leq\\
&&
C'  \fd^{q_1,q_1}_{\L\times\L}(x_1,x_2)\,\fd_\L^{N}(xx_1)\,\fd^{N}_\L(xx_2)\,||F_1||^N_{q_1+q_2}||F_2||^N_{q_1+q'_2}\leq\\
&&
C" \fd^{q_1+N,q_1+N}_{\L\times\L}(x_1,x_2)\,\fd_\L^{2N}(x)\,||F_1||^N_{q_1+q_2}||F_2||^N_{q_1+q'_2}\;.
\end{eqnarray*}
Therefore
\begin{eqnarray*}
&&|\tilde{Y}_{(x_1,x_2)}(\CR(F_1,F_2))|^{N_1}_{q_1}\\
&\leq&C^{\prime\prime\prime}
||F_1||^N_{q_1+q_2}||F_2||^N_{q_1+q'_2}\fd^{q_1+N,q_1+N}_{\L\times\L}(x_1,x_2)\;\sup_x\fd_\L^{2N-N_1}(x)
\end{eqnarray*}
which leads to the estimate:
\begin{eqnarray*}
&&\sup_{(x_1,x_2)}\{\fd_{\L\times\L}^{-N_2,-N'_2}(x_1,x_2)|\tilde{Y}_{(x_1,x_2)}(\CR(F_1,F_2))|^{N_1}_{q_1}\}\\
&\leq&C_j
||F_1||^N_{q_1+q_2}||F_2||^N_{q_1+q'_2}\times\\
&\times&\sup_{(x_1,x_2)}\{\fd_{\L\times\L}^{(-N_2+q_1+N,-N'_2+q_1+N)}(x_1,x_2)\}\;\sup_x\fd_\L^{2N-N_1}(x)
\end{eqnarray*}
which leads to the condition of finiteness:
$$
\begin{array}{ccc}
q_1+q_2&\leq&p\\
q_1+q'_2&\leq&p\\
q_1+N&\leq&N_{2}\\
q_1+N&\leq&N'_{2}\\
2N&\leq&N_1
\end{array}
$$
which setting 
$$
q\;:=\;\min\{q_2,q_2'\}
$$
yields the assertion. \EPf

\noindent We then obtain
\begin{cor}\label{STARS} Let $A\in\CB^{N_0}(\L^2,\C)$. Let $q\in\N$. Then,
provided $p$ is large enough, one has that,
to all $F_1,F_2$ in $\CB^N_p(\L,\fB)$, the oscillating integral  associates an element of $\CB^{2N}_{q}(\L,\fB)$:
$$\star_{S}\;=\;\widetilde{\int_{\L\times\L}A\,e^{-iS}}:\CB^N_p(\L,\fB)\times\CB^N_p(\L,\fB)\to\CB^{2N}_{q}(\L,\fB):(F_1,F_2)\mapsto F_1\star_S F_2\;.$$
Moreover, the product $\star_S$ is bilinear and separately continuous.

\noindent In particular, for $p=\infty$, one has 
$$\star_{S}\;=\;\widetilde{\int_{\L\times\L}A\,e^{-iS}}:\CB^0_\infty(\L,\fB)\times\CB^0_\infty (\L,\fB)\to\CB^0_\infty(\L,\fB)\;.$$
\end{cor}

\begin{dfn}
\label{wass}
We say that the product $\star_S$, given in Corollary \ref{STARS}, is {\bf weakly associative} when for all $\psi_1,\psi_2,\psi_3\in\CD(\L,\fB)$, one has
$(\psi_1\star_S\psi_2)\star_S\psi_3=\psi_1\star_S(\psi_2\star_S\psi_3)$ in $\CB^0_\infty(\L,\fB)$.
\end{dfn}

\begin{prop}\label{WASS}
Weak associativity implies strong associativity in 
the sense that, when weakly associative, 
 for all elements $(F,F',F'')\in\CB^N_p(\L,\fB)\times\CB^N_p(\L,\fB)
\times\CB^N_p(\L,\fB)$,
one has the equality $(F\star_SF')\star_S F''=F\star_S (F'\star_S F'')$ in 
$\CB^{4N+1}_q(\L,\fB)$ for some $q$.
\end{prop}
\Pf
This follows by density of the test functions in $\CB$-type spaces. Indeed, let 
$\{\phi_n\}_{n\in\N}$, $\{\phi'_{n'}\}_{n'\in\N}$ and $\{\phi''_{n''}\}_{n''\in\N}$ be sequences 
of test functions converging to $F,F'$ and $F''$ respectively in $\CB^{N+\epsilon}_p(\L,\fB)$.
By separate continuity, we have that 
$$
F\star_SF'\;=\;\lim_{n,n'}(\phi_n\star_S\phi_{n'})
$$
where the limit on the RHS is taken in $\CB_{q}^{2(N+\epsilon)}(\L,\fB)$ with $q<<p$. Now, since 
$$
\CB^N_p(\L,\fB)\;\subset\;\CB_{q}^{2(N+\epsilon)}(\L,\fB)\;,
$$
the convergence $\lim_{n''}\phi_{n''}\;=\;F''$ holds in $\CB_{q}^{2(N+\epsilon)}(\L,\fB)$ as well.
Hence, for every $q'<<q<<p$, one has, again by separate continuity
$$
(F\star_SF')\star_S F''\;=\;\lim_{n,n',n''}(\phi_n\star_S\phi'_{n'})\star(\phi''_{n''})
$$
where the limit in the RHS is taken in $\CB_{q'}^{4(N+\epsilon)}(\L,\fB)$. 

\noindent The same argument for $F\star_S (F'\star_S F'')$, yields 
$$
F\star_S(F'\star_S F'')\;=\;\lim_{n,n',n''}\phi_n\star_S(\phi'_{n'}\star\phi''_{n''})\;.
$$
And we conclude by weak associativity.
\EPf
We now introduce a notion of Schwartz space which will be of importance
in the sequel.
\begin{dfn}\label{SCHWARTZSPACE}
We define the \emph{Schwartz space} of $\L$ as the following function space:
$$
\CS(\L)\;=\;\{\varphi\in C^\infty(\L)\;|\;\forall X\in\CU(\fL), N\in\R: \;|\fd^N\tilde{X}(\varphi)|<\infty\}\;.
$$
\end{dfn}
Of course the Schwartz space is a Fr\'echet space for the following semi-norms:
$$
|\varphi|_{J,N}|\;:=\;|\varphi|^{-N}_J
$$
with the notation of Definition \ref{SEMINORM}. Note that $N$ is variable here, which was not the case
for the spaces of $\CB$-type. 

\noindent Note also that
$$
\CS(\L)\;=\;\cap_{N,q\in\N}\;\CB_q^{-N}(\L)
$$
and that one has an obvious $\fB$-valued analogue.

\begin{prop}
The Schwartz space $\CS(\L)$ is a two-sided ideal in $(\CB^0_\infty(\L),\star_S)$.
\end{prop}
\Pf The proof is similar to the one of Formula \ref{RAB}. 
We consider $F\in\CB^0_\infty(\L)$, which we view as an element of $\CB^0_(\L)$, and $\varphi\in\CS(\L)$. As previously, we consider basic elements $Y\in\CU_{q_1}(\fL)\otimes\CU_{q_2}(\fL)$ and $X\in\CU_q(\fL)$. We then bound the quantity:
\begin{equation}\label{PSI}
\sup_{(x_1,x_2)\in\L^2}\{\left|\tilde{Y}_{(x_1,x_2)}\CR(F\otimes\varphi)\right|_{X,N}\fd^{-N_1}_\L(x_1)\fd^{-N_2}_\L(x_2)\}
\end{equation}
where we have first observed that $\CR(F\otimes\varphi)$ belongs to $\CS(\L)$.

\noindent Again, following the same lines as before, we have that
$$
\tilde{X}_{x}\left(\tilde{Y}_{(x_1,x_2)}\CR(F\otimes\varphi)\right)\;=\;\CR
\left[(\Ad_{x_1^{-1}}X\otimes\Ad_{x_2^{-1}}X).Y\right]^\sim_{(xx_1,xx_2)}(F\otimes\varphi)\;.
$$
The module of the above expression is therefore bounded by elements of the form
$$
C\fd^q_\L(x_1)\fd^q_\L(x_2)|F|_{q+q_1}|\varphi|_{q+q_2,N'}\fd^{-N'}_\L(x_2)\fd^{-N'}_\L(x)
$$
for any $N'\in\N_0$.

\noindent In order to bound (\ref{PSI}), we therefore need to bound expressions of the form
$$
\fd^{-N_1+q}_\L(x_1)|F|_{q+q_1}|\varphi|_{q+q_2,N}\;.
$$
We therefore have that the restriction of $\CR$ to $\CB^0_p(\L)\times\CS(\L)$ is valued in
$$
\CB^{q,0}_{p-q}(\L^2,\CB^{-N'}_{q'}(\L))
$$
for all $q, q',N'\in\N_0$ ($N'$ large enough).
\EPf

\chapter{Left-invariant non-formal star-products and Drinfel'd-Iwasawa twists}

\section{Unitary representations of symmetric spaces}
\begin{dfn}
Let $(M,s)$ be a symmetric space with transvection involutive Lie group $(G,\sigma)$. A {\bf unitary representation} of $(M,s)$ is a 
triple $(\CH,\bU,\bSigma)$ where $\CH$ is a separable Hilbert space, $\bU$ is a unitary (strongly continuous) representation $\bU$ of the Lie group $G$ on $\CH$:
$$
\bU:G\to GL(\CH)\;,
$$
and where $\bSigma$ is a unitary involution:
$$
\bSigma:\CH\to\CH
$$
such that, for every $g\in G$:
$$
\bU(\sigma g)\;=\;\bSigma\,\bU(g)\,\bSigma\;.
$$
\noindent The representation $(\CH,\bOmega)$ is called {\bf irreducible} whenever the underlying group
representation $(\CH,\bU)$ is.
\end{dfn}
\begin{prop} Let $(\CH,\bU,\bSigma)$ be a unitary representation of the symmetric space $(M,s)$.
Let $K$ be isotropy subgroup in $G$ of a base point $o\in M$, and consider the $G$-homogeneous
spaces isomorphism:
$$
G/K\to M:gK\mapsto g.o\;.
$$
Then,

(i) The map:
$$
\bOmega:G\to U(\CH): g\mapsto\bOmega(g)\;:=\;\bU(g)\bSigma\bU(g^{-1})
$$
is constant each left-lateral class of $K$ in $G$. It therefore descends
to $M=G/K$ as a map again denoted by $\bOmega$:
$$
\bOmega:M\to U(\CH): g.o\mapsto\bOmega(g)\;.
$$

(ii) The map $\bOmega$
satisfies the following properties:
$$
\bOmega(x)\bOmega(y)\bOmega(x)\;=\;\bOmega(\bs_xy)
$$
for all $x,y\in M$.
\end{prop}

\subsection{Representations of group type symmetric spaces }\label{REP}

\noindent The notion of representation of a symmetric space is defined in a natural way: it is a representation of its transvection group that is compatible
with the symmetric space structure as homogeneous space.

\begin{dfn}\label{SIR} A \emph{square integrable irreducible unitary representation} of a group type symmetric space $M$ is a triple $(\CH,U,\Sigma)$ where 

\begin{enumerate}
\item[$\bullet$] $\CH$ is a separable Hilbert space,

\item[$\bullet$] $U$ is a  unitary representation of $G$ on $\CH$, 

\item[$\bullet$] 
There exists a $G$-invariant dense Fr\'echet subspace $\CH^\infty$ of $\CH$  and an  automorphism
$\Sigma$  of $\CH^\infty$ that commutes with $U(K)$.
\end{enumerate}

such that
\begin{enumerate}
\item[(i)] the following map
$$
\Omega:G/K\to \Aut(\CH^\infty):gK\mapsto U(g)\Sigma U(g^{-1}) \;.
$$
is continuous and injective.
\item[(ii)]The restriction of the representation $U$ to $\S$ is  irreducible and square integrable.
\end{enumerate}
\end{dfn}

\begin{rmk} Note that the map $\Omega$ is necessarily $G$-equivariant.
\end{rmk}

\noindent Given such a represented group type symmetric space $M$, since $\Omega$ is continuous, 
we consider the map
$$
\Omega:\CD(M)\to \CB(\CH^\infty):f\mapsto\int_\S f(x)\,\Omega(x)\,{\rm d}x
$$
where the measure on $M$ is the transportation of the left-invariant Haar measure
on $\S$ and where $\CB(\CH^\infty)$ denotes the algebra of 
continuous linear operators in $\CH^\infty$.

\begin{dfn}
The representation is called \emph{tracial} if the map $\Omega$ extends to  a unitary ($G$-equivariant)
isomorphism from $L^2(M)$ onto the algebra $\CL^2(\CH)$ of Hilbert-Schmidt operators on $\CH$:
$$
\Omega : L^2(M)\stackrel{\sim}{\longrightarrow}\CL^2(\CH)\;.
$$
\end{dfn}
\begin{rmk} In the tracial case,
the inverse $\sigma_\Omega$ of $\Omega$ is given by
$$
\sigma_\Omega(A)(x)\;=\;\Tr(A\Omega(x))
$$
understood as a weak trace distribution.
\end{rmk}

\subsection{Retract pairs}\label{RETP}

\noindent The notion of retract pair concerns a pair of group type symmetric spaces that share a common symmetry. Generally, we think of the situation where one of them is a curvature contraction of the other.
They then form a retract pair when they share this common symmetry at the representation level as well.

\noindent In the sequel, we will use the notion of retract pair to somehow invert the contraction level at the 
functional level.

\vspace{3mm}

\noindent Let $M_1=G_0/K_1$ and $M_2=G_1/K_2$ be group type symmetric spaces on the same abstract Lie group $\S$.
Let $(\CH_j,U_j,\Sigma_j)$ ($j=1,2$) be a square integrable irreducible unitary representation
of $M_j$.
\begin{lem}
Let $\phi^{(1)}_e\in\CH_1\backslash\{0\}$ and $\phi^{(2)}_e\in\CH_2\backslash\{0\}$ be both in the domains of the square-root of the Duflo-Moore operators respectively associated to both representation. Set for every $x\in\S$:
$$
\phi^{(1)}_x\;:=\;U_1(x)\phi^{(1)}_e\quad\phi^{(2)}_x\;:=\;U_2(x)\phi^{(2)}_e\;.
$$
Then, for all $\varphi\in\CH_1$ and $\psi\in\CH_2$, the function
\begin{equation}\label{INT}
\S\to\C:x\,\mapsto\,<\varphi,\phi^{(1)}_x>_1\,<\phi^{(2)}_x,\psi>_2
\end{equation}
(where $<\,,\,>_j$ denotes the inner product on $\CH_j$)
is integrable  w.r.t. a left-invariant Haar measure ${\rm d}x$ on $\S$.
\end{lem}
\Pf For square-integrability of the representations (restricted to $\S$), Cauchy-Schwartz inequality yields the assertion.
\EPf
\noindent We denote by $T$ the map from $\CH_2$ to $\CH_1$ defined by the above matrix coefficients,
which wihtin the bra-cket notations is expressed by the weak integrals:
$$
[T:\CH_2\to\CH_1]\;:=\;\int_\S|\phi^{(1)}_x>\,<\phi^{(2)}_x|\,{\rm d}x\;.
$$

\begin{prop}\label{EQUIV}
The map $T$ is either trivial either a unitary equivalence of $\S$-representations (up to a nonzero constant 
factor).
\end{prop}
\Pf
Consider the  map
$$
[T^\ast :\CH_1\to\CH_2]\;:=\;\int_\S|\phi^{(2)}_x>\,<\phi^{(1)}_x|\,{\rm d}x\;.
$$
One readily checks that $T$ and $T^\ast $ are intertwiners:
$$
U_1T\;=\;TU_2\quad\mbox{\rm and}\quad U_2T^\ast \;=\;T^\ast U_1
$$
Therefore we have
$$
 T^\ast TU_2\;=\;T^\ast U_1T\;=\;U_2 T^\ast T\;.
$$
The self-adjoint element $T^\ast T$ therefore commutes with the 
irreducible representation $U_2$. Hence Schur's lemma concludes.
\EPf
\begin{dfn}\label{RETP}
The represented  group type symmetric spaces $M_1$ and $M_2$ constitute  a \emph{retract pair}
$(M_1,M_2)$ if the restrictions of the representations $U_0$ and $U_1$ of $\S$ are unitarily equivalent.

\noindent Based on Proposition \ref{EQUIV}, which provides an explicit intertwiner, we will in this case therefore assume that 
$$\CH_1\;=\;\CH_2\;=:\;\CH\;.$$
\end{dfn}

\noindent We will see that in many cases the geometry of a group type
symmetric space is quite well encoded by a representation. 
The idea of the retract pair then, as we will see later on, is that from a ``nice" 
group type symmetric space, one can reach the geometry of a richer and more complicated  one. 
By ``reaching", we mean finding some kind of equivalence, at the operator algebraic level,
that passes from the ``nice"one to the other one.

\section{Equivariant quantization of Iwasawa factors $\S=AN$}\label{SAN}
The results exposed in this section are extracted from \cite{B1} and \cite{BG}, except for the first subsection 
which is classical. 
\subsection{Weyl's quantization}\label{WQ}
We consider the plane $\R^2\;=\;\{x=(q,p)\}$. For every $\tilde{u}\in\CD(\R)$, we set
$$
\Omega_W(x)\tilde{u}(q_0)\;=\;e^{2i/\theta(q-q_0).p}\tilde{u}(2q-q_0)
$$
where $\mu$ is a non-zero real parameter.
For every $\varphi\in\CD(\R^2)$, we consider
\begin{equation}\label{WEYL}
\Omega_W(\varphi)\tilde{u}(q_0)\;:=\;\int\varphi(x)\Omega_W(x)\tilde{u}(q_0){\rm d}x\;=\;
\frac{1}{2}\int \,e^{i/\theta(q-q_0).p}\varphi(\frac{q_0+q}{2},p)\,\tilde{u}(q)\,{\rm d}p\,{\rm d}q\;.
\end{equation}
We denote by $\CS(\R^2)$ the usual Schwartz space of rapidly decreasing smooth functions on $\R^2$.
One then has the following classical result:
\begin{thm} Let $\theta\in\R_0$. Then,

(i) Formula (\ref{WEYL}) defines a unitary isomorphism (up to constant):
$$
\Omega_W:L^2(\R^2)\to\CL^2(L^2(\R)):\varphi\mapsto\Omega_W(\varphi)\;.
$$

(ii) For all $u,v\in\CS(\R^2)$, one has the symbol composition formula:
$$
\left(\Omega^W_\theta\right)^{\ast}(\Omega_{\theta,\lambda}^W(u)\Omega_{\theta,\lambda}^W(v))(x)\;=:\;u\star_\theta^Wv(x)\;=\;
\frac{1}{\theta^2}\int_{\R^2\times\R^2}e^{i/\theta S_W(x,y,z)}u(y)\,v(z)\,{\rm d}y\,{\rm d}z
$$
where, setting $\omega\;:=\;{\rm d}p\wedge{\rm d}q$:
$$
S_W(x,y,z)\;:=\;\omega(x,y)+\omega(y,z)+\omega(z,x)\;.
$$

(iii) The product $\star^W_\theta$ closes on the Schwartz space $\CS(\R^2)$. The algebra $(\CS(\R^2),\star_\theta^W)$ is then 
Fr\'echet.
\end{thm}

\subsection{Quantization of $\S$}\label{QAN}
We now consider the affine group $\S $ as realized in Section \ref{PLP}. 
\begin{nt} Let $\bm\in\CB^N_\infty(\S,\R)$, $\bm>0$ and $\tilde{E}\bm=0$ (i.e. $\bm$ depends on the variable $a\in A$ only).

1. We denote by 
$$
\CF\varphi(a,\xi)\;:=\;\int_Ne^{-in\xi}\varphi(a,n)\,{\rm d}n
$$
the partial Fourier transform in the $N$-variable.

2. We consider the one-parameter family of diffeomorphisms of $\R^2_{(a,\xi)}$:
$$
\phi_\theta(a,\xi)\;:=\;(a,\frac{2}{\theta}\sinh(\frac{\theta}{2}\xi))\;.
$$

3. We set
$$
\bm_\theta(\xi)\;:=\;\bm(\frac{\theta}{2}\arcsinh(2\xi))\;.
$$

4. For every $\varphi\in\CD(\S)$, we define
$$
T_{\theta,\bm}\varphi\;:=\;\CF^{-1}\circ\,\bm_\theta^{-1}\circ\,(\phi^{-1}_\theta)^\star\circ\,\CF(\varphi)\;.
$$

5. We denote by $L^2_{\theta,\bm}(\S)$ the Hilbert completion of $\CD(\S)$ under the inner product:
$$
<u,v>_{\theta,\bm}\;:=\;\int(\phi^\star_\theta\bm_\theta\CF\overline{u})\,.\,(\phi^\star_\theta\bm_\theta\CF v)\,{\rm d}\xi\,{\rm d}a\;.
$$

6. For every $\lambda\in\R_0$, we set 
$$
\CH_\lambda\;:=\;L^2(A,e^{2\lambda a}\,{\rm d}a)\;.
$$

\end{nt}

\begin{lem}\label{AR} Let $r\in\R$.
Then, for all $k\geq2(1+r)$ and $N\in\Z$, the operator
$$
\bA_r\;:=\;\CF^{-1}\circ\,\left(1+\frac{\theta^2\xi^2}{4}\right)^{r}\circ\,\CF
$$
extends  from $\CD(\S)$ to a continuous linear map:
$$
\bA_r\;:\;\CB^N_{k}(\S)\to\CB^{N+2}_{k-2(1+r)}(\S).
$$
In particular, the map $\bA_r$ continuously preserves the Schwartz space $\CS(\S)$.

\end{lem}
\Pf
We have, for $\varphi\in\CD(\S)$:
$$
\bA_r\varphi(a_0,n_0)\;:=\;\int\,e^{in_0\xi}\left(1+\frac{\theta^2\xi^2}{4}\right)^{r}\,e^{-in\xi}\varphi(a_0,n)\,{\rm d}n\,{\rm d}\xi\;.
$$
Observing that, for every $m\in\N$, one has:
$$
(1-\partial_n^2)^m\frac{1}{(1+\xi^2)^{m}}\,e^{-in\xi}\;=\;e^{-in\xi}\;,
$$
an integration by parts yields:
$$
\bA_r\varphi(a_0,n_0)\;:=\;\int\,e^{in_0\xi}\chi_m(\xi)\,e^{-in\xi}\varphi_m(a_0,n)\,{\rm d}n\,{\rm d}\xi
$$
with
$$
\chi_m(\xi)\:=\;\frac{\left(1+\frac{\theta^2\xi^2}{4}\right)^{r}}{(1+\xi^2)^{m}}\quad\mbox{\rm and}\quad
\varphi_m\;:=\;((1-\tilde{E}^2)^m\varphi)\;.
$$
Since
$$
\tilde{H}\;=\;\partial_{a_0}-2n_0\partial_{n_0} \;,
$$
one has
\begin{eqnarray*}
&\tilde{H}\bA_r\varphi(a_0,n_0)&\;=\;\int\,e^{in_0\xi}\chi_m(\xi)\,e^{-in\xi}
\partial_{a_0}\varphi_m(a_0,n)\,{\rm d}n\,{\rm d}\xi-2in_0\int\,\xi\,e^{in_0\xi}\chi_m(\xi)\,e^{-in\xi}\varphi_m(a_0,n)\,{\rm d}n\,{\rm d}\xi\\
&=&
\int\,e^{in_0\xi}\chi_m(\xi)\,e^{-in\xi}
\partial_{a_0}\varphi_m(a_0,n)\,{\rm d}n\,{\rm d}\xi+2n_0\int\,e^{in_0\xi}\chi_m(\xi)\,\partial_{n}(e^{-in\xi})\varphi_m(a_0,n)\,{\rm d}n\,{\rm d}\xi\\
&=&
\int\,e^{in_0\xi}\chi_m(\xi)\,e^{-in\xi}
\partial_{a_0}\varphi_m(a_0,n)\,{\rm d}n\,{\rm d}\xi-2n_0\int\,e^{in_0\xi}\chi_m(\xi)\,(e^{-in\xi})\tilde{E}\varphi_m(a_0,n)\,{\rm d}n\,{\rm d}\xi\\
&=&
\int\,e^{in_0\xi}\chi_m(\xi)\,e^{-in\xi}
\partial_{a_0}\varphi_m(a_0,n)\,{\rm d}n\,{\rm d}\xi+2i\int\,\partial_\xi(e^{in_0\xi})\chi_m(\xi)\,e^{-in\xi}\tilde{E}\varphi_m(a_0,n)\,{\rm d}n\,{\rm d}\xi\\
&=&
\int\,e^{in_0\xi}\chi_m(\xi)\,e^{-in\xi}
\partial_{a_0}\varphi_m(a_0,n)\,{\rm d}n\,{\rm d}\xi-2i\int\,e^{in_0\xi}\chi_m'(\xi)\,e^{-in\xi}\tilde{E}\varphi_m(a_0,n)\,{\rm d}n\,{\rm d}\xi\\
&&
-2\int\,e^{in_0\xi}\chi_m(\xi)\,ne^{-in\xi}\tilde{E}\varphi_m(a_0,n)\,{\rm d}n\,{\rm d}\xi\\
&=&
\int\,e^{in_0\xi}\chi_m(\xi)\,e^{-in\xi}
\tilde{H}\varphi_m(a_0,n)\,{\rm d}n\,{\rm d}\xi-2i\int\,e^{in_0\xi}\chi_m'(\xi)\,e^{-in\xi}\tilde{E}\varphi_m(a_0,n)\,{\rm d}n\,{\rm d}\xi\;.
\end{eqnarray*}

Similarly, one has
\begin{eqnarray*}
&\tilde{E}\bA_r\varphi(a_0,n_0)&\;=\;i\int\,\xi\,e^{in_0\xi}\chi_m(\xi)\,e^{-in\xi}
\varphi_m(a_0,n)\,{\rm d}n\,{\rm d}\xi\\
&=&
-\int\,e^{in_0\xi}\chi_m(\xi)\,\partial_n(e^{-in\xi})
\varphi_m(a_0,n)\,{\rm d}n\,{\rm d}\xi\\
&=&
\int\,e^{in_0\xi}\chi_m(\xi)\,e^{-in\xi}
\tilde{E}\varphi_m(a_0,n)\,{\rm d}n\,{\rm d}\xi\;.
\end{eqnarray*}
Iterating, we conclude that for every $X\in\CU_k(\s)$, $\tilde{X}\bA_r\varphi$ is a linear 
combination of elements of the form:
$$
\Xi(a_0,n_0)\;:=\;\int\,e^{in_0\xi}\chi_{{}_{X'}}(\xi)\,e^{-in\xi}\tilde{X'}\varphi(a_0,n)\,{\rm d}n\,{\rm d}\xi
$$
with
$$
|X'|= k+2m
$$
and $\chi_{{}_{X'}}$ of the form
$$
\chi_{{}_{X'}}(\xi)\;=\;\frac{b_{{}_{X'}}(\xi)}{(1+\xi^2)^{m-r}}
$$
where $b_{{}_{X'}}\in\CB_\infty^0(\R)$.

\noindent Now, we estimate the growth in the variable $n$.  Since 
$$
(1-\partial_\xi^2)^j\frac{1}{(1+n^2)^{j}}\,e^{-in\xi}\;=\;e^{-in\xi}\;,
$$
we have for every $j\in\N$:
$$
\Xi(a_0,n_0)\;:=\;\int\,(1-\partial_\xi^2)^j\left(e^{in_0\xi}\chi_{{}_{X'}}(\xi)\right)\,\frac{1}{(1+n^2)^{j}}\,e^{-in\xi}\tilde{X'}\varphi(a_0,n)\,{\rm d}n\,{\rm d}\xi
$$
which is a linear combination of elements of the form:
$$
\Xi_\ell(a_0,n_0)\;:=\;n_0^\ell\int\, e^{in_0\xi}\chi_{{}_{X'}}(\xi)\,\frac{1}{(1+n^2)^{j}}\,e^{-in\xi}\tilde{X'}\varphi(a_0,n)\,{\rm d}n\,{\rm d}\xi
$$
where 
$$
\ell\leq 2j
$$
and where $b_{{}_{X'}}$ in $\chi_{X'}$ has possibly been redefined. One concludes observing that
$$
\fd_\S(a,n)\;\sim\;\sqrt{n^2+\cosh(4a)}\;.
$$

\EPf

\begin{lem}\label{INVOLUTIONSM}

\noindent (i) Let $\varphi\in\CH_\lambda$. The formula
$$
(U_{\theta,\lambda}(a,n)\varphi)(a_0)\;:=\;e^{-\lambda a\,+\,\frac{i}{\theta} e^{2(a-a_0)}n}\varphi(a_0-a)\;
$$
defines an irreducible unitary representation of $\S$ on $\CH_\lambda$.

\noindent (ii) Let 
$$
\Lambda(a)\;:=\;e^{-\lambda a}\;.
$$
Then the multiplication operator by $\Lambda$:
$$
m_\lambda:\varphi\mapsto\Lambda\varphi
$$
defines a unitary equivalence of $\S$-representation:
\begin{eqnarray*}
m_\lambda:\CH_0\to\CH_{\lambda}\\
m_\lambda U_{\theta,0}m_{\lambda}^{-1}\;=\;U_{\theta,\lambda}\;.
\end{eqnarray*}
\noindent (iii) The formula
$$
\Sigma_{\bm,\lambda}\varphi(a_0)\;:=\;\bm(a_0)\Lambda^2(a_0)\varphi(-a_0)\;=\;\bm(a_0) e^{-2\lambda a_0}\varphi(-a_0)
$$
defines a (possibly unbounded) invertible operator:
$$
\Sigma_{\bm,\lambda}:\CH_\lambda\to\CH_\lambda
$$
that underlies, together with $(\CH_\lambda,U_{\theta,\lambda})$, a square integrable unitary irreducible representation:
$$
\Omega_{\theta, \bm,\lambda}:\M\to U(\CH_\lambda)
$$
of the symmetric space $\M$.
\end{lem}
The following proposition describes 
\begin{prop}\label{CL2}
(i) The map $T_{\theta,\bm}^{-1}$ stabilises the Schwartz space:
$$
T_{\theta,\bm}^{-1}\:\;\CS(\R^2)\to\CS(\R^2)\;.
$$

(ii) The map 
$$
\varphi\mapsto\Omega_{\theta,\bm,\lambda}(\varphi)\;:=\;m_\Lambda \left(\Omega^W_{\theta}T_{\theta,\bm}^{-1}(\varphi)\right)m_{\Lambda^{-1}}
$$
extends from $\CD(\S)$ to an $\S$-equivariant unitary isomorphism between $L^2_{\theta,\bm}(\S)$ and the Hilbert-Schmidt 
operators on $\CH_\lambda$:
$$
\Omega_{\theta,\bm,\lambda}:L^2_{\theta,\bm}(\S)\to\CL^2(\CH_\lambda)\;.
$$

(iii) In terms of Lemma (\ref{INVOLUTIONSM}) (iii), one has for every $\varphi\in\CD(\S)$:
$$
\Omega_{\theta,\bm,\lambda}(\varphi)\;=\;\int_\S\varphi(x)\,\Omega_{\theta,\bm,\lambda}(x)\,{\rm d}x\;.
$$
\end{prop}
We now have analogous results as in the flat Heisenberg case.
\begin{prop}\label{STARM} 
Set
$$
\bm^0(a,n)\;:=\;\sqrt{\cosh(2a)}\;.
$$
Then,

(i) $$ L^2_{\theta,\bm^0}(\S)\;=\;L^2(\S)\;.$$

(ii) For every $u,v\in\CD(\S)$, the function
$$
\Omega_{\theta,\bm,\lambda}^{-1}(\Omega_{\theta,\bm,\lambda}(u)\Omega_{\theta,\bm,\lambda}(v))
$$
does not depend on $\lambda$.

(iii) Setting 
$$
u\star_{\theta, \bm}v\;:=\;\Omega_{\theta,\bm,\lambda}^{-1}(\Omega_{\theta,\bm,\lambda}(u)\Omega_{\theta,\bm,\lambda}(v))\;,
$$
one has
$$
u\star_{\theta, \bm}v(x)\;=\;\frac{1}{\theta^2}\int_{\L} A_{\bm}(x,y,z)\,e^{\frac{1}{\theta}S_\M(x,y,z)}u(y)\,v(z)\,{\rm d}y\,{\rm d}z
$$
where
$$
A_{\bm}((a_0,n_0),(a_1,n_1),(a_2,n_2))\;:=\;\bm(a_0-a_1)\frac{(\bm^0)^2(a_1-a_2)}{\bm(a_1-a_2)}\bm(a_2-a_0)
$$
and
$$
S_\M((a_0,n_0),(a_1,n_1),(a_2,n_2))\;:=\;\sinh(a_0-a_1)n_2+\sinh(a_2-a_0)n_1+\sinh(a_1-a_2)n_0\;.
$$

(iv) The product $\star_{\theta, \bm}$ closes on the Schwartz space $\CS(\S)$. The algebra $(\CS(\S),\star_{\theta, \bm})$ is then 
Fr\'echet.
\end{prop}

\begin{nt}
The case $\bm=1$ is of particular importance for us. We therefore adopt the following notations:
\begin{eqnarray*}
\Sigma_{1,\lambda}&=:&\Sigma_\lambda\\
\Omega_{\theta,1,\lambda}&=:&\Omega_{\theta,\lambda}\\
\star_{\theta,1}&=:&\star_\theta\;.
\end{eqnarray*}
\end{nt}
We have now three results concerning extensions of $\star_{\theta, \bm}$ to $\CB$-type spaces. The following proposition directely
follows from the proofreading of Theorem 2.43 and Proposition 2.47 of \cite{BG}.

\begin{prop}\label{ASBG}
Let $t,s\geq0$, $p\in\N\cup\{\infty\}$ with $p>4+2N+t+s$ ($\bm\in\CB^N_\infty(\S)$). Then, 

(i) the product $\star_{\theta, \bm}$ canonically extends form $\CD(\S)\times\CD(\S)$ to a bi-linear separately continuous 
mapping:
$$
\star_{\theta, \bm}:\CB_p^t(\S)\times \CB_p^s(\S)\to \CB_{p-4-2N-2t}^{t+s}(\S)\;.
$$
(ii) The above extension is associative in the sense that for all $F_1,F_2,F_3\in\CB^t_\infty(\S)$, 
one has 
$$
(F_1\star_{\theta, \bm}F_2)\star_{\theta, \bm}F_3\;=\;F_1\star_{\theta, \bm}(F_2\star_{\theta, \bm}F_3)
$$
in $\CB^{3t}_\infty(\S)$.
\end{prop}

\begin{lem}
The product $\star_{\theta}$ associated with $\bm=1$ corresponds to 
$$
T_{\theta}\varphi\;:=\;\CF^{-1}\circ\,(\phi^{-1}_\theta)^\star\circ\,\CF(\varphi)\;.
$$
One then has
$$
\left(T_{\theta}^{-1}\right)^\ast T_{\theta}^{-1}\;=\;\bA_{-1/2}.
$$
\end{lem}
\Pf
From \cite{B1}, one has
$$
\mbox{\rm Jac}_{\phi_\theta^{-1}}(\xi)\;=\;(1+\frac{\theta^2\xi^2}{4})^{-1/2}\;.
$$
The lemma then follows from a direct computation. \EPf
\begin{cor}\label{BS}
Setting
$$
\bS_0\;:=\;\bA_{1/4}\;,
$$
the product  law on $\CS(\S)$
$$
u\star^2_\theta v\;:=\;\bS_0^{-1}(\bS_0 u\star_\theta\bS_0 v)
$$
is tracial in the sense that
$$
\int\overline{u}\star^2_\theta v\;=\;\int\overline{u}v\;.
$$
\end{cor}
\Pf
We have
\begin{eqnarray*}
&&
\Tr(\Omega_{\theta,\lambda}(\varphi_1)^\ast\Omega_{\theta,\lambda}(\varphi_2))\;=\;\Tr(\Omega_WT_{\theta}^{-1}(\varphi_1)^\ast\Omega_WT_{\theta}^{-1}(\varphi_2))\;=\;\int T_{\theta}^{-1}(\overline{\varphi_1})\star_WT_{\theta}^{-1}(\varphi_2)\\
&=&
\int T_{\theta}^{-1}(\overline{\varphi_1})T_{\theta}^{-1}(\varphi_2)\;=\;
\int \left(T_{\theta}^{-1}\right)^\ast T_{\theta}^{-1}(\overline{\varphi_1})\varphi_2\;=\;
\int \bA_{-1/2}(\overline{\varphi_1})\varphi_2\;.
\end{eqnarray*}
Setting 
$$
\Omega^{(2)}_{\theta,\lambda}\;:=\;\Omega_{\theta,\lambda}\bS_0
$$
one then has
\begin{eqnarray*}
&&
\Tr(\Omega^{(2)}_{\theta,\lambda}(\varphi_1)^\ast\Omega^{(2)}_{\theta,\lambda}(\varphi_2))\;=\;\int\overline{\varphi_1}\star^2\varphi_2\;=\;
\Tr(\Omega_{\theta,\lambda}\bS_0(\varphi_1)^\ast\Omega_{\theta,\lambda}\bS_0(\varphi_2))\\
&=&
\int\overline{\bS_0\varphi_1}\star_1\bS_0\varphi_2\;=\;\int \bA_{-1/2}(\overline{\bS_0\varphi_1})\bS_0\varphi_2
\;=\;\int \bS_0^\ast \bA_{-1/2}\bS_0(\overline{\varphi_1})\bS_0\varphi_2\;=\;
\int \overline{\varphi_1}\varphi_2
\end{eqnarray*}
because $\bA$ hence $\bS_0$ are self-adjoint.
\EPf

\noindent At last, we have a non-Abelian version \cite{BG} of the classical Calder\`on-Vaillancourt theorem.
\begin{thm}\label{CVBG}
There exists $p_0\in\N$ such that the map $\Omega_{\theta,\bm,\lambda}$ canonically extends from
$\CD(\S)$ to a continuous mapping from $\CB^0_{p_0}(\S)$ valued in the bounded operators:
$$
\Omega_{\theta,\bm,\lambda}:\CB^0_{p_0}(\S)\to\CB(\CH_\lambda)\;.
$$
\end{thm}

\begin{prop} The mapping 
$$
\Omega_{\theta,\bm,\lambda}:\CB^0_{p_0}(\S)\to\CB(\CH_\lambda)\;.
$$
is injective.
\end{prop}
\Pf
For every $F\in\CB^0_{p_0}(\S)$ and every $\varphi\in\CS(\S)$, the operator 
$\Omega_{\theta,\bm,\lambda}(F)\Omega_{\theta,\bm,\lambda}(\varphi)$ is Hilbert-Schmidt.
If $\Omega_{\theta,\bm,\lambda}(F)=0$, one therefore has $F\star_{\theta, \bm}\varphi=0$.
 Differentiating w.r.t $\theta$ around zero then yields
$F\varphi=0$ for every $\varphi\in\CS(\S)$. Hence $F=0$.
\EPf
\noindent Combining these two last results, one gets
\begin{prop}\label{CSTARBG}\cite{BG} Denoting by $||\;.\;||_\lambda$ the operator norm on $\CH_\lambda$,
the map
$$
||\;.\;||_{\theta,\bm,\lambda}:\CB^0_{p_0}(\S)\to\R:F\mapsto||\Omega_{\theta,\bm,\lambda}(F)||_\lambda
$$
restricts to a $C^\star$-norm on $(\CB^0_\infty(\S),\star_{\theta, \bm})$.
\end{prop}
In fact, in terms of the regular action one has:
\begin{thm}
The regular action 
$$
\CB^0_\infty(\S)\times\CS(\S)\to\CS(\S):(F,\varphi)\mapsto F\star_{\theta, \bm}\varphi\;=:\;L_{\star_{\theta, \bm}}(F)\varphi
$$
extends to a bounded action of $ \CB^0_\infty(\S)$ on $L^2_{\theta,\bm}(\S)$, and one has
$$
||L_{\star_{\theta, \bm}}(F)||_{\theta,\bm}\;=\;||\Omega_{\theta,\bm,\lambda}(F)||_\lambda
$$
where $||\;.\;||_{\theta,\bm}$ denotes the operator norm of operators on $L^2_{\theta,\bm}(\S)$.
In particular, as expected, $||\;.\;||_{\theta,\bm,\lambda}$ does not depend on $\lambda$.
\end{thm}
\chapter{Quantization of the hyperbolic plane and Drinfel'd twists}\label{HYPSCHW}
We start from the Moyal $\star$-product on the plane $\R^2=\{(a,n)\}$:
\begin{eqnarray}\label{asympt}
 u \star^0_\hbar v &:=& u.v + i\frac \hbar 2\{ u,v \} +
 \underset{k=2}{\overset{\infty}\sum}\frac{(i\,\hbar/2)^k}{k!} \Omega^{i_1
 j_1}\cdots \Omega^{i_k j_k} \partial_{i_1\cdots i_k}u \partial_{j_1\cdots
 j_k}v \quad,\\
 \label{MoyalAsympt}
& = &\sum_{k=0}^{\infty} \, \frac{(i\,\hbar/2)^k}{k!}
\sum_{p=0}^{k}(-1)^{p} \,
\frac{k!}{p!(k-p)!}\, \partial_a^{k-p} \partial_n^p u \,
\partial_a^{p} \partial_n^{k-p} v  \quad ,
\end{eqnarray}
where \begin{equation}\label{PoissonBrackAN}
 \Omega_{ik}\Omega^{kj}=-\delta_i^k\qquad,\qquad \{ u,v \}=\Omega^{ij}\partial_iu\partial_jv: = \partial_a u \partial_n v - \partial_n u \partial_a
  v \quad.
  \end{equation}
  
\noindent It turns out that in coordinates $(a,n)$ the Moyal product (\ref{asympt}) is $\mathfrak{sl}_2(\R)$-covariant
with respect to the hyperbolic action of $G:=SL_2(\R)$ on the hyperbolic plane $\D=\S$ \cite{Bithese}. Precisely, presenting the Lie algebra $\g:=\mathfrak{sl}_2(\R)$ as generated over $\R$ by $H, E$ and $F$ satisfying:
\begin{equation}
 [H,E] = 2E \quad,\quad [H,F] = -2F \quad,\quad [E,F] = H,\label{HEF}
 \end{equation}
the moment map associated with the action of $G$ on $\D=G/K=\S$ reads:
\begin{equation}\label{moments}
\lambda_H=\sqrt{\kappa}\,n\,;\,\lambda_E=\frac{\sqrt{\kappa}}{2}e^{-2a}\,;\,\lambda_F=-\frac{\sqrt{\kappa}}{2}e^{2a}(1+n^2)\;.
\end{equation}
The associated fundamental vector fields are given by:
\begin{equation}\label{fundamental}
 H^\star = -\partial_a; \, E^\star = -e^{-2a}\partial_n;  \,
 F^\star= e^{2a}(n \partial_a - (\kappa+n^2) \partial_n)\;.
 \end{equation}
At last, the representation of $\g$ by derivations of $\star^0_\hbar$ admits the expression:
\begin{eqnarray}
\rho_{\hbar } (H) &=& -\partial_a \\
\rho_{\hbar }(E) &=& -\frac{e^{-2a}}{ \hbar } \sin ( \hbar \partial_n) \\
\rho_{\hbar }(F) &=&
 e^{2a}\left(\frac{\hbar  }4\sin ( \hbar \partial_n)\partial_a^2+
  n\cos ( \hbar \partial_n)\partial_a-(\kappa+n^2)\frac{\sin ( \hbar \partial_n)}\hbar
\right)\;.\label{lambdaF} 
\end{eqnarray}
We now consider the partial Fourier transform in the $n$-variable
\begin{equation}
{\cal F}(\varphi)(a,\zeta) = \int_{-\infty}^\infty e^{-i \zeta n} \varphi(a, n) {\rm d}n \quad,
\end{equation}
and denote by $\tilde{\S}:=\{(a,\zeta)\}$ the space where the Fourier transformed $\CF(g)$ is defined on.
\noindent Defining the following one-parameter family of diffeomorphisms of $\tilde{\S}$:
\begin{equation}
\phi_\hbar  (a,\zeta) = \left(a,\frac{\sinh (\hbar\, \zeta)}{\hbar}\right) \quad ,
\end{equation}
and denoting by $\tilde{\CS}$ (resp. $\CS$) the space of Schwartz test functions $\CS(\tilde{\S})$
(resp. $\CS({\S})$) on $\tilde{\S}=\{(a,\zeta)\}$ (resp. on $\CS(\S)$ of $\S=\{(a,n)\}$), one observes the following inclusions:
\begin{equation}
\phi_\hbar^\star\,\tilde{\CS}\,\subset\,\tilde{\CS}\;\mbox{\rm and }\,\tilde{\CS}\,\subset\,(\phi_\hbar^{-1})^\star\tilde{\CS}\,\subset\,\tilde{\CS}'\;.
\end{equation}
Therefore, every data of (reasonable) one parameter smooth family of invertible functions $\CP_\hbar=\CP_\hbar(b)$ yields an operator on the Schwartz space 
$\CS(\S)$ of $\S=\{(a,n)\}$:
\begin{equation}
\bT^{-1}:\CS(\S)\longrightarrow\CS(\S)
\end{equation}
defined as 
\begin{equation}\label{OpT}
\bT^{-1} \varphi(a_0,n_0)= \frac{1}{2\,\pi } \int \, e^{i\zeta n_0}\, \CP_\hbar(
\zeta)\, e^{\frac{-i}{\hbar}n \sinh(\hbar \zeta )} \, \varphi(a_0,n) \, {\rm d}n\,
{\rm d}\zeta\;.
\end{equation}
More generally, denoting by $\CM_f$ the pointwise multiplication operator by $f$, the  operator:
\begin{equation}
\bT:\CS(\S)\longrightarrow\CS'(\S)
\end{equation}
defined as
\begin{equation}\label{equivmoyal}
\bT\;:=\;\CF^{-1}\,\circ\,(\phi_\hbar^{-1})^\star\,\circ\CM_{\frac{1}{\CP}}\,\circ\CF
\end{equation}
is a left-inverse of $T^{-1}$.
Intertwining Moyal-Weyl's product by $\bT$ yields the above  product as $\star_{\hbar,{\cal P}} =\bT(\star^0_\hbar)$. Observe that the latter closes on the range space $\CE_\hbar:=\bT\CS$
yielding a (non-formal) one parameter family of associative  function algebras: $(\CE_\hbar,\star_{\hbar,{\cal P}} )$. Asymptotic expansions of these non-formal products produce
genuine Poincar\'e-invariant formal star products on $\M$ \cite{B1}.

\vspace{2mm}

\noindent For generic $\cal P$ the above oscillatory integral product defines, for all real {\sl value} of $\hbar$, an associative product law on some function space (as opposed to formal power series space) on $\M$. The most remarkable case being probably the one where $\cal P$ is pure phase. In the latter case, the above
product formula extends (when $\hbar\neq0$) to the space $L^2(\M)$ of square integrable functions as a Poincar\'e invariant Hilbert associative algebra. The star product there appears to be {\sl strongly closed}: for all
$u$ and $v$ in $L^2(\M)$, $u\,\star_{\hbar,{\cal P}} \,v$ belongs to $L^1(\M)$, and one has:
\begin{equation*}
\int_\M u\,\star_{\hbar,{\cal P}} \,v\;=\;\int_\M u\,v\;.
\end{equation*}

\noindent Combining \cite{BG} page 111 and \cite{BDS} Section 3, one finds that the relation with the star-product defined in \ref{STARM} is given by
$$
A_{\mbox{can}}(a_{1},a_{2},a_{3})\,\frac{\CP(a_{1}-a_{2})\,\CP(a_{3}-a_{1})}{\CP(a_{2}-a_{3})}\;=\;\bm(a_{1}-a_{2})\frac{\bm_{0}^{2}(a_{2}-a_{3})}{\bm(a_{2}-a_{3})}\,\bm(a_{3}-a_{1})\;.
$$

\noindent Now we consider the following derivation $D_0$ of $\star:=\bT(\star^0)$:
the operator 
\begin{equation}\label{DERIV}
D_0\;:=\;\bT\circ\rho_\hbar(F)\circ \bT^{-1}\;.
\end{equation}
The particular choice of 
\begin{equation}
\CP(\zeta)\;\equiv1
\end{equation}
yields the following 
expression for $i\,\Box:=\CF\,\circ\,D_0\,\circ\,\CF^{-1}$ that appears to be a second order differential operator:
\begin{eqnarray}\label{box1}
 \Box_{(a,\zeta)} = e^{2a} \left[   \frac{\hbar^2}4\,\zeta\,
\partial_a^2 + \zeta\, (1+\hbar^2\, \zeta^2) \partial_\zeta^2 + (1+\hbar ^2
\, \zeta^2) \partial_a \partial_\zeta + \hbar^2 \,\zeta\, \partial_a  +  (2 + 3 \hbar^2 \, \zeta^2)
\partial_\zeta - \zeta\,(k-\hbar^2) \right] \;.\label{Box0}
\end{eqnarray}
\begin{prop}
\noindent Let $\star:=\star_{\hbar,\CP}$ be a (formal) invariant star-product on $\M$ whose we denote by $\Der(\star)$ the algebra of derivations.
Under the condition
$$
\CP'(0)\;=\;0\;,
$$
\noindent (i) the classical moment 
$$
D^1:\g\to\Der(\star):X\mapsto \frac{1}{i\hbar}[\lambda_X\;,\;.\;]_\star
$$
is a Lie map.

\noindent (ii) The restriction of $D^1$ to $\s$ coincides with the geometric action by fundamental vector fields:
$$
D^1_X\;=\;X^\star\quad\forall X\in\s\;.
$$
\noindent (iii) Every Lie map
$$
D:\g\to\Der(\star): X\mapsto D_X
$$
that extends $D^1|_{\s}$ is of the form
$$
D_F\;=\;D^1_F\;+\;c_\hbar\,[e^{2a}\;,\;.\;]_\star
$$
where $c_\hbar\in\C[\hbar^{-1},\hbar]]$ is a formal constant.
\end{prop}

\noindent \Pf We  seek for a Lie algebra homomorphism
$$
D:\g\to\Der(\star)
$$
that extends the mapping $\s\to\Gamma^\infty(T\M):X\mapsto X^\star$.

\noindent Since $\star$ is formally equivalent to Moyal's star product, every of its derivation must be 
interior. In particular, there exists a linear map
$$
\Lambda:\g\to C^\infty(\M)[\hbar^{-1},\hbar]]:X\mapsto\Lambda_X
$$
such that
$$
D_X\varphi\;=:\;[\Lambda_X,\varphi]_\star\quad\quad(\forall\varphi\in C^\infty(\M))\;.
$$
Denoting 
$$
\CA_\hbar\;:=\;(C^\infty(\M)[\hbar^{-1},\hbar]],\star)\;,
$$
Jacobi's identity implies that the homomorphism defect
$$
\Gamma_\hbar:\bigwedge^2\g\to\CA_\hbar:X\wedge Y\mapsto\Lambda_{[X,Y]}-[\Lambda_X,\Lambda_Y]_\star
$$
is actually valued in the center of $\CA_\hbar$ i.e. in  the formal constants $\C[\hbar^{-1},\hbar]]$. Observing that it is also a Chevalley two-cocycle,
it must be exact for $\g$ is simple:
$$
\Gamma_\hbar\;=:\;\delta\gamma_\hbar\quad\quad\gamma_\hbar\in\g^\star\otimes\C[\hbar^{-1},\hbar]]\;.
$$
Possibly modifying it by subtracting $\gamma_\hbar$, we may therefore assume that the map $\Lambda:\g\to\CA_\hbar$ is Lie. 

\noindent Using the condition $\CP'(0)=0$,  the expression (\ref{equivmoyal})  yields
$$
\bT\lambda_X\;=\;\lambda_X
$$
for every $X\in\s$.
The map 
$$
\Lambda^0:\s\to\CA_\hbar:X\mapsto\frac{1}{i\hbar}\lambda_X\;=:\;\Lambda^0_X
$$
is therefore Lie.
Moreover, by a similar computation as for $T(n)=n$, we observe that there is a (generally non-zero) formal constant
$$\kappa_\hbar\;:=\;\frac{{\rm d}^2}{{\rm d}\xi^2}\left(\frac{\cosh(\hbar\xi)}{\CP(\xi)}\right)\left|_{\xi=0}\right.$$ such that
$$
T(n^2)\;=\;\kappa_\hbar + n^2\;.
$$
Since the operators $T$ and $T^{-1}$ commute with every multiplication operator $\CM_f$ with $\partial_nf\equiv0$, 
we get
$$
T\lambda_F\;=\;\lambda_F-\frac{\sqrt{\kappa}}{2}\kappa_\hbar\,e^{2a}\;.
$$
Setting
$$
\Lambda^0_F\;:=\;\frac{1}{i\hbar}\left(\lambda_F-\frac{\sqrt{\kappa}}{2}\kappa_\hbar\,e^{2a}\right)\;,
$$
we therefore get the Lie maps
$$
\Lambda^0:\g\to\CA_\hbar
$$
and 
$$
D^0:\g\to\Der(\star):X\mapsto D^0_X\;:=\;[\Lambda^0_X\;,\;.\;]_\star\;.
$$
Given another extension $D$, we consider the defect
$$
D_F-D^0_F\;=\;[\gamma_F\;,\;.\;]_\star\quad\mbox{\rm with}\quad\gamma_F\;:=\;\Lambda_F-\Lambda^0_F\;.
$$
Since for every $X\in\s$, we have $D_X\;=\;D^0_X\;=\;X^\star$, we observe:
$$
D_X\gamma_F\;=\;[\Lambda_X,\Lambda_F]_\star-[\Lambda^0_X,\Lambda^0_F]_\star\;=\;\Lambda_{[X,F]}-\Lambda^0_{[X,F]}\;=\;X^\star\gamma_F\;.
$$
Hence
$ E^\star\gamma_F\;=\;0\quad\mbox{\rm and}\quad H^\star\gamma_F\;=\;-2\gamma_F\;.$ \EPf
\begin{rmk}
Note that varying $D$ only corresponds to redefining the curvature element $\kappa$.
\end{rmk}
\begin{prop}
Let $U$ be a $\S$-equivariant equivalence between $\star$ and a $G$-invariant $\star$-product $\sharp$.
Then,

\noindent (i) There exists a formal distribution $v\in\CD'(\S)[[\hbar]]$ whose associated left-invariant convolution operator on $\S$ coincides with $U^{-1}$:
$$
U^{-1}\;=\;\int_\S v(y)\,R^\star_y\,{\rm d}y\;.
$$
\noindent (ii) The Fourier transformed convolution kernel $\hat{v}\;:=\;\CF v$ is a weak solution of second order PDE:
$$
\Box(P\hat{v})\;=\;\zeta(P\hat{v})
$$
where the multiplicative factor $P$ is defined by
$$
P\;:=\;\left(\phi_\hbar^{-1}\right)^\star\CP\;.
$$
\noindent (iii) The $\star$-products $\star_{\hbar,\CP}$ for which there exists a $\zeta$-independent solution $\hat{v}$ correspond, up to a formal constant multiple, either
to
$$
\CP\;=\;1\quad\mbox{\rm or}\quad\CP\;=\;\cosh(\hbar\zeta)\;.
$$
\end{prop}

\section{Retracting the hyperbolic plane: non formal context}\label{RHP}
In this section, we will first prove that the pair $(\M\;:=\;(\S,s^{(0)},\omega)\;,\;\D)$ constitute a retract pair on the  group manifold $\S=AN$. We will then describe how both symmetric space structures lift at the representation
level as ``phase" symmetries. At last we will express an operator that realizes a passage (a curvature deformation process) from one phase symmetry 
to the other and vice versa. This last operator will define the intertwiner $\fW$ (see the Introduction) that will 
transport the star product on $\M$ onto a star product on $\D$ with full $\SLdr$-symmetry.
\begin{nt}
For every $\lambda\in\R_0$, we set 
$$
\CH_\lambda\;:=\;L^2(A,e^{2\lambda a}\,{\rm d}a)\;.
$$
\end{nt}
From direct computations, we obtain
\begin{prop}\label{INVOLUTIONSM}
\noindent (i) Let $\theta$ be a positive real. Consider the upper half-plane 
$$
\D\;:=\;\{\,z\;:=\;p+iq\,\in\C\;:\;q>0\,\}\;.
$$
and let 
$$\CD_{{\left(\frac{\theta+2}{\theta}\right)}}\,:=\,\{h\in\mbox{\rm hol}(\D)\;|\;\int_{\D}\frac{1}{q^{2-2{\left(\frac{\theta+2}{\theta}\right)}}}|h(z)|^2{\rm d}q{\rm d}p<\infty\}$$
where $\mbox{\rm hol}(\D)$ stands for the space of holomorphic functions on $\D$.

\noindent Denote by $V_{\theta}$ the usual projective holomorphic discrete series of $\SLdr$ on $\CD_{\left(\frac{\theta+2}{\theta}\right)}$.

\noindent Then,

(i) The expression 
$$
T_{\theta,\lambda}\tilde{u}(z_0)\;:=\;\int_\R e^{(\lambda-2{\left(\frac{\theta+2}{\theta}\right)}) q_0}e^{\frac{i}{\theta}e^{-2q_0}z_0}\,\tilde{u}(q_0)\,{\rm d}q_0
$$
defines a unitary isomorphism
$$
T_{\theta,\lambda}\;:\;\CH_\lambda\to\CD_{\left(\frac{\theta+2}{\theta}\right)}
$$
that intertwines $V_{\theta}|_\S$ and $U_{\theta,\lambda}|_\S$.

\noindent We denote by ${\bf U}_{\theta, \lambda}$ the corresponding representation of $\SLdr$
on $\CH_\lambda$:
$$
{\bf U}_{\theta, \lambda}\;:=\;T_{\theta,\lambda}^{-1}V_{\theta} T_{\theta,\lambda}\;.
$$

(ii) Consider $J_\mu$ the \emph{first type Bessel function} with parameter $\mu$. Then
the formula:
$$
{\vartheta}_{\theta,\lambda}\varphi(a_0)\;:=\;\frac{4}{\theta}e^{-2\lambda a_0+i\pi{\left(\frac{\theta+2}{\theta}\right)}}\int_{-\infty}^\infty
J_{2{\left(\frac{\theta+2}{\theta}\right)}-1}\left(\frac{2}{\theta}e^{-2a}\right)\,e^{2(\lambda-1)a}\varphi(2a-a_0)\,{\rm d}a
$$
defines a unitary involution:
$$
{\vartheta}_{\theta,\lambda}:\CH_\lambda\to\CH_\lambda\;.
$$

(iii) The triple $(\CH_\lambda,{\bf U}_{\theta, \lambda},{\vartheta}_{\theta,\lambda})$ is a square integrable unitary irreducible representation of the hyperbolic plane.

\noindent We denote by
$$
{\bf \Omega}_{\theta, \lambda}: \D\to \Aut(\CH_\lambda)
$$
the corresponding mapping (see Definition \ref{SIR}).
 
\end{prop}
\Pf 
The restriction to $\S$ of the holomorphic discrete series reads:
$$V_{\theta}(a_0,n_0)h(z)\;:=\;e^{-2{\left(\frac{\theta+2}{\theta}\right)} a_0}h((a_0,n_0)^{-1}.z)$$
where 
$$(a_0,n_0).z\;:=\;e^{2a_0}(n_0+z)\;.$$
We will now write an explicit expression for an intertwiner $T$ from $(\CH_0,U_{{\theta},0})$ ($\theta\neq0$) to $(\CD_{\left(\frac{\theta+2}{\theta}\right)},V_{{\left(\frac{\theta+2}{\theta}\right)}})$ of the form (cf. Proposition \ref{EQUIV}):
$$
T\tilde{u}\;:=\;\int_\S\,<U_{{\theta},0}(g)\tilde{u}^0,\tilde{u}>\,V_{\theta}(g)h^0\,{\rm d}g
$$
where $\tilde{u}^0\in\CH_0$ and $h^0\in\CD_{\left(\frac{\theta+2}{\theta}\right)}$ are fixed non-zero elements in the domain of the Duflo-Moore operator. This amounts to:
$$
T\tilde{u}(z_0)=\int\,e^{-2{\left(\frac{\theta+2}{\theta}\right)} a+\frac{i}{\theta} e^{2(a-q_0)}n}\,h^0(e^{-2a}z_0-n)\,\overline{\tilde{u}^0}(q_0-a)\tilde{u}(q_0)\,{\rm d}q_0{\rm d}a{\rm d}n\;.
$$
Choosing
$$
h^0(z):=(i+z)^{-2{\left(\frac{\theta+2}{\theta}\right)}}\;,
$$
the formulae 
\begin{equation}\label{Gr}
\underline{\CF}(u)(\tau):=\int_\R u(t)\,e^{-it\tau}{\rm d}t\quad \mbox{\rm and}\quad\underline{\CF}^{-1}\left((\alpha+i\tau)^{-k}\right)(t)=\frac{t^{k-1}}{(k-1)!}e^{-\alpha t}H(t)
\end{equation}
($H=H(t)$ denotes the characteristic function of $]0,\infty[$ and $\Re(\alpha)>0$) lead to the following
expression:
$$
T\tilde{u}(z_0)=\frac{(-i/\theta)^{2{\left(\frac{\theta+2}{\theta}\right)}-1}}{2\pi(\Gamma\left(2{\left(\frac{\theta+2}{\theta}\right)}-1\right)}\int\,e^{-2{\left(\frac{\theta+2}{\theta}\right)} a}\,e^{-\frac{1}{\theta}(1-iz_0e^{-2a})e^{2(a-q_0)}}\,e^{2(2{\left(\frac{\theta+2}{\theta}\right)}-1)(a-q_0)}\,\overline{\tilde{u}^0}(q_0-a)\tilde{u}(q_0)\,{\rm d}q_0{\rm d}a
$$
or
$$
T\tilde{u}(z_0)=\frac{(-i/\theta)^{2{\left(\frac{\theta+2}{\theta}\right)}-1}}{2\pi(\Gamma\left(2{\left(\frac{\theta+2}{\theta}\right)}-1\right)}\left(\int_\R\,e^{2(1-{\left(\frac{\theta+2}{\theta}\right)}) q-\frac{1}{\theta}e^{-2q}}\,\overline{\tilde{u}^0}(q)\,{\rm d}q\right)\,\int_\R e^{-2{\left(\frac{\theta+2}{\theta}\right)} q_0}e^{+\frac{i}{\theta}e^{-2q_0}z_0}\,\tilde{u}(q_0)\,{\rm d}q_0\;.
$$
In terms of the Laplace transform 
$$
\CL[f](w)\;:=\;\int_0^\infty\,e^{-rw}\,f(r)\,{\rm d}r\quad(\Re(w)>0)\;,
$$
we have:
$$
T\tilde{u}(z_0)\;=\;C_{\theta,\tilde{u}^0}\,\CL\left[r^{{\left(\frac{\theta+2}{\theta}\right)}-1}\tilde{v}(r)\right](-\frac{i}{\theta}z_0)
$$
where 
$$\tilde{v}(e^{-2q_0}):=\tilde{u}(q_0)\quad \mbox{\rm and}\quad C_{\theta,\tilde{u}^0}:=\frac{(-i/\theta)^{2{\left(\frac{\theta+2}{\theta}\right)}-1}}{4\pi(\Gamma\left(2{\left(\frac{\theta+2}{\theta}\right)}-1\right)}\int_\R\,e^{2(1-{\left(\frac{\theta+2}{\theta}\right)}) q-\frac{1}{\theta}e^{-2q}}\,\overline{\tilde{u}^0}(q)\,{\rm d}q\;.$$
The symmetry $s_K\;=\;-1/z$ centered at $i$ naturally ``quantizes" or lifts to the state space $\CD_{\left(\frac{\theta+2}{\theta}\right)}$ as (see also \cite{BG}):
$$\vartheta(h)(z)\;:=\;i^{2{\left(\frac{\theta+2}{\theta}\right)}}z^{-2{\left(\frac{\theta+2}{\theta}\right)}}h(s_K(z))\;.$$
Using the intertwiner $T$, we can therefore transport it at the level of $\CH_0$:
$${\vartheta}_{\theta,0}(\tilde{u})\,:=\,T^{-1}\vartheta\left(T(\tilde{u})\right)\;.$$
In order to do so, we start by observing that setting
$$
h_{r_0}(z)\;:=\;e^{\frac{i}{\theta}r_0z}\quad(z\in\D)\;,
$$
we have
$$
\vartheta(h_{r_0})(z)\;:=\;i^{2{\left(\frac{\theta+2}{\theta}\right)}}z^{-2{\left(\frac{\theta+2}{\theta}\right)}}e^{\frac{-ir_0}{\theta z}}\;,
$$
which we view as a function $\eta_{r_0}$ on the domain $\Re(w)>0$:
$$
\eta_{r_0}(w)\;:=\;i^{2{\left(\frac{\theta+2}{\theta}\right)}}(-\frac{i}{\theta})^{2{\left(\frac{\theta+2}{\theta}\right)}}w^{-2{\left(\frac{\theta+2}{\theta}\right)}}e^{-\frac{r_0}{\theta ^2w}}\quad\mbox{\rm with}\quad w\;:=\;\frac{-i}{\theta}z\;.
$$
From the formula $$\frac{1}{2\pi i}\int_{\epsilon+i\R}e^{Z-\frac{\zeta^2}{4Z}}Z^{-\mu-1}\,{\rm d}Z=\left(\frac{2}{\zeta}\right)^\mu J_{\mu}(\zeta)\quad(\epsilon, \mu,\zeta>0)\;,$$ we get
$$\CL^{-1}(\eta_{r_0})(s)\;=\;\frac{1}{\theta}\left(\frac{s}{r_0}\right)^{-\frac{-2{\left(\frac{\theta+2}{\theta}\right)}+1}{2}}J_{-(1-2{\left(\frac{\theta+2}{\theta}\right)})}\left(\frac{2}{\theta}\sqrt{r_0s}\right)\quad(s>0)\;.$$
We now express $$T\;=\;D_{\frac{-i}{\theta}}\circ\CL\circ\CM_{Cr^{{\left(\frac{\theta+2}{\theta}\right)}-1}}\circ\psi^{-1\star}\quad\mbox{\rm and}\quad T^{-1}\;=\;\psi^\star\CM_{\frac{1}{C}r^{1-{\left(\frac{\theta+2}{\theta}\right)}}}\CL^{-1}D_{i\theta}\quad (C\;:=\;C_{\theta,\tilde{u}^0})$$ where $D$
stands for the variable dilation, $\CM$ for functional multiplication and $\psi(q)\;:=\;e^{-2q}$. 
Observing that $$T\tilde{u}(z)\;=\;C_{\theta,\tilde{u}^0}\,\int_0^\infty r_0^{{\left(\frac{\theta+2}{\theta}\right)}-1}\tilde{u}(\psi^{-1}(r_0))h_{r_0}(z){\rm d}r_0\;,$$ a tedious
computation yields:
$$
{\vartheta}_{\theta,0}(\tilde{u})(q)\;=\;\int_0^\infty J_{2{\left(\frac{\theta+2}{\theta}\right)}-1}(\sigma)\,\tilde{u}\left(\ln\left(\frac{2}{\theta \sigma}\right)\,-\,q\right)\,{\rm d}\sigma\;.
$$
 The case $\lambda\neq0$ is obtained from the above one by considering the unitary equivalence mentioned in  (ii) of Lemma \ref{INVOLUTIONSM}.
 \EPf
\noindent As corollary of the proof of the above proposition, we have

\begin{lem}
Explicitly, one has for every $\varphi\in\CH_\lambda$:
$$
{\vartheta}_{\theta,\lambda}\Sigma_\lambda\varphi(a_0)\;=\;\frac{4e^{i\pi{\left(\frac{\theta+2}{\theta}\right)}}}{\theta}\int_{-\infty}^\infty
J_{2{\left(\frac{\theta+2}{\theta}\right)}-1}\left(\frac{2}{\theta}e^{-2a}\right)\,e^{2(\lambda-1)a}
\left(U_{\theta,\lambda}((2a,0))\varphi\right)(a_0)\,{\rm d}a\;.
$$
Formally, setting
$$
\CJ_{\theta,\lambda}(a)\;:=\;\frac{4e^{i\pi{\left(\frac{\theta+2}{\theta}\right)}}}{\theta}\,J_{2{\left(\frac{\theta+2}{\theta}\right)}-1}\left(\frac{2}{\theta}e^{-2a}\right)\,e^{2(\lambda-1)a}
$$
and
$$
(a,0)^2\;=:\;a^2\,\in\,A
$$
one writes
$$
{\vartheta}_{\theta,\lambda}\Sigma_\lambda\;=\;\int_A
\CJ_{\theta,\lambda}(a)\,U_{\theta,\lambda}(a^2)\,{\rm d}a\;.
$$
\end{lem}
We will now give a rigorous meaning to the above formula and extend the notion of retract 
to any, sufficiently regular, action of the group $A\simeq\R$ on a Fr\'echet space.

\noindent We let $\fF$ be a Fr\'echet space equipped with an increasing family of semi-norms 
$\{|\;|\}_{j\in\N}$. We consider a sub-isometric strongly continuous representation $\rho$ of $A$ on
$\fF$ and equip the space $\fF_\infty$ of smooth vectors with the usual refined topology associated 
with the action that turns it into a smooth Fr\'echet $A$-module for the restriction of $\rho$ to $\fF_\infty$.
Our aim in a first place, is to define the following ``improper integral" which to every $v\in\fF_\infty$, associates the quantity (in $\fF_\infty$):
$$
\Xi^{(\rho)}_{\theta,\lambda}(v)\;:=\;\int_A
\CJ_{\theta,\lambda}(a)\,\rho(a)v\,{\rm d}a\;.
$$
For this, we observe the following proposition whose proof is postponed to the proof of Proposition \ref{LAMBDA}.
\begin{prop}\label{LAMBDA0}
Let $s\in\R$ and ${\left(\frac{\theta+2}{\theta}\right)}>0$. Then there exists $\lambda_{\theta,s}$ such that the element
$$
\CD(A,\fF_\infty)\to\fF_\infty:\phi\mapsto
\int_A
\CJ_{\theta,\lambda_{\theta,s}}(a)\,\phi(a)\,{\rm d}a
$$
continuously extends from  $\CD(A,\fF_\infty)$ to $\CB^s_\infty(A,\fF_\infty)$:
$$
\widetilde{\int_A}\CJ_{\theta,\lambda_{\theta,s}}\;:\;\CB^s_\infty(A,\fF_\infty)\to\fF_\infty\;.$$
\end{prop}
\begin{cor} Let $v\in\fF_\infty$, then the element
$$
A\to\fF_\infty:a \mapsto\rho(a)v
$$
belongs to $\CB^0_\infty(A,\fF_\infty)$. 
In particular, Proposition \ref{LAMBDA} implies that the integral
$$
\Xi^{(\rho)}_{\theta,\lambda}(v)\;:=\;\int_A
\CJ_{\theta,\lambda}(a)\,\rho(a)v\,{\rm d}a\;.
$$
is a well-defined element of $\fF_\infty$.
\end{cor}
\Pf
The first assertion is obvious. For the second one, we use the fact that the closure of $\CD(\S,\fF_\infty)$ in $\CB^\epsilon_\infty(A,\fF_\infty)$ 
contains $\CB^0_\infty(A,\fF_\infty)$ (c.f. \cite{BG}).
\EPf

\begin{prop}\label{OMEGABAR}
For all $\lambda\in\R$ and $x\in\S$, one has
$$
{\bf \Omega}_{\theta,\lambda}(x)\;=\;\int_{-\infty}^\infty\CJ_{\theta,\lambda}(a)\,\Omega_{\theta,\lambda}(xa)\,{\rm d}a
$$
when acting in $\CH_\lambda^\infty$.

\noindent Also, for every $\phi\in\CD(\S)$, one has 
$$
{\bf \Omega}_{\theta,\lambda}(\phi)\;=\;\Omega_{\theta,\lambda}\left(\Xi^{(R^\vee)}_{\theta,\lambda+1}(\phi)\right)
$$
where $R^\vee_a\phi(x)\;:=\;\phi(xa^{-1})$\;.
\end{prop}
\Pf We will use the fact that
$$
{\vartheta}_{\theta,\lambda} U_{\theta,\lambda}(a){\vartheta}_{\theta,\lambda}\;=\;\Sigma_\lambda U_{\theta,\lambda}(a)\Sigma_{\lambda}\;=\;U_{\theta,\lambda}(a^{-1})\;.
$$
One has
\begin{eqnarray*}
{\bf \Omega}_{\theta,\lambda}(x)&=&U_{\theta,\lambda}(x){\vartheta}_{\theta,\lambda} U_{\theta,\lambda}(x^{-1})\\
&=& U_{\theta,\lambda}(x){\vartheta}_{\theta,\lambda}\Sigma_\lambda^2 U_{\theta,\lambda}(x^{-1})\\
&=&U_{\theta,\lambda}(x)\,\left(\int_A\CJ_{\theta,\lambda}(a)\,U_{\theta,\lambda}(a^2){\rm d}a\right)\,\Sigma_\lambda U_{\theta,\lambda}(x^{-1})\\
&=&\,\int_A\CJ_{\theta,\lambda}(a)\,U_{\theta,\lambda}(xa)\Sigma^2_\lambda U_{\theta,\lambda}(a)
\Sigma_\lambda U_{\theta,\lambda}(x^{-1})\,{\rm d}a\\
&=&\,\int_A\CJ_{\theta,\lambda}(a)\,U_{\theta,\lambda}(xa)\Sigma_\lambda
U_{\theta,\lambda}(a^{-1}x^{-1})\,{\rm d}a\;,
\end{eqnarray*}
which proves the first assertion noticing that $\Sigma_\lambda$ acts as a continuous linear operator on $\CH^\infty_\lambda$.

\noindent The second assertion is a direct application of Lemma \ref{OMEGABAR}:
 \begin{eqnarray*}
 {\bf \Omega}_{\theta,\lambda}(\phi)&=&\int_\S\phi(y){\bf \Omega}_{\theta,\lambda}(y)\,{\rm d}y\\
 &=&\int_{\S\times A}\phi(y)\CJ_{\theta,\lambda}(a)\,\Omega_{\theta,\lambda}(ya)\,{\rm d}a\,{\rm d}y\\
 &=&\int_{\S\times A}\Delta(a)\phi(ya^{-1})\CJ_{\theta,\lambda}(a)\,\Omega_{\theta,\lambda}(y)\,{\rm d}a\,{\rm d}y
 \end{eqnarray*}
 where 
 $$
 \Delta(a){\rm d}(ya)\;:=\;{\rm d}y\;.
 $$
 Observing that $\Delta(a)\;=\;e^{2a}$, we find
 \begin{eqnarray*}
 {\bf \Omega}_{\theta,\lambda}(\phi)&=&\frac{4e^{i\pi{\left(\frac{\theta+2}{\theta}\right)}}}{\theta}\int_{\S\times A}\left(R^\star_{a^{-1}}\phi\right)(y)\,e^{2\lambda a}\,J_{2{\left(\frac{\theta+2}{\theta}\right)}-1}(\frac{2}{\theta}e^{-2a})\,\Omega_{\theta,\lambda}(y)\,{\rm d}a\,{\rm d}y
 \\
 &=&\int_{\S\times A}\left(R^\star_{a^{-1}}\phi\right)(y)\,\CJ_{\theta,\lambda+1}(a)\,\Omega_{\theta,\lambda}(y)\,{\rm d}a\,{\rm d}y\;.
 \end{eqnarray*}
\EPf
\begin{dfn}\label{WW} Let $C_b(\S)$ denote the space of right-uniformly continuous functions on $\S$.
Denoting by $C_b(\S)^\infty$ the space of $C^\infty$-vectors for the action $R^\vee$ of $A$ on $C_b(\S)$, we consider 
the following continuous linear operator on $C_b(\S)^\infty$:
$$
\frac{-1}{4}e^{\frac{2i}{\theta}\pi}\;\Xi^{R^\vee}_{\theta,\lambda +1}\;=:\;\fW_{\theta,\lambda}\;:=\;\frac{4e^{i\pi{\left(\frac{\theta+2}{\theta}\right)}}}{\theta}\int_{A}e^{2\lambda a}\,J_{2{\left(\frac{\theta+2}{\theta}\right)}-1}(\frac{2}{\theta}e^{-2a})\,R^\star_{a^{-1}}\,{\rm d}a\;.
$$
We call the above element the \emph{inverse retract operator}.
\end{dfn}
\begin{lem}\label{INVR}
We have

\begin{eqnarray*}
\left(\Xi^{R^\vee}_{\theta, \lambda +1}\right)^{-1}&=&-e^{\frac{2i}{\theta}\pi}\,\Xi^{R}_{\theta,1-\lambda}\;=\;
-\frac{4e^{i\pi{\left(\frac{\theta+4}{\theta}\right)}}}{\theta}\int_{A}e^{-2\lambda a}\,J_{2{\left(\frac{\theta+2}{\theta}\right)}-1}(\frac{2}{\theta}e^{-2a})\,R^\star_{a}\,{\rm d}a\;.
\end{eqnarray*}
\end{lem}
\Pf
We start with the well known inversion formula:
$$
\int\int J_{\mu}(\eta)\,J_{\mu}(\xi)\,\,
\phi(\frac{\eta}{\xi})\,{\rm d}\eta\,{\rm d}\xi\;=\;\phi(1)\;.
$$
In other words, for every $x$, the operator:
$$
A\phi(x)\;:=\;\int J_\mu(\xi)\,\phi(x\xi^{-1})\,{\rm d}\xi
$$
has inverse:
$$
A^{-1}\psi(x)\:=\;\int J_\mu(\eta)\,\psi(x\eta){\rm d}\eta\;.
$$
Now we consider the operators
$$
A_\lambda\phi(x)\;:=\;\int J_\mu(\xi)\,\xi^\lambda\phi(x/\xi){\rm d}\xi
$$
and
$$
B_\tau\phi(x)\;:=\;\int J_\mu(\eta)\eta^\tau\phi(x\eta){\rm d}\eta\;.
$$
We have then
$$
B_\tau A_\lambda(\phi)(x)\;=\;\int J_\mu(\eta)J_\mu(\xi)\eta^\tau\xi^\lambda\phi(x\eta/\xi){\rm d}\eta\,{\rm d}\xi\;.
$$
Setting $\Psi(\eta/\xi)\;:=\;(\eta\xi^{-1})^{-\lambda}\phi(x\eta/\xi)$ and $\tau=-\lambda$, the above inversion formula leads to
$$
B_{-\lambda} A_\lambda(\phi)(x)\;=\;\Psi(1)\;=\;\phi(x)\;.
$$
The lemma then follows when setting, for $x\in\S$:
$$
\xi\;:=\;\frac{2}{\theta}\,e^{-2a}\quad\mbox{\rm and}\quad\phi(\xi^{-1})\;:=:\;\varphi_x(\xi^{-1})\;:=\;\varphi(x.\xi^{-1})\;.
$$
\EPf
\section{Fundamental property of the retract}\label{FUND}

We now pass to extensions of the retract operators to spaces of type $\CB$. The following proposition involves the parameter $\lambda$. It is the only reason why this parameter appears in this paper.
\begin{prop}\label{LAMBDA}
Let $\theta$ be such that ${\left(\frac{\theta+2}{\theta}\right)}>0$.

(i) For all $s$ and $k\in\N$, there exists $\lambda_{\theta,k,s}$ such that  the inverse retract operator $\fW_{\theta,\lambda_{\theta,s,k}}$ extends 
from $\CD(\S)$ to $\CB^s_\infty(\S)$ as a
continuous injection $\CB^s_\infty(\S)$ into $\CB^s_k(\S)$:
$$
\fW_{\theta,\lambda_{\theta,k,s}}:\CB^s_\infty(\S)\to\CB^s_k(\S)\;.
$$
(ii) And similarly for the direct retract: for every $p\in\N$ with $2{\left(\frac{\theta+2}{\theta}\right)}-1/2\leq p\leq k$, 
the direct  retract operator associated with the same parameter $\lambda_{\theta,k,s}$ as in (i)
continuously injects $\CB^s_p(\S)$ into $\CB^s_p(\S)$:
$$
\fW^{-1}_{\theta,\lambda_{\theta,k,s}}:\CB^s_p(\S)\to\CB^s_p(\S)\;.
$$
\end{prop}
\Pf
Recall that considering $x\in\S$ and $a\;:=\;(a,0)\in A$, we have for every $k$:
\begin{equation}\label{ZETA}
\tilde{E}^k_x(R^\star_{a^{-1}}\varphi)\;=\;e^{2ka}\tilde{X}_{x.a^{-1}}(\varphi)\quad\mbox{\rm and}\quad\tilde{E}^k_x(R^\star_{a}\varphi)\;=\;e^{-2ka}\tilde{X}_{x.a}(\varphi)\;.
\end{equation}
 \noindent We will use an integration by parts technique in order to estimate the retracts. \noindent For convenience, we will express the group $A$ under its multiplicative form:
$$
A\;=\;\R^+_0\quad\mbox{\rm through}\quad\zeta\;:=\;e^{-2a}
$$
For every $\varphi\in\CD(\S)$, the inverse retract essentially turns into the form
$$
\fW\varphi(x)\;:=\;\int_0^\infty \zeta^{-\lambda}J_\mu(\zeta)\,\varphi(x.a^{-1})\,{\rm d}\zeta
$$
where we adopt the same notations as in the proof of Lemma \ref{INVR} {and where we redefined 
$\lambda\leftarrow\lambda+1$ (due to the change of variable $\zeta\;:=\;e^{-2a}$, ${\rm d}a$ becomes 
$\zeta^{-1}{\rm d}\zeta$)}.

\noindent Since
$$
\tilde{H}|_{C^\infty(A)}\;=\;\partial_a\;:=\;-2\zeta\partial_\zeta\;,
$$
we observe:
$$
e^{-\frac{\zeta^2}{4Z}}\;=\;\frac{1}{1+\zeta^2}\,(1+Z\,\tilde{H})e^{-\frac{\zeta^2}{4Z}}\;.
$$
The adjoint operator w.r.t. {$\zeta$ (and not to $a$)}:
$$
\bO\;:=\;\left[\frac{1}{1+\zeta^2}\,(1+Z\,\tilde{H})\right]^\ast\;=\;(1{+2Z}-Z\tilde{H})\,\circ\,\frac{1}{1+\zeta^2}
$$
applied to $\psi\in\CD(\S)$ involves  multiplications by $\frac{1}{1+\zeta^2}$ and  by $Z$ of derivatives of $\psi$. An induction
yields
$$
\bO^N\psi\;=\;\frac{1}{\left[1+\zeta^2\right]^N}\sum_{n,m=0}^Nb_{m,n}(\zeta)\,Z^n\,(\tilde{H})^m.\psi
$$
where the $b_{m,n}$ are smooth bounded functions with all derivatives bounded i.e. belonging to $\CB^0_\infty(A)$.

 \noindent We now use  the expression:
$$
\frac{1}{2\pi i}\int_{\epsilon+i\R}Z^{n-(1+\mu)}e^{Z-\frac{\zeta^2}{4Z}}\,{\rm d}Z\;=\;\left(\frac{2}{\zeta}\right)^{\mu-n} J_{\mu-n}(\zeta)\;.
$$An integration by parts then yields:
\begin{eqnarray*}
&&
\fW\varphi(x)\;=\;\sum_{m,n}^N\int_0^\infty\int_{\epsilon+i\R}
\frac{Z^{n-(1+\mu)}e^{Z-\frac{\zeta^2}{4Z}}b_{m,n}(\zeta)}{\left[1+\zeta^2\right]^N}\,(\tilde{H})_{\zeta}^m\left(\,\zeta^{\mu-\lambda} L_x^\star\varphi\right)\,{\rm d}Z\,{\rm d}\zeta\;.
\end{eqnarray*}
Observing that $\zeta^{\mu-\lambda}$ is an eigen-function of $\tilde{H}$, we are led to bound  expressions of the following type:
$$
\CI\;=\;\int_0^\infty\int_{\epsilon+i\R}
\frac{\zeta^{\mu-\lambda}\,Z^{n-(1+\mu)}e^{Z-\frac{\zeta^2}{4Z}}b(\zeta)}{\left[1+\zeta^2\right]^N}\,L_x^\star\left[(\tilde{H})^m\varphi\right]\,{\rm d}Z\,{\rm d}\zeta\;.
$$
The latter entails for every $X\in\CU_k(\s)$:
$$
\tilde{X}_x.\fW\varphi\;=\;\sum_{m,n}^N\int_0^\infty
\frac{\zeta^{n-\lambda} \,b^{(2)}_{m,n}(\zeta)}{\left[1+\zeta^2\right]^N}\,J_{\mu-n}(\zeta)\,\left(\Ad_{\zeta^{-1}}(X)\,H^m\right)^\sim_{xa^{-1}}.\varphi\,{\rm d}\zeta\quad(b^{(2)}_{m,n}\in\CB^1_\infty(A))\;.
$$
Therefore, we get positive constants $\{c_{m,n}\}$ and $X'\in\CU_k(\s)$ {(at worse c.f. (\ref{ZETA}))} such that, for every $s\in\R$:
$$
\sup\{\frac{1}{\fd^s}|\tilde{X}(\fW\varphi)|\}\;\leq\;\sum_{m,n=0}^Nc_{m,n}\sup_x\left\{
\frac{1}{\fd^s(x)}\int_{0}^\infty\frac{\zeta^{n-k-\lambda}}{[1+\zeta^2]^N}\,|J_{\mu-n}(\zeta)|\,\left|\widetilde{H^mX'}_{x.a^{-1}}(\varphi)\right|\,{\rm d}\zeta\right\}\;.
$$
Hence, because $\fd|_A(\zeta)\;=\;\left(\zeta+1/\zeta\right)$ and submultiplicativity of the
weight $\fd$, one has the following expression in terms of the input semi-norms associated to
$s'\in\R$:
$$
|\fW\varphi|_{k,s}\,\leq\,\sum_{m,n=0}^Nc_{m,n}\,|\varphi|_{k+m,s'}\,\sup_x\left\{
\frac{1}{\fd^{s-s'}(x)}\right\}\,\int_{0}^\infty\frac{\zeta^{n-k-\lambda}}{[1+\zeta^2]^N}\,|J_{\mu-n}(\zeta)|\,
\left(\zeta+1/\zeta\right)^{s'}{\rm d}\zeta\;.
$$
Now, 

$\bullet$ for $\zeta$ near zero the integrand behaves, up to multiplicative constant as
$$
\frac{\zeta^{n-k-\lambda}}{[1+\zeta^2]^N}\,|J_{\mu-n}(\zeta)|\,\left(\zeta+1/\zeta\right)^{s'}
\,\sim\,\zeta^{\mu-k-s'-\lambda}\;,$$
while 

$\bullet$ for $\zeta$ large, we have
$$
\frac{\zeta^{n-k-\lambda}}{[1+\zeta^2]^N}\,|J_{\mu-n}(\zeta)|\left(\zeta+1/\zeta\right)^{s'}\,\sim\,\zeta^{n-k-2N-1/2+s'-\lambda}\;.
$$
For given $k$ and $s$, finiteness of $|\fW\varphi|_{k,s}$ therefore amounts to satisfying the following inequalities:
$$
\left\{
\begin{array}{ccc}
s-s'&\geq&0\\
0&<&\mu-k-s'+1-\lambda\\
n-k-2N-1/2+s'-\lambda&<&-1\;.
\end{array}
\right.
$$
In other words ($n\leq N$):
$$
(1)\quad\left\{
\begin{array}{ccc}
s'\leq&s\\
k-1+s'+\lambda&<&\mu\\
-k-N+s'-\lambda&<&-1/2\;.
\end{array}
\right.
$$
The inequality $k-1+s'+\lambda<\mu$ imposes a strong condition on $\lambda$: in order 
to hold, one is constrained to choose $\lambda$ smaller than the limiting quantity $\lambda_{\theta,k,s'}$ given by
$$
\lambda_{\theta,k,s'}\;:=\;\mu-s'-k+1\;.
$$
Which for large $\mu$ is generally negative.

\noindent The function $\fW_{m,{\left(\frac{\theta+2}{\theta}\right)},\lambda_{\theta,k,s'}}\varphi$ is therefore $C^k$.

\noindent Now we analyse the situation with the order of derivation $\ell$ is smaller 
than $k$:
$$
\ell\leq k\;.
$$
The analyse therefore leads us to (${\left(\frac{\theta+2}{\theta}\right)}$ being large):
$$
(1')\quad
\begin{array}{ccc}
\ell-1+s'+\lambda_{\theta,k,s'}&\leq&\mu\;.
\end{array}
$$
Which amounts to 
$$
\ell-1+s'+\mu-s'-k+1\leq\mu\quad\mbox{\rm i.e.}\quad\ell-k\leq0
$$
which holds. Therefore the element $\fW_{m,{\left(\frac{\theta+2}{\theta}\right)},k,s',\lambda_{\theta,k,s'}}\varphi$ belongs to $\CB_k^s(\S)$\;.

\vspace{3mm}

\noindent Concerning the assertion relative to the direct retract and 
within obvious notations omitting indices, we have the expression
$$
\fW^{-1}\varphi(x)\;:=\;\int_0^\infty \zeta^{\lambda}J_\mu(\zeta)\,\varphi(x.a)\,{\rm d}\zeta
\quad (\varphi\in\CD(\S))\;.
$$
The same discussion as for the inverse retract leads to the following estimate:
$$
|\fW\varphi|_{k,s,p}\,\leq\,\sum_{m,n=0}^Nc_{m,n}\,|\varphi|_{k+m,s'}\,\sup_x\left\{
\frac{1}{\fd^{s-s'}(x)}\right\}\,\int_{0}^\infty\frac{\zeta^{n+k+\lambda}}{[1+\zeta^2]^{N'}}\,|J_{\mu-n}(\zeta)|\,
\left(\zeta+1/\zeta\right)^{s'}{\rm d}\zeta\;.
$$
Now, 

$\bullet$ for $\zeta$ near zero the integrand behaves, up to multiplicative constant as
$$
\frac{\zeta^{n+k+\lambda}}{[1+\zeta^2]^{N'}}\,|J_{\mu-n}(\zeta)|\,\left(\zeta+1/\zeta\right)^{s'}
\,\sim\,\zeta^{\mu+k-s'+\lambda}\;,$$
while 

$\bullet$ for $\zeta$ large, we have
$$
\frac{\zeta^{n+k+\lambda}}{[1+\zeta^2]^{N'}}\,|J_{\mu-n}(\zeta)|\left(\zeta+1/\zeta\right)^{s'}\,\sim\,\zeta^{n+k-2N-1/2+s'+\lambda}\;.
$$
For given $k$ and $s$, finiteness of $|\fW\varphi|_{k,s}$ therefore amounts to satisfying the following inequalities:
$$
\left\{
\begin{array}{ccc}
s-s'&\geq&0\\
0&<&\mu+k-s'+1+\lambda\\
k-{N'}-1/2+s'+\lambda&<&-1\;.
\end{array}
\right.
$$
In other words
$$
(2)\quad\left\{
\begin{array}{ccc}
s'&\leq&s\\
-k-1+s'-\lambda&<&\mu\\
k-{N'}+s'+\lambda&<&-1/2\;.
\end{array}
\right.
$$
Which, together with the system (1) gives
$$
\left\{
\begin{array}{ccc}
s'\leq&s\\
k-1+s'+\lambda&<&\mu\\
-k-N+s'-\lambda&<&-1/2\\
-k-1+s'-\lambda&<&\mu\\
k-{N'}+s'+\lambda&<&-1/2
\end{array}
\right.
$$
where $N$ and $N'$ are independant. They can therefore be chosen such that the
inequalities $-k-N+s'-\lambda<-1/2$ and $k-{N'}+s'+\lambda<-1/2$ are 
satisfied for any choice of $\lambda$.
These conditions in the limit case $\lambda=\lambda_{\mu,k,s'}-\epsilon$ amount to
$$
\left\{
\begin{array}{ccc}
s'&\leq&s\\
s'-1+\epsilon/2&<&\mu\\
\end{array}
\right.
$$
But now the function $\varphi$ is only $C^p$ with possibly $p\leq k$. Which leaves us with an output
$\fW^{-1}_{\theta, \lambda_{\theta,k,s'}}\varphi$ at best in $C^{p-N'}$, with $N'>\mu+3/2$ i.e.
$$
\mu+3/2<N'\leq p\;.
$$
Imposing
$$
k\geq p\geq \mu+3/2\;.
$$
Now, the set of inequalities (2), using the density of $\CD$-type symbols in $\CB$, leads to the fact that
$\fW^{-1}$ sends $\CB^s_\infty(\S)$ continuously into $\CB^{s}_\infty(\S)$, for every $s\in\R$.

\noindent While the second one (2) implies the same result for $\fW^{-1}$ since $N'$ is arbitrarily large.

\vspace{3mm}

\noindent A similar argument (but much simpler) as above applies in the case of the improper 
integral involved in the statement of Proposition \ref{LAMBDA0}.
\EPf
\noindent As a direct consequence of the above discussion we have
\begin{thm}\label{STABSCHWARTZ}
The retract operators stabilize the Schwartz space $\CS(\S)\;=\;\CS(\M)\;=\;\CS(\D)$ (c.f. Definition \ref{SCHWARTZSPACE}).
\end{thm}

\section{Explicit formula for the non-formal star-product and Drinfel'd twist}

\noindent We realize the hyperbolic plane as an adjoint orbit in $\g$ of a compact element $Z_0$. The dual picture would be more intrinsical but the computations are
simpler within the adjoint orbit setting. We consider the moment map
$$
J:\S\to\g:x\mapsto\Ad_xZ_0
$$
where $\S$ is realized as a sub-group of $G$.

\noindent At every $x\in\S$, we realize the tangent space as the sub-space of $\g$ Killing-orthogonal to the line, denoted hereafter $\k_x$, supported by $J(x)$:
$$
\p_x\;:=\;J_{\star x}(T_x\S)\;=\;J(x)^\perp\;.
$$
The (partial) {\bf \v{S}evera map} is then defined at $x$ as 
$$
\bsigma_x:\S\to\p_x\,\subset\,\g:y\mapsto\pi_x(J(y))
$$
where $\pi_x:\g\to\p_x$ denotes the Killing-orthogonal projection. We then observe
\begin{lem} The \v{S}evera map is $\S$-equivariant:
$$
\sigma_{gx}\;=\;\Ad_g\circ\sigma_x\circ L_{g^{-1}}\;.
$$
For $Z_0\;:=\;\frac{1}{2}(E-F)$, one has
$$\bsigma_e(y)\;=\;[Z_0,[\Ad_{y}Z_0,Z_0]]$$
and more generally:
$$
\sigma_x(y)\;=\;[J(x),[J(y),J(x)]]\;.
$$
\end{lem}
A simple computation then yields
\begin{prop}\label{PSISEVERA}
Let us normalize the Killing form $\beta$ on $\g$ in such a way that
$$
\beta(Z_0,Z_0)\;=-1
$$
and endow the vector space $\p:=\p_e$ with the symplectic structure
$$
\omega_0(X,Y)\;:=\;\beta(Z_0,[X,Y])\;.
$$
Then the vectors 
$$
e_q\;:=\;\frac{1}{2}(E+F)\quad \quad e_p\;:=\;\frac{1}{2}H
$$
constitute an orthnormal symplectic basis of $\p$ ($\omega_0(e_q,e_p)=1$) and through the symplectic identification
$$
\R^2\to\p:(q,p)\mapsto qe_q+pe_p
$$
the \v{S}evera map reads
$$
\sigma_e(\exp(aH)\exp(nE))\;=\;\Psi^{-1}Z(a,n)
$$
where $Z$ is the natural map
$$
Z\;:\;\S\to\;\Pi\;:\;an\mapsto\;an(1)\;.
$$
\end{prop}
\begin{dfn}
The {\bf canonical area} around point $x$ is defined as
$$
S(x,y,z)\;:=\;\beta(J(x),[\sigma_x(y),\sigma_x(z)])\;.
$$
We set 
$$
S(x,y)\;:=\;S(e,x,y)\;.
$$
\end{dfn}
\begin{rmk}
We note the relevance of the function
$$
\hat{S}:\L\times\S\to\g:(x,y,z)\mapsto[\sigma_x(y),\sigma_x(z)]\;.
$$
which, in the present two-dimensional and normalized setting, is related to the canonical area through the identity:
$$
\hat{S}(x,y,z)\;=\;S(x,y,z)\,J(x)\;.
$$
\end{rmk}

\noindent We now give an geometrical explicit formula (i.e. coordinate free) for the $\mbox{SL}_{2}(\R)$-invariant star product on the hyperbolic plane.
\begin{thm}\label{FINALTHM}
Let 
$$
\D\;:=\;\mbox{SL}_{2}(\R)/\mbox{SO}(2)
$$
be the hyperbolic plane.

\noindent Consider the following special function of the real variable:
$$
{\bf E}_{\lambda}(t)\;:=\;\int_{0}^{\infty}\,s^{2}
e^{-\,is\,t}
\;\,J_{\lambda}(s)\;{\rm d}s\;=\;-\ddtd\CL(J_{\lambda})(-it)\;.
$$

\noindent For a suitable choice of a constant $c_{\theta}$, the product
$$
\varphi\sharp_{\theta}\psi(x)\;:=\;c_{\theta}\int_{\D^{2}}A(x,y,z)\,{\bf E}_{2\frac{\theta+2}{\theta}-1}\left(S(x,y,z)\right)\,
\varphi(y)\,\psi(z)\;{\rm d}y\;{\rm d}z
$$
where 
\begin{enumerate}
\item ${\rm d}x$ denotes (a normalisation of) the $\mbox{SL}_{2}(\R)$-invariant measure on $\D$,
\item $S(x,y,z)$ denotes the \v{S}evera area on $\D$ (c.f. \ref{SEVERAAREA}),
\item and $A$ is defined in terms of the Jacobians of the \v{S}evera projections (c.f. \ref{SEVERAPROJ}) as
$$
A(x,y,z)\;:=\;\mbox{Jac}_{\sigma_{x}}(y)\,\mbox{Jac}_{\sigma_{x}}(z)
$$
\end{enumerate}
is such that 
\begin{enumerate}
\item[(i)] it extends from $C^{\infty}_{c}(\D)$ to the Schwartz space $\CS(\D)$ (c.f. Definition \ref{SCHWARTZSPACE}) as an $\mbox{SL}_{2}(\R)$-equivariant associative Fr\'echet algebra structure.
\item[(ii)] It admits an asymptotic expansion as an $\mbox{SL}_{2}(\R)$-invariant formal star-product on the hyperbolic plane equipped with (a normalization) of its $\mbox{SL}_{2}(\R)$-invariant symplectic Kahler two-from.
\item[(iii)] Identifying the hyperbolic plane $\D$ with the Iwasawa factor $\S:=AN$ of $\mbox{SL}_{2}(\R)=ANK$, the formula
$$
\CF_{\theta}\;:=\;c_{\theta}\,\int_{\S^{2}}\mbox{Jac}_{\sigma_{e}}(y)\,\mbox{Jac}_{\sigma_{e}}(z)\,{\bf E}_{2\frac{\theta+2}{\theta}-1}\left(S(y,z)\right)\,
R_{y}^{\star}\otimes\,R_{z}^{\star}\;{\rm d}y\;{\rm d}z
$$
defines a  Drinfel'd twist based on the Iwasawa factor $\S:=AN$.
\end{enumerate}
\end{thm}
\Pf
\noindent We first notice that
\begin{eqnarray*}
&&
\int_{\p^{2}}\CG(\varphi_{1})(x^{\prime})\,\CG(\varphi_{2})(x^{\prime\prime})\;\int_{0}^{\infty}\,J_{\lambda}(s)\;e^{i/s\,\omega^{\p}(x^{\prime}\,,\,x^{\prime\prime})}\;{\rm d}x^{\prime}\;{\rm d}x^{\prime\prime}\\
&=&
\int_{0}^{\infty}\,J_{\lambda}(s)\int_{\p^{4}}\varphi_{1}(x_{1})\,\varphi_{2}(x_{2})\;e^{i\,(\omega^{\p}(x_{1},x^{\prime})\,+\,\omega^{\p}(x_{2},x^{\prime\prime}))}\;e^{i/s\,\omega^{\p}(x^{\prime}\,,\,x^{\prime\prime})}\;{\rm d}s\;{\rm d}x^{\prime}\;{\rm d}x^{\prime\prime}\;{\rm d}x_{1}\;{\rm d}x_{2}\\
&=&
\int_{\p^{4}}\int_{0}^{\infty}
e^{i\,(\omega^{\p}(x_{1},x^{\prime})\,+\,\omega^{\p}(x_{2},x^{\prime\prime})\;+\;1/s\,\omega^{\p}(x^{\prime}\,,\,x^{\prime\prime}))}
\;\,J_{\lambda}(s)\varphi_{1}(x_{1})\,\varphi_{2}(x_{2})\;{\rm d}s\;{\rm d}x^{\prime}\;{\rm d}x^{\prime\prime}\;{\rm d}x_{1}\;{\rm d}x_{2}\\
&=&
\int_{\p^{4}}\int_{0}^{\infty}
e^{i\,(\omega^{\p}(x_{1}\,-\,1/s\,x^{\prime\prime},x^{\prime})\,+\,\omega^{\p}(x_{2},x^{\prime\prime}))}
\;\,J_{\lambda}(s)\varphi_{1}(x_{1})\,\varphi_{2}(x_{2})\;{\rm d}s\;{\rm d}x^{\prime}\;{\rm d}x^{\prime\prime}\;{\rm d}x_{1}\;{\rm d}x_{2}\\
&=&
\int_{\p^{2}}\int_{0}^{\infty}
e^{i\,\omega^{\p}(x_{2},x^{\prime\prime})}
\;\,J_{\lambda}(s)\varphi_{1}(1/s\,x^{\prime\prime})\,\varphi_{2}(x_{2})\;{\rm d}s\;{\rm d}x_{2}\;{\rm d}x^{\prime\prime}\\
&=&
\int_{\p^{2}}\int_{0}^{\infty}\,s^{2}
e^{-\,is\,\omega^{\p}(x_{1},x_{2})}
\;\,J_{\lambda}(s)\;\varphi_{1}(x_{1})\,\varphi_{2}(x_{2})\;{\rm d}s\;{\rm d}x_{1}\;{\rm d}x_{2}\;.
\end{eqnarray*}
\noindent \noindent The above Proposition \ref{PSISEVERA} shows that the symplectic structure $\omega^{\p}$ implied in the symplectic Fourier
transform $\CG$ is the \v{S}evera area element at $e$. Therefore, setting
$$
x_{j}\;=\;\Psi^{-1}(\zeta_{j})\;,
$$
we have (with a slight abuse of notation : identifying $\S$ with $\Pi$ through the map $Z$):
$$
\omega^{\p}(x_{1},x_{2})\;=\;S(o\,,\,\zeta_{1}\,,\,\zeta_{j})\;.
$$
Hence
$$
\varphi_{1}\bullet\;\varphi_{2}(0)\,=\,\int_{\Pi^{2}}\hspace{-1,5mm}\left(\int_{0}^{\infty}\,s^{2}
e^{-\,is\,S(\zeta_{1},\zeta_{2})}
J_{\lambda}(s)\;\;{\rm d}s\right)\varphi_{1}(\Psi^{-1}(\zeta_{1}))\,\varphi_{2}(\Psi^{-1}(\zeta_{2}))\,\mbox{Jac}_{\Psi^{-1}}(\zeta_{1})\,\mbox{Jac}_{\Psi^{-1}}(\zeta_{2})\,{\rm d}\zeta_{1}{\rm d}\zeta_{2}
$$
where ${\rm d}\zeta$ denotes the hyperbolic measure on the hyperbolic plane $\Pi$ (realized as a Poincar\'e half plane).

\noindent Therefore, for all smooth and compactly supported fonctions $f_{1}$ and $f_{2}$ on the hyperbolic plane $\Pi$, we have
$$
(\Psi^{-1\,\star})(f_{1})\bullet(\Psi^{-1\,\star})(f_{2})(0)=\int_{\Pi^{2}}\hspace{-1,5mm}\left(\int_{0}^{\infty}\,s^{2}
e^{-\,is\,S(\zeta_{1},\zeta_{2})}
J_{\lambda}(s){\rm d}s\right)\mbox{Jac}_{\Psi^{-1}}(\zeta_{1})\mbox{Jac}_{\Psi^{-1}}(\zeta_{2})f_{1}(\zeta_{1})f_{2}(\zeta_{2}){\rm d}\zeta_{1}{\rm d}\zeta_{2}
$$
where $\mbox{Jac}_{\Psi^{-1}}(\zeta)$ denotes the Jacobian of the map $\Psi^{-1}$ relating the symplectic (Lebesgue) measure on $\p$ with the hyperbolic measure on $\Pi$.

\noindent Now we observe that for every $\zeta\;=\;x\;+\;i\xi$ in the half plane $\Pi$, the change of coordinates
$$
\left\{\begin{array}{ccc}x&=&e^{2a}\\ \xi&=&-\,e^{2a}n\end{array}\right.
$$
identifies the product $f\bullet g$ with the product $\star$ defined by
$$
\left(Z^{\star}(f)\star Z^{\star}(g)\right)(Z(x_{0})):=\int_{\S\times\S}\cosh(2(a_{0}\,-\,a_{2})\cosh(2(a_{1}\,-\,a_{0})\,e^{i\,S(x_{0},x_{2},x_{2})}f(Z(x_{1}))\,g(Z(x_{2})){\rm d}x_{1}{\rm d}x_{2}.
$$
Within this setting, we also observe that for every $x\;+\;i\xi\;=\;\in\Pi$ and $t:=e^{2b}\in A$, one has
$$
tx\;+\;i\xi\;=\;Z(anb)\;.
$$
Therefore Unterberger's map $G_{\lambda}$ is a solution of the retract. Precisely we have
$$
G_{\lambda}g(Z(an))\;=\;2\pi\,\int_{A}g(Z(anb))\;J_{\lambda}(4\pi\,e^{2b})\;e^{2b}\;{\rm d}b\;,
$$
or equivalently:
$$
G_{\lambda}\circ Z^{\star}\;=\;2\pi\,\int_{A}\;J_{\lambda}(4\pi\,e^{2a})\;e^{2a}R_{a}^{\star}\,{\rm d}a\;.
$$
That is
$$
G_{2\frac{\theta+2}{\theta}-1}\;=\;\frac{1}{4}\,e^{-\,i\,\pi\frac{\theta+2}{\theta}}\,\fW_{-1\,,\,\theta}\,\circ\,R_{-\frac{1}{2}\log(2\pi\theta)}^{\star}\;.
$$
Note that since the star-product $\star$ is $AN$-invariant, the $A$-translation $R_{-\frac{1}{2}\log(2\pi\theta)}^{\star}$ won't be
seen when intertwining $\bullet$ with $G_{\lambda}^{-1}(\bullet)$.

\noindent As already observed by A. and J. Unterberger, the product $G_{\lambda}^{-1}(\bullet)$ is therefore $\mbox{SL}_{2}(\R)$-equivariant.

\noindent Again identifying $\S$ with $\Pi$ through the map $Z$, we have that
$$
\sigma_{o}(\eta)\;=\;\Psi^{-1}(\eta)\quad(\eta\in\Pi)\;.
$$
Hence, by the ${\mbox{SL}}_{2}(\R)$-equivariance of the  \v{S}evera map, we get, for every $g$ in ${\mbox{SL}}_{2}(\R)$:
\begin{eqnarray*}
&&
f_{1}\sharp f_{2}(g.o):=\int_{\Pi^{2}}\hspace{-1,5mm}\left(\int_{0}^{\infty}\,s^{2}
e^{-\,is\,S(\zeta_{1},\zeta_{2})}
J_{\lambda}(s){\rm d}s\right)\mbox{Jac}_{\Psi^{-1}}(\zeta_{1})\mbox{Jac}_{\Psi^{-1}}(\zeta_{2})f_{1}(g.\zeta_{1})f_{2}(g\zeta_{2}){\rm d}\zeta_{1}{\rm d}\zeta_{2}\\
&=&
\int_{\Pi^{2}}\hspace{-1,5mm}\left(\int_{0}^{\infty}\,s^{2}
e^{-\,is\,S(g^{-1}\zeta_{1},g^{-1}\zeta_{2})}
J_{\lambda}(s){\rm d}s\right)\mbox{Jac}_{\Psi^{-1}}(g^{-1}\zeta_{1})\mbox{Jac}_{\Psi^{-1}}(g^{-1}\zeta_{2})f_{1}(\zeta_{1})f_{2}(\zeta_{2}){\rm d}\zeta_{1}{\rm d}\zeta_{2}\\
&=&
\int_{\Pi^{2}}\hspace{-1,5mm}\left(\int_{0}^{\infty}\,s^{2}
e^{-\,is\,S(g.o\,,\,\zeta_{1},\zeta_{2})}
J_{\lambda}(s){\rm d}s\right)\mbox{Jac}_{\sigma_{g.o}}(\zeta_{1})\mbox{Jac}_{\sigma_{g.o}}(\zeta_{2})f_{1}(\zeta_{1})f_{2}(\zeta_{2}){\rm d}\zeta_{1}{\rm d}\zeta_{2}
\end{eqnarray*}
i.e., for every $\zeta$ in $\Pi$:
$$
f_{1}\sharp f_{2}(\zeta)=\int_{\Pi^{2}}\hspace{-1,5mm}\left(\int_{0}^{\infty}\,s^{2}
e^{-\,is\,S(\zeta\,,\,\zeta_{1},\zeta_{2})}
J_{\lambda}(s){\rm d}s\right)\mbox{Jac}_{\sigma_{\zeta}}(\zeta_{1})\mbox{Jac}_{\sigma_{\zeta}}(\zeta_{2})\;f_{1}(\zeta_{1})\;f_{2}(\zeta_{2})\;{\rm d}\zeta_{1}{\rm d}\zeta_{2}\;.
$$
We conclude by Theorem \ref{STABSCHWARTZ}.
\EPf
\begin{rmks}
\begin{enumerate}
\item[(i)] In Theorem \ref{FINALTHM} item (iii), despite the fact that the formula is non-formal, it is not clear that the  Drinfel'd twist defines a strict universal deformation formula (in the sense of Rieffel) for actions of $AN$ on Fr\'echet $AN$-module algebras. Indeed, an analysis analogue to what has been done in Chapter \ref{OSCINT} still needs to be done. So, in the present state, the Drinfel'd twist $\CF_{\theta}$ should be considered as formal (through its asymptotic expansion).
\end{enumerate}
\end{rmks}

\section*{Appendix: Andr\'e and Julianne Unterberger's symbolic calculus}

In \cite{UU88}, A. and J. Unterberger consider the open half plane
$$
\Pi\;:=\;\{x+i\xi\;|\;x>0\;;\;\xi\in\R\}\;\subset\;\C
$$
equipped with the action of the group $\mbox{SL}_{2}(\R)$ given by
$$
\left(\begin{array}{ccc}a&&b\\c&&d\end{array}\right).\zeta\;:=\;\frac{a\zeta\;-\;ib}{ic\zeta\;+\;d}\;.
$$
This is an equivalent realization of the hyperbolic plane equipped with the natural isometric action of $\mbox{SL}_{2}(\R)$.   

\noindent The so-called ``composition formula for Fuchs symbols'' $\bff$ and $\bg\in C^\infty_c(\Pi)$ is expressed as
\begin{eqnarray*}
\bff\bullet\bg(x+i\xi)&:=&\int_{\Pi\times\Pi}\bff(y'+i\eta')\,\bg(y"+i\eta")\\
&&
\exp-2i\pi\left\{
\left(\frac{y'}{y"}-\frac{y"}{y'} \right)\frac{\xi}{x}+ \left(\frac{y"}{x}-\frac{x}{y"} \right)\frac{\eta'}{y'}+ \left(\frac{x}{y'}-\frac{y'}{x} \right)\frac{\eta"}{y"}
\right\}\\
&&
\left(\frac{x}{y'}+\frac{y'}{x} \right)\,\left(\frac{x}{y"}+\frac{y"}{x} \right)\,
{y'}^{-2}\,{y"}^{-2}\,{\rm d}y'\,{\rm d}\eta'\,{\rm d}y"\,{\rm d}\eta"\;.
\end{eqnarray*}
Setting, for every $\lambda>0$:
$$
G_\lambda\bg(y+i\eta)\;:=\;4\pi\,\int_0^\infty\bg(ty+i\eta)\,J_\lambda(4\pi\,t)\,{\rm d}t\;,
$$
they prove
\begin{thm}
\begin{eqnarray}\label{UUFORMULA}
G_\lambda^{-1}\left(G_\lambda(\bff)\bullet G_\lambda(\bg)\right)(1)&=&
\int_{\R^2\times\R^2}\CG(\bff\circ\Psi)(q',p')\,\CG(\bg\circ\Psi)(q",p")\\
&&
4\pi\,\int_0^\infty\,J_\lambda(4\pi\,s)\,\exp\frac{i\pi}{s}(-q'p"+q"p')\,{\rm d}s\\
&&
{\rm d}q'\,{\rm d}p'\,{\rm d}q"\,{\rm d}p"
\end{eqnarray}
where
$$
\Psi:\R^2\to\Pi:(q,p)\mapsto\frac{\sqrt{1+q^2+p^2}+q}{1-ip}
$$
and 
$$
\CG\varphi(\alpha,\beta)\;:=\;\int_{\R^2}\varphi(q,p)\,e^{2i\pi(-\alpha p+\beta q)}\,{\rm d}q\,{\rm d}p\quad(\varphi\in C^\infty_c(\R^2))\;.
$$
\end{thm}

\end{document}